# Florentin Smarandache
*editor*

# Proceedings of the First International Conference on Neutrosophy, Neutrosophic Logic, Neutrosophic Set, Neutrosophic Probability and Statistics

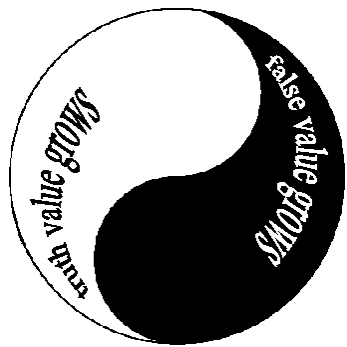
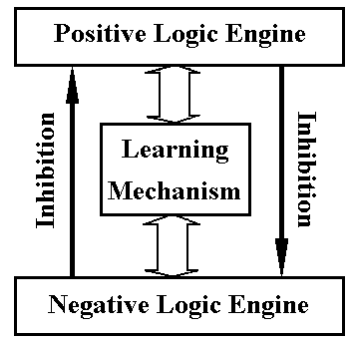

*University of New Mexico - Gallup*
*1-3 December 2001*

# FLORENTIN SMARANDACHE
*editor*

# PROCEEDINGS OF THE FIRST INTERNATIONAL CONFERENCE ON NEUTROSOPHY, NEUTROSOPHIC LOGIC, NEUTROSOPHIC SET, NEUTROSOPHIC PROBABILITY AND STATISTICS

*2002*

# *First International Conference on Neutrosophy, Neutrosophic Logic, Set, Probability and Statistics*

*December 1-3, 2001*
*University of New Mexico*

*Contributed papers were sent, by November 30, 2001, to the organizer:*
*Florentin Smarandache, University of New Mexico,*
*200 College Road, Gallup, NM 87301, USA.*
*Tel.: (505) 863-7647, Fax: (505) 863-7532 (Attn: Neutrosophic Conference). E-mail: smarand@unm.edu,*
*http://www.gallup.unm.edu/~smarandache/FirstNeutConf.htm.*
*A selection of them is being published in these Proceedings of the Conference.*

*Abstracts of papers were submitted to the following web site, provided by The York University, from Toronto, Canada, at: http://at.yorku.ca/cgi-bin/amca/submit/cagu-01, and they can be viewed at http://at.yorku.ca/cgi-bin/amca/cagu-01.*

*Invited Speakers: J. Dezert (France), Charles Le (US), I. Stojmenovic (Canada).*

*For more information on neutrosophics see the below links:*
- *0. Introduction:*
http://www.gallup.unm.edu/~smarandache/Introduction.PDF
- *1. Neutrosophy - a new branch of philosophy:*
http://www.gallup.unm.edu/~smarandache/Neutrosophy.PDF
- *2. Neutrosophic Logic - a unifying field in logics:*
http://www.gallup.unm.edu/~smarandache/NeutrosophicLogic.PDF
- *3. Neutrosophic Set - a unifying field in sets:*
http://www.gallup.unm.edu/~smarandache/NeutrosophicSet.PDF
- *4. Neutrosophic Probability - a generalization of classical and imprecise probabilities - and Neutrosophic Statistics:*
http://www.gallup.unm.edu/~smarandache/NeutrosophicProbStat.PDF

# PREFACE
# (An Introduction to Neutrosophy, Neutrosophic Logic, Neutrosophic Set, and Neutrosophic Probability and Statistics)

Florentin Smarandache
Department of Mathematics
University of New Mexico
Gallup, NM 87301, USA

## 0.1 Introduction to Non-Standard Analysis.

In 1960s Abraham Robinson has developed the non-standard analysis, a formalization of analysis and a branch of mathematical logic, that rigorously defines the infinitesimals. Informally, an infinitesimal is an infinitely small number. Formally, x is said to be infinitesimal if and only if for all positive integers n one has $|x| < 1/n$. Let $\varepsilon > 0$ be a such infinitesimal number. The hyper-real number set is an extension of the real number set, which includes classes of infinite numbers and classes of infinitesimal numbers. Let's consider the non-standard finite numbers $1^+ = 1+\varepsilon$, where "1" is its standard part and "$\varepsilon$" its non-standard part, and $^-0 = 0-\varepsilon$, where "0" is its standard part and "$\varepsilon$" its non-standard part.

Then, we call $]^-0, 1^+[$ a non-standard unit interval. Obviously, 0 and 1, and analogously non-standard numbers infinitely small but less than 0 or infinitely small but greater than 1, belong to the non-standard unit interval. Actually, by "$^-a$" one signifies a monad, i.e. a set of hyper-real numbers in non-standard analysis:

$\mu(^-a) = \{a-x : x \in \mathbb{R}^*, x \text{ is infinitesimal}\}$,

and similarly "$b^+$" is a monad:

$\mu(b^+) = \{b+x : x \in \mathbb{R}^*, x \text{ is infinitesimal}\}$.

Generally, the left and right borders of a non-standard interval $]^-0, 1^+[$ are vague, imprecise, themselves being non-standard (sub)sets $\mu(^-a)$ and $\mu(b^+)$ as defined above.

Combining the two before mentioned definitions one gets, what we would call, a binad of "$^-c^+$":

$\mu(^-c^+) = \{c-x : x \in \mathbb{R}^*, x \text{ is infinitesimal}\} \cup \{c+x : x \in \mathbb{R}^*, x \text{ is infinitesimal}\}$, which is a collection of open punctured neighborhoods (balls) of c.

Of course, $^-a < a$ and $b^+ > b$. No order between $^-c^+$ and c.

Addition of non-standard finite numbers with themselves or with real numbers:

$^-a + b = ^-(a + b)$
$a + b^+ = (a + b)^+$
$^-a + b^+ = ^-(a + b)^+$
$^-a + ^-b = ^-(a + b)$ (the left monads absorb themselves)
$a^+ + b^+ = (a + b)^+$ (analogously, the right monads absorb themselves)

Similarly for subtraction, multiplication, division, roots, and powers of non-standard finite numbers with themselves or with real numbers.

By extension let $\inf ]^-0, 1^+[ = ^-a$ and $\sup ]^-0, 1^+[ = b^+$.

## 0.2 Definition of Neutrosophic Components.

Let T, I, F be standard or non-standard real subsets of $]^-0, 1^+[$,

with    sup T = t_sup, inf T = t_inf,
        sup I = i_sup, inf I = i_inf,
        sup F = f_sup, inf F = f_inf,
and     n_sup = t_sup+i_sup+f_sup,
        n_inf = t_inf+i_inf+f_inf.

The sets T, I, F are not necessarily intervals, but may be any real sub-unitary subsets: discrete or continuous; single-element, finite, or (countably or uncountably) infinite; union or intersection of various subsets; etc.
They may also overlap. The real subsets could represent the relative errors in determining t, i, f (in the case when the subsets T, I, F are reduced to points).
Statically T, I, F are subsets.

*But dynamically, looking therefore from another perspective, the components T, I, F are at each instance dependant on many parameters, and therefore they can be considered set-valued vector functions or even operators.* The parameters can be: time, space, etc. (some of them are hidden/unknown parameters): T(t, s, …), I(t, s, …), F(t, s, …), where t=time, s=space, etc., that's why the neutrosophic logic can be used in quantum physics.
The Dynamic Neutrosophic Calculus can be used in psychology.
Neutrosophics try to reflect the dynamics of things and ideas.
See an example:

The proposition "Tomorrow it will be raining" does not mean a fixed-valued components structure; this proposition may be say 40% true, 50% indeterminate, and 45% false at time $t_1$; but at time $t_2$ may change at 50% true, 49% indeterminate, and 30% false (according with new evidences, sources, etc.); and tomorrow at say time $t_{145}$ the same proposition may be 100%, 0% indeterminate, and 0% false (if tomorrow it will indeed rain). This is the dynamics: the truth value changes from a time to another time.
In other examples:
In other examples: the truth value of a proposition may change from a place to another place, for example: the proposition "It is raining" is 0% true, 0% indeterminate, and 100% false in Albuquerque (New Mexico), but moving to Las Cruces (New Mexico) the truth value changes and it may be (1, 0, 0).
Also, the truth value depends/changes with respect to the observer (subjectivity is another parameter of the functions/operators T, I, F). For example: "John is smart" can be (.35, .67, .60) according to his boss, but (.80, .25, .10) according to himself, or (.50, .20, .30) according to his secretary, etc.

In the this book T, I, F, called *neutrosophic components*, will represent the truth value, indeterminacy value, and falsehood value respectively referring to neutrosophy, neutrosophic logic, neutrosophic set, neutrosophic probability, neutrosophic statistics.

This representation is closer to the human mind reasoning. It characterizes/catches the imprecision of knowledge or linguistic inexactitude received by various observers (that's why T, I, F are subsets - not necessarily single-elements), uncertainty due to incomplete knowledge or acquisition errors or stochasticity (that's why the subset I exists), and vagueness due to lack of clear contours or limits (that's why T, I, F are subsets and I exists; in particular for the appurtenance to the neutrosophic sets).
One has to specify the superior (x_sup) and inferior (x_inf) limits of the subsets because in many problems arises the necessity to compute them.

### 0.3. Operations with Sets.

Let $S_1$ and $S_2$ be two (unidimensional) real standard or non-standard subsets, then one defines:

*Addition of Sets:*
$S_1 \oplus S_2 = \{x \mid x=s_1+s_2, \text{ where } s_1 \in S_1 \text{ and } s_2 \in S_2\}$,
with inf $S_1 \oplus S_2$ = inf $S_1$ + inf $S_2$, sup $S_1 \oplus S_2$ = sup $S_1$ + sup $S_2$;
and, as some particular cases, we have
$\{a\} \oplus S_2 = \{x \mid x=a+s_2, \text{ where } s_2 \in S_2\}$
with inf $\{a\} \oplus S_2$ = a + inf $S_2$, sup $\{a\} \oplus S_2$ = a + sup $S_2$.

*Subtraction of Sets:*
$S_1 \ominus S_2 = \{x \mid x=s_1-s_2, \text{ where } s_1 \in S_1 \text{ and } s_2 \in S_2\}$.
For real positive subsets (most of the cases will fall in this range) one gets
inf $S_1 \ominus S_2$ = inf $S_1$ - sup $S_2$, sup $S_1 \ominus S_2$ = sup $S_1$ - inf $S_2$;
and, as some particular cases, we have
$\{a\} \ominus S_2 = \{x \mid x=a-s_2, \text{ where } s_2 \in S_2\}$,
with inf $\{a\} \ominus S_2$ = a - sup $S_2$, sup $\{a\} \ominus S_2$ = a - inf $S_2$;
also $\{1^+\} \ominus S_2 = \{x \mid x=1^+-s_2, \text{ where } s_2 \in S_2\}$,
with inf $\{1^+\} \ominus S_2 = 1^+$ - sup $S_2$, sup $\{1^+\} \ominus S_2$ = 100 - inf $S_2$.

*Multiplication of Sets:*
$S_1 \odot S_2 = \{x \mid x=s_1 \cdot s_2, \text{ where } s_1 \in S_1 \text{ and } s_2 \in S_2\}$.
For real positive subsets (most of the cases will fall in this range) one gets
inf $S_1 \odot S_2$ = inf $S_1 \cdot$ inf $S_2$, sup $S_1 \odot S_2$ = sup $S_1 \cdot$ sup $S_2$;
and, as some particular cases, we have
$\{a\} \odot S_2 = \{x \mid x=a \cdot s_2, \text{ where } s_2 \in S_2\}$,
with inf $\{a\} \odot S_2$ = a * inf $S_2$, sup $\{a\} \odot S_2$ = a $\cdot$ sup $S_2$;
also $\{1^+\} \odot S_2 = \{x \mid x=1 \cdot s_2, \text{ where } s_2 \in S_2\}$,
with inf $\{1^+\} \odot S_2 = 1^+ \cdot$ inf $S_2$, sup $\{1^+\} \odot S_2 = 1^+ \cdot$ sup $S_2$.

*Division of a Set by a Number:*
Let $k \in \mathbb{R}^*$, then $S_1 \oslash k = \{x \mid x=s_1/k, \text{ where } s_1 \in S_1\}$.

1. **NEUTROSOPHY:**

   **1.1 Definition:**
   Neutrosophy is a new branch of philosophy that studies the origin, nature, and scope of neutralities, as well as their interactions with different ideational spectra.

It is the base of *neutrosophic logic*, a multiple value logic that generalizes the fuzzy logic and deals with paradoxes, contradictions, antitheses, antinomies.

**1.2 Characteristics** of this mode of thinking:
- proposes new philosophical theses, principles, laws, methods, formulas, movements;
- reveals that world is full of indeterminacy;
- interprets the uninterpretable;
- regards, from many different angles, old concepts, systems:
showing that an idea, which is true in a given referential system, may be false in another one, and vice versa;
- attempts to make peace in the war of ideas,
and to make war in the peaceful ideas;
- measures the stability of unstable systems,
and instability of stable systems.
Let's note by <A> an idea, or proposition, theory, event, concept, entity, by <Non-A> what is not <A>, and by <Anti-A> the opposite of <A>. Also, <Neut-A> means what is neither <A> nor <Anti-A>, i.e. neutrality in between the two extremes. And <A'> a version of <A>.
   <Non-A> is different from <Anti-A>.

**1.3 Main Principle:**
Between an idea <A> and its opposite <Anti-A>, there is a continuum-power spectrum of neutralities <Neut-A>.

**1.4 Fundamental Thesis of Neutrosophy:**
Any idea <A> is T% true, I% indeterminate, and F% false, where $T, I, F \subset ]^-0, 1^+[$.

**1.5 Main Laws of Neutrosophy:**
Let <α> be an attribute, and $(T, I, F) \subset ]^-0, 1^+[^3$. Then:
- There is a proposition <P> and a referential system {R}, such that <P> is T% <α>, I% indeterminate or <Neut-α>, and F% <Anti-α>.
- For any proposition <P>, there is a referential system {R}, such that <P> is T% <α>, I% indeterminate or <Neut-α>, and F% <Anti-α>.
- <α> is at some degree <Anti-α>, while <Anti-α> is at some degree <α>.

## 2. NEUTROSOPHIC LOGIC:

**2.1 Introduction.**
As an alternative to the existing logics we propose the Neutrosophic Logic to represent a mathematical model of uncertainty, vagueness, ambiguity, imprecision, undefined, unknown, incompleteness, inconsistency, redundancy, contradiction. It is a non-classical logic.
Eksioglu (1999) explains some of them:
 "Imprecision of the human systems is due to the imperfection of knowledge that humain
 receives (observation) from the external world. Imperfection leads to a doubt about the
 value of a variable, a decision to be taken or a conclusion to be drawn for the actual
 system. The sources of uncertainty can be stochasticity (the case of intrinsic imperfection
 where a typical and single value does not exist), incomplete knowledge (ignorance of the
 totality, limited view on a system because of its complexity) or the acquisition errors
 (intrinsically imperfect observations, the quantitative errors in measures)."

"Probability (called sometimes the objective probability) process uncertainty of random type (stochastic) introduced by the chance. Uncertainty of the chance is clarified by the time or by events' occurrence. The probability is thus connected to the frequency of the events' occurrence."
"The vagueness which constitutes another form of uncertainty is the character of those which contours or limits lacking precision, clearness. […]
For certain objects, the fact to be in or out of a category is difficult to mention. Rather, it is possible to express a partial or gradual membership."
Indeterminacy means degrees of uncertainty, vagueness, imprecision, undefined, unknown, inconsistency, redundancy.

A question would be to try, if possible, to get an axiomatic system for the neutrosophic logic. Intuition is the base for any formalization, because the postulates and axioms derive from intuition.

### 2.2 Definition:
A logic in which each proposition is estimated to have the percentage of truth in a subset T, the percentage of indeterminacy in a subset I, and the percentage of falsity in a subset F, where T, I, F are defined above, is called *Neutrosophic Logic*.

We use a subset of truth (or indeterminacy, or falsity), instead of a number only, because in many cases we are not able to exactly determine the percentages of truth and of falsity but to approximate them: for example a proposition is between 30-40% true and between 60-70% false, even worst: between 30-40% or 45-50% true (according to various analyzers), and 60% or between 66-70% false.
The subsets are not necessary intervals, but any sets (discrete, continuous, open or closed or half-open/half-closed interval, intersections or unions of the previous sets, etc.) in accordance with the given proposition.
A subset may have one element only in special cases of this logic.

Constants: (T, I, F) truth-values, where T, I, F are standard or non-standard subsets of the non-standard interval $]^-0, 1^+[$, where $n_{inf} = \inf T + \inf I + \inf F \geq {}^-0$, and $n_{sup} = \sup T + \sup I + \sup F \leq 3^+$.
Atomic formulas: a, b, c, … .
Arbitrary formulas: A, B, C, … .

The neutrosophic logic is a formal frame trying to measure the truth, indeterminacy, and falsehood.
My hypothesis is that **no theory is exempted from paradoxes**, because of the language imprecision, metaphoric expression, various levels or meta-levels of understanding/interpretation which might overlap.

### 2.3 Definition of Neutrosophic Logical Connectives:
The connectives (rules of inference, or operators), in any non-bivalent logic, can be defined in various ways, giving rise to lots of distinct logics. For example, in three-valued logic, where three possible values are possible: true, false, or undecided, there are 3072 such logics!

(Weisstein, 1998) A single change in one of any connective's truth table is enough to form a (completely) different logic.

> The rules are hypothetical or factual. How to choose them? The philosopher Van Fraassen (1980) [see Shafer, 1986] commented that such rules may always be controvertible "for it always involves the choice of one out of many possible but nonactual worlds". There are general rules of combination, and ad hoc rules.

For an applied logic to artificial intelligence, a better approach, the best way would be to define the connectives recursively (Dubois, Prade), changing/adjusting the definitions after each step in order to improve the next result. This might be comparable to approximating the limit of a convergent sequence, calculating more and more terms, or by calculating the limit of a function successively substituting the argument with values closer and closer to the critical point. The recurrence allows evolution and self-improvement.

Or to use *greedy algorithms*, which are combinatorial algorithms that attempt at each iteration as much improvement as possible unlike myopic algorithms that look at each iteration only at very local information as with steepest descent method.

As in non-monotonic logic, we make assumptions, but we often err and must jump back, revise our assumptions, and start again. We may add rules that don't preserve monotonicity.

In bio-mathematics Heitkoetter and Beasley (1993-1999) present the *evolutionary algorithms* which are used "to describe computer-based problem solving systems which employ computational models of some of the known mechanisms of evolution as key elements in their design and implementation". They simulate, via processes of selection, mutation, and reproduction, the evolution of individual structures. The major evolutionary algorithms studied are: genetic algorithm (a model of machine learning based on genetic operators), evolutionary programming (a stochastic optimization strategy based on linkage between parents and their offspring; conceived by L. J. Fogel in 1960s), evolution strategy, classifier system, genetic programming.

Pei Wang devised a Non-Axiomatic Reasoning System as an intelligent reasoning system, where intelligence means working and adopting with insufficient knowledge and resources.

The inference mechanism (endowed with rules of transformation or rules of production) in neutrosophy should be non-monotonic and should comprise ensembles of recursive rules, with preferential rules and secondary ones (priority order), in order to design a good expert system. One may add new rules and eliminate old ones proved unsatisfactory. There should be strict rules, and rules with exceptions. Recursivity is seen as a computer program that learns from itself. The statistical regression method may be employed as well to determine a best algorithm of inference.

Non-monotonic reasoning means to make assumptions about things we don't know. Heuristic methods may be involved in order to find successive approximations.

In terms of the previous results, a default neutrosophic logic may be used instead of the normal inference rules. The distribution of possible neutrosophic results serves as an orientating frame for the new results. The flexible, continuously refined, rules obtain iterative and gradual approaches of the result.

A comparison approach is employed to check the result (conclusion) p by studying the opposite of this: what would happen if a non-p conclusion occurred? The inconsistence of information shows up in the result, if not eliminated from the beginning. The data bases should be stratified. There exist methods to construct preferable coherent sub-bases within incoherent bases. In Multi-Criteria Decision one exploits the complementarity of different criteria and the complementarity of various sources.

For example, the Possibility Theory (Zadeh 1978, Dubois, Prade) gives a better approach than the Fuzzy Set Theory (Yager) due to self-improving connectives. The Possibility Theory is proximal to the Fuzzy Set Theory, the difference between these two theories is the way the fusion operators are defined.

One uses the definitions of neutrosophic probability and neutrosophic set operations.
Similarly, there are many ways to construct such connectives according to each particular problem to solve; here we present the easiest ones:

One notes the neutrosophic logical values of the propositions $A_1$ and $A_2$ by
$NL(A_1) = (T_1, I_1, F_1)$ and $NL(A_2) = (T_2, I_2, F_2)$.
For all neutrosophic logical values below: if, after calculations, one obtains numbers $< 0$ or $> 1$, one replaces them by $^-0$ or $1^+$ respectively.

### 2.3.1 Negation:
$NL(\neg A_1) = (\{1^+\} \ominus T_1, \{1^+\} \ominus I_1, \{1^+\} \ominus F_1)$.

### 2.3.2 Conjunction:
$NL(A_1 \wedge A_2) = (T_1 \odot T_2, I_1 \odot I_2, F_1 \odot F_2)$.
(And, in a similar way, generalized for n propositions.)

### 2.3.3 Weak or inclusive disjunction:
$NL(A_1 \vee A_2) = (T_1 \oplus T_2 \ominus T_1 \odot T_2, I_1 \oplus I_2 \ominus I_1 \odot I_2, F_1 \oplus F_2 \ominus F_1 \odot F_2)$.
(And, in a similar way, generalized for n propositions.)

### 2.3.4 Strong or exclusive disjunction:
$NL(A_1 \underline{\vee} A_2) =$
$(T_1 \odot (\{1\} \ominus T_2) \oplus T_2 \odot (\{1\} \ominus T_1) \ominus T_1 \odot T_2 \odot (\{1\} \ominus T_1) \odot (\{1\} \ominus T_2),$
$I_1 \odot (\{1\} \ominus I_2) \oplus I_2 \odot (\{1\} \ominus I_1) \ominus I_1 \odot I_2 \odot (\{1\} \ominus I_1) \odot (\{1\} \ominus I_2),$
$F_1 \odot (\{1\} \ominus F_2) \oplus F_2 \odot (\{1\} \ominus F_1) \ominus F_1 \odot F_2 \odot (\{1\} \ominus F_1) \odot (\{1\} \ominus F_2))$.
(And, in a similar way, generalized for n propositions.)

### 2.3.5 Material conditional (implication):
$NL(A_1 \mapsto A_2) = (\{1^+\} \ominus T_1 \oplus T_1 \odot T_2, \{1^+\} \ominus I_1 \oplus I_1 \odot I_2, \{1^+\} \ominus F_1 \oplus F_1 \odot F_2)$.

### 2.3.6 Material biconditional (equivalence):
$NL(A_1 \leftrightarrow A_2) = ((\{1^+\} \ominus T_1 \oplus T_1 \odot T_2) \odot (\{1^+\} \ominus T_2 \oplus T_1 \odot T_2),$
$(\{1^+\} \ominus I_1 \oplus I_1 \odot I_2) \odot (\{1^+\} \ominus I_2 \oplus I_1 \odot I_2),$
$(\{1^+\} \ominus F_1 \oplus F_1 \odot F_2) \odot (\{1^+\} \ominus F_2 \oplus F_1 \odot F_2))$.

### 2.3.7 Sheffer's connector:
$NL(A_1 | A_2) = NL(\neg A_1 \vee \neg A_2) = (\{1^+\} \ominus T_1 \odot T_2, \{1^+\} \ominus I_1 \odot I_2, \{1^+\} \ominus F_1 \odot F_2)$.

### 2.3.8 Peirce's connector:
$NL(A_1 \downarrow A_2) = NL(\neg A_1 \wedge \neg A_2) =$
$= ( (\{1^+\} \ominus T_1) \odot (\{1^+\} \ominus T_2), (\{1^+\} \ominus I_1) \odot (\{1^+\} \ominus I_2), (\{1^+\} \ominus F_1) \odot (\{1^+\} \ominus F_2) ).$

### 2.4 Generalizations:
When the sets are reduced to an element only respectively, then
$t\_sup = t\_inf = t$, $i\_sup = i\_inf = i$, $f\_sup = f\_inf = f$,
and $n\_sup = n\_inf = n = t+i+f$)
   Hence, the neutrosophic logic generalizes:
   - the *intuitionistic logic*, which supports incomplete theories (for $0 < n < 1$ and $i=0$, $0 \le t, i, f \le 1$);
   - the *fuzzy logic* (for $n = 1$ and $i = 0$, and $0 \le t, i, f \le 1$);
   from "CRC Concise Concise Encyclopedia of Mathematics", by Eric W. Weisstein, 1998, the fuzzy logic is "an extension of two-valued logic such that statements need not to be True or False, but may have a degree of truth between 0 and 1";
   - the *Boolean logic* (for $n = 1$ and $i = 0$, with t, f either 0 or 1);
   - the *multi-valued logic* (for $0 \le t, i, f \le 1$);
   definition of <many-valued logic> from "The Cambridge Dictionary of Philosophy", general editor Robert Audi, 1995, p. 461: "propositions may take many values beyond simple truth and falsity, values functionally determined by the values of their components"; Lukasiewicz considered three values (1, 1/2, 0). Post considered m values, etc. But they varied in between 0 and 1 only. In the neutrosophic logic a proposition may take values even greater than 1 (in percentage greater than 100%) or less than 0.
   - the *paraconsistent logic* (for $n > 1$ and $i = 0$, with both t, f < 1);
   - the *dialetheism*, which says that some contradictions are true (for $t = f = 1$ and $i = 0$; some paradoxes can be denoted this way too);
   - the *faillibilism*, which says that uncertainty belongs to every proposition (for $i > 0$);
   - the *paradoxist logic*, based on paradoxes ($i > 1$);
   - the *pseudoparadoxist logic*, based on pseudoparadoxes ($0 < i < 1$, $t + f > 1$);
   - the *tautologic logic*, based on tautologies ($i < 0$, $t > 1$).
   Compared with all other logics, the neutrosophic logic and intuitionistic fuzzy logic introduce a percentage of "indeterminacy" - due to unexpected parameters hidden in some propositions, or unknowness, or God's will, but only neutrosophic logic let each component t, i, f be even boiling *over 1* (overflooded) or freezing *under 0* (underdried): to be able to make distinction between relative truth and absolute truth, and between relative falsity and absolute falsity.
   For example: in some tautologies $t > 1$, called "overtrue". Similarly, a proposition may be "overindeterminate" (for $i > 1$, in some paradoxes), "overfalse" (for $f > 1$, in some unconditionally false propositions); or "undertrue" (for $t < 0$, in some unconditionally false propositions), "underindeterminate" (for $i < 0$, in some unconditionally true or false propositions), "underfalse" (for $f < 0$, in some unconditionally true propositions).
   This is because we should make a distinction between unconditionally true ($t > 1$, and $f < 0$ or $i < 0$) and conditionally true propositions ($t \le 1$, and $f \le 1$ or $i \le 1$).

### 3  NEUTROSOPHIC SET:

### 3.1 Definition:
Let T, I, F be real standard or non-standard subsets of $]^-0, 1^+[$,
  with   sup T = t_sup, inf T = t_inf,
         sup I = i_sup, inf I = i_inf,
         sup F = f_sup, inf F = f_inf,
  and    n_sup = t_sup+i_sup+f_sup,
         n_inf = t_inf+i_inf+f_inf.

Let U be a universe of discourse, and M a set included in U. An element x from U is noted with respect to the set M as x(T, I, F) and belongs to M in the following way:
it is t% true in the set, i% indeterminate (unknown if it is) in the set, and f% false, where t varies in T, i varies in I, f varies in F.

*Statically T, I, F are subsets, but dynamically the components T, I, F are set-valued vector functions/operators depending on many parameters*, such as: time, space, etc. (some of them are hidden parameters, i.e. unknown parameters).

### 3.2 General Examples:
Let A and B be two neutrosophic sets.
One can say, by language abuse, that any element neutrosophically belongs to any set, due to the percentages of truth/indeterminacy/falsity involved, which varies between 0 and 1 or even less than 0 or greater than 1.
Thus: x(50,20,30) belongs to A (which means, with a probability of 50% x is in A, with a probability of 30% x is not in A, and the rest is undecidable); or y(0,0,100) belongs to A (which normally means y is not for sure in A); or z(0,100,0) belongs to A (which means one does know absolutely nothing about z's affiliation with A).
More general, x( (20-30), (40-45)∪[50-51], {20,24,28} ) belongs to the set A, which means:
- with a probability in between 20-30% x is in A (one cannot find an exact approximate because of various sources used);
- with a probability of 20% or 24% or 28% x is not in A;
- the indeterminacy related to the appurtenance of x to A is in between 40-45% or between 50-51% (limits included);
The subsets representing the appurtenance, indeterminacy, and falsity may overlap, and n_sup = 30+51+28 > 100 in this case.

### 3.3 Physics Examples:
a)  For example the Schrodinger's Cat Theory says that the quantum state of a photon can basically be in more than one place in the same time, which translated to the neutrosophic set means that an element (quantum state) belongs and does not belong to a set (one place) in the same time; or an element (quantum state) belongs to two different sets (two different places) in the same time. It is a question of "alternative worlds" theory very well represented by the neutrosophic set theory.
In Schroedinger's Equation on the behavior of electromagnetic waves and "matter waves" in quantum theory, the wave function Psi that describes the superposition of possible states may be simulated by a neutrosophic function, i.e. a function whose values are not unique for each argument from the domain of definition (the vertical line test fails, intersecting the graph in more points).

Don't we better describe, using the attribute "neutrosophic" than "fuzzy" or any others, a quantum particle that neither exists nor non-exists?
b)  How to describe a particle $\zeta$ in the infinite micro-universe that belongs to two distinct places $P_1$ and $P_2$ in the same time? $\zeta \in P_1$ and $\zeta \notin P_1$ as a true contradiction, or $\zeta \in P_1$ and $\zeta \in \neg P_1$.

### 3.4 Philosophical Examples:
Or, how to calculate the truth-value of Zen (in Japanese) / Chan (in Chinese) doctrine philosophical proposition: the present is eternal and comprises in itself the past and the future?
In Eastern Philosophy the contradictory utterances form the core of the Taoism and Zen/Chan (which emerged from Buddhism and Taoism) doctrines.
How to judge the truth-value of a metaphor, or of an ambiguous statement, or of a social phenomenon which is positive from a standpoint and negative from another standpoint?

There are many ways to construct them, in terms of the practical problem we need to simulate or approach.  Below there are mentioned the easiest ones:

### 3.5 Application:
A cloud is a neutrosophic set, because its borders are ambiguous, and each element (water drop) belongs with a neutrosophic probability to the set (e.g. there are a kind of separated water drops, around a compact mass of water drops, that we don't know how to consider them: in or out of the cloud).
Also, we are not sure where the cloud ends nor where it begins, neither if some elements are or are not in the set.  That's why the percent of indeterminacy is required and the neutrosophic probability (using subsets - not numbers - as components) should be used for better modeling: it is a more organic, smooth, and especially accurate estimation.  Indeterminacy is the zone of ignorance of a proposition's value, between truth and falsehood.

### 3.6 Neutrosophic Set Operations:
One notes, with respect to the sets A and B over the universe U,
   $x = x(T_1, I_1, F_1) \in A$ and $x = x(T_2, I_2, F_2) \in B$, by mentioning x's *neutrosophic membership appurtenance*.
And, similarly,  $y = y(T', I', F') \in B$.
For all neutrosophic set operations: if, after calculations, one obtains numbers $< 0$ or $> 1$, one replaces them by $^-0$ or $1^+$ respectively.

### 3.6.1  Complement of A:
If $x( T_1, I_1, F_1 ) \in A$,
then $x( \{1^+\}\ominus T_1, \{1^+\}\ominus I_1, \{1^+\}\ominus F_1 ) \in C(A)$.

### 3.6.2  Intersection:
If $x( T_1, I_1, F_1 ) \in A$, $x( T_2, I_2, F_2 ) \in B$,
then $x( T_1\odot T_2, I_1\odot I_2, F_1\odot F_2 ) \in A \cap B$.

### 3.6.3  Union:
If $x( T_1, I_1, F_1 ) \in A$, $x( T_2, I_2, F_2 ) \in B$,
then $x( T_1\oplus T_2\ominus T_1\odot T_2, I_1\oplus I_2\ominus I_1\odot I_2, F_1\oplus F_2\ominus F_1\odot F_2 ) \in A \cup B$.

### 3.6.4 Difference:
If $x(T_1, I_1, F_1) \in A$, $x(T_2, I_2, F_2) \in B$,
then $x(T_1 \ominus T_1 \odot T_2, I_1 \ominus I_1 \odot I_2, F_1 \ominus F_1 \odot F_2) \in A \setminus B$,
because $A \setminus B = A \cap C(B)$.

### 3.6.5 Cartesian Product:
If $x(T_1, I_1, F_1) \in A$, $y(T', I', F') \in B$,
then $(x(T_1, I_1, F_1), y(T', I', F')) \in A \times B$.

### 3.6.6 M is a subset of N:
If $x(T_1, I_1, F_1) \in M \Rightarrow x(T_2, I_2, F_2) \in N$,
where $\inf T_1 \leq \inf T_2$, $\sup T_1 \leq \sup T_2$, and $\inf F_1 \geq \inf F_2$, $\sup F_1 \geq \sup F_2$.

### 3.7 Generalizations and Comments:

From the intuitionistic logic, paraconsistent logic, dialetheism, faillibilism, paradoxes, pseudoparadoxes, and tautologies we transfer the "adjectives" to the sets, i.e. to intuitionistic set (set incompletely known), paraconsistent set, dialetheist set, faillibilist set (each element has a percenatge of indeterminacy), paradoxist set (an element may belong and may not belong in the same time to the set), pseudoparadoxist set, and tautologic set respectively.

Hence, the neutrosophic set generalizes:
- the *intuitionistic set*, which supports incomplete set theories (for $0 < n < 1$ and $i = 0$, $0 \leq t, i, f \leq 1$) and incomplete known elements belonging to a set;
- the *fuzzy set* (for $n = 1$ and $i = 0$, and $0 \leq t, i, f \leq 1$);
- the *classical set* (for $n = 1$ and $i = 0$, with t, f either 0 or 1);
- the *paraconsistent set* (for $n > 1$ and $i = 0$, with both $t, f < 1$);
- the *faillibilist set* ($i > 0$);
- the *dialetheist set*, which says that the intersection of some disjoint sets is not empty (for $t = f = 1$ and $i = 0$; some paradoxist sets can be denoted this way too);
- the *paradoxist set* ($i > 1$);
- the *pseudoparadoxist set* ($0 < i < 1$, $t + f > 1$);
- the *tautological set* ($i < 0$).

Compared with all other types of sets, in the neutrosophic set each element has three components which are subsets (not numbers as in fuzzy set) and considers a subset, similarly to intuitionistic fuzzy set, of "indeterminacy" - due to unexpected parameters hidden in some sets, and let the superior limits of the components to even boil *over 1* (overflooded) and the inferior limits of the components to even freeze *under 0* (underdried).

For example: an element in some tautological sets may have $t > 1$, called "overincluded". Similarly, an element in a set may be "overindeterminate" (for $i > 1$, in some paradoxist sets), "overexcluded" (for $f > 1$, in some unconditionally false appurtenances); or "undertrue" (for $t < 0$, in some unconditionally false appurtenances), "underindeterminate" (for $i < 0$, in some unconditionally true or false appurtenances), "underfalse" (for $f < 0$, in some unconditionally true appurtenances).

This is because we should make a distinction between unconditionally true (t > 1, and f < 0 or i < 0) and conditionally true appurtenances (t ≤ 1, and f ≤ 1 or i ≤ 1).

In a *rough set* RS, an element on its boundary-line cannot be classified neither as a member of RS nor of its complement with certainty. In the neutrosophic set a such element may be characterized by x(T, I, F), with corresponding set-values for T, I, F ⊆ ]$^-$0, 1$^+$[.

## 4. NEUTROSOPHIC PROBABILITY:

### 4.1 Definition:
Let T, I, F be real standard or non-standard subsets included in ]$^-$0, 1$^+$[,
with sup T = t_sup, inf T = t_inf,
sup I = i_sup, inf I = i_inf,
sup F = f_sup, inf F = f_inf,
and n_sup = t_sup+i_sup+f_sup,
n_inf = t_inf+i_inf+f_inf.

The *neutrosophic probability* is a generalization of the classical probability and imprecise probability in which the chance that an event A occurs is t% true - where t varies in the subset T, i% indeterminate - where i varies in the subset I, and f% false - where f varies in the subset F. *Statically T, I, F are subsets, but dynamically the components T, I, F are set-valued vector functions/operators depending on many parameters*, such as: time, space, etc. (some of them are hidden parameters, i.e. unknown parameters).

In classical probability n_sup ≤ 1, while in neutrosophic probability n_sup ≤ 3$^+$.

In imprecise probability: the probability of an event is a subset T ⊂ [0, 1], not a number p ∈ [0, 1], what's left is supposed to be the opposite, subset F (also from the unit interval [0, 1]); there is no indeterminate subset I in imprecise probability.

One notes NP(A) = (T, I, F), a triple of sets.

### 4.2 Neutrosophic Probability Space:
The universal set, endowed with a neutrosophic probability defined for each of its subset, forms a neutrosophic probability space.

Let A and B be two neutrosophic events, and NP(A) = ($T_1$, $I_1$, $F_1$), NP(B) = ($T_2$, $I_2$, $F_2$) their neutrosophic probabilities. Then we define:
($T_1$, $I_1$, $F_1$) ⊞ ($T_2$, $I_2$, $F_2$) = ($T_1 \oplus T_2$, $I_1 \oplus I_2$, $F_1 \oplus F_2$),
($T_1$, $I_1$, $F_1$) ⊟ ($T_2$, $I_2$, $F_2$) = ($T_1 \ominus T_2$, $I_1 \ominus I_2$, $F_1 \ominus F_2$),
($T_1$, $I_1$, $F_1$) ⊡ ($T_2$, $I_2$, $F_2$) = ($T_1 \odot T_2$, $I_1 \odot I_2$, $F_1 \odot F_2$).

NP(A∩B) = NP(A) ⊡ NP(B);
NP(¬A) = {1$^+$} ⊟ NP(A), [this second axiom may be replaced, in specific applications, with NP(¬A) = ($F_1$, $I_1$, $T_1$)];
NP(A∪B) = NP(A) ⊞ NP(B) ⊟ NP(A) ⊡ NP(B).

Neutrosophic probability is also a non-additive probability, like the classical probability, but even for independent events, i.e. $P(A \cup B) \neq P(A)+P(B)$. We have equality only when A or B are impossible events.

A probability-function P is called additive if $P(A \cup B) = P(A)+P(B)$, sub-additive if $P(A \cup B) \leq P(A)+P(B)$, and super-additive if $P(A \cup B) \geq P(A)+P(B)$.

In the Dempster-Shafer Theory $P(A) + P(\neg A)$ may be $\neq 1$, in neutrosophic probability almost all the time $P(A) + P(\neg A) \neq 1$.

1. NP(impossible event) = $(T_{imp}, I_{imp}, F_{imp})$,
where sup $T_{imp} \leq 0$, inf $F_{imp} \geq 1$; no restriction on $I_{imp}$.
   NP(sure event) = $(T_{sur}, I_{sur}, F_{sur})$,
where inf $T_{sur} \geq 1$, sup $F_{sur} \leq 0$.
   NP(totally indeterminate event) = $(T_{ind}, I_{ind}, F_{ind})$;
where inf $I_{ind} \geq 1$; no restrictions on $T_{ind}$ or $F_{ind}$.
2. NP(A) $\in$ {(T, I, F), where T, I, F are real standard or non-standard subsets included in $]^-0, 1^+[$ that may overlap}.
3. NP(A$\cup$B) = NP(A) ⊞ NP(B) ⊟ NP(A$\cap$B).
4. NP(A) = {1} ⊟ NP($\neg$A).

### 4.3 Applications:

#1. From a pool of refugees, waiting in a political refugee camp in Turkey to get the American visa, a% have the chance to be accepted - where a varies in the set A, r% to be rejected - where r varies in the set R, and p% to be in pending (not yet decided) - where p varies in P.

Say, for example, that the chance of someone Popescu in the pool to emigrate to USA is (between) 40-60% (considering different criteria of emigration one gets different percentages, we have to take care of all of them), the chance of being rejected is 20-25% or 30-35%, and the chance of being in pending is 10% or 20% or 30%. Then the neutrosophic probability that Popescu emigrates to the Unites States is

   NP(Popescu) = ( (40-60), (20-25)U(30-35), {10,20,30} ), closer to the life.

This is a better approach than the classical probability, where $40 \leq P(Popescu) \leq 60$, because from the pending chance - which will be converted to acceptance or rejection - Popescu might get extra percentage in his will to emigration,
and also the superior limit of the subsets sum
   60+35+30 > 100
and in other cases one may have the inferior sum < 0,
while in the classical fuzzy set theory the superior sum should be 100 and the inferior sum $\geq 0$.
In a similar way, we could say about the element Popescu that
Popescu( (40-60), (20-25)U(30-35), {10,20,30} ) belongs to the set of accepted refugees.

#2. The probability that candidate C will win an election is say 25-30% true (percent of people voting for him), 35% false (percent of people voting against him), and 40% or 41% indeterminate (percent of people not coming to the ballot box, or giving a blank vote - not selecting anyone, or giving a negative vote - cutting all candidates on the list).
Dialectic and dualism don't work in this case anymore.

#3. Another example, the probability that tomorrow it will rain is say 50-54% true according to meteorologists who have investigated the past years' weather, 30 or 34-35% false according to today's very sunny and droughty summer, and 10 or 20% undecided (indeterminate).

#4. The probability that Yankees will win tomorrow versus Cowboys is 60% true (according to their confrontation's history giving Yankees' satisfaction), 30-32% false (supposing Cowboys are actually up to the mark, while Yankees are declining), and 10 or 11 or 12% indeterminate (left to the hazard: sickness of players, referee's mistakes, atmospheric conditions during the game).  These parameters act on players' psychology.

### 4.4 Remarks:

Neutrosophic probability are useful to those events which involve some degree of indeterminacy (unknown) and more criteria of evaluation - as above.  This kind of probability is necessary because it provides a better approach than classical probability to uncertain events.

This probability uses a subset-approximation for the truth-value (like *imprecise probability*), but also subset-approximations for indeterminacy- and falsity-values.

Also, it makes a distinction between "*relative sure event*", event which is sure only in some particular world(s): NP(*rse*) = 1, and "*absolute sure event*", event which is sure in all possible worlds: NP(*ase*) = $1^+$; similarly for "*relative impossible event*" / "*absolute impossible event*", and for "*relative indeterminate event*" / "*absolute indeterminate event*".

In the case when the truth- and falsity-components are complementary, i.e. no indeterminacy and their sum is 100, one falls to the classical probability.  As, for example, tossing dice or coins, or drawing cards from a well-shuffled deck, or drawing balls from an urn.

### 4.5 Generalizations:

An interesting particular case is for n = 1, with $0 \le t, i, f \le 1$, which is closer to the classical probability.

For n = 1 and i = 0, with $0 \le t, f \le 1$, one obtains the classical probability.

If I disappear and F is ignored, while the non-standard unit interval $]^-0, 1^+[$ is transformed into the classical unit interval [0, 1], one gets the imprecise probability.

From the intuitionistic logic, paraconsistent logic, dialetheism, faillibilism, paradoxism, pseudoparadoxism, and tautologism we transfer the  "adjectives" to probabilities, i.e. we define the ***intuitionistic probability*** (when the probability space is incomplete), ***paraconsistent probability***, ***faillibilist probability***, ***dialetheist probability***, ***paradoxist probability***, ***pseudoparadoxist probability***, and ***tautologic probability*** respectively.

Hence, the neutrosophic probability generalizes:
- the *intuitionistic probability*, which supports incomplete (not completely known/determined) probability spaces (for 0 < n < 1 and i = 0, $0 \le t, f \le 1$) or incomplete events whose probability we need to calculate;
- the *classical probability* (for n = 1 and i = 0, and $0 \le t, f \le 1$);
- the *paraconsistent probability* (for n > 1 and i = 0, with both t, f < 1);
 - the *dialetheist probability*, which says that intersection of some disjoint probability spaces is not empty (for t = f = 1 and i = 0; some paradoxist probabilities can be denoted this way);
 - the *faillibilist probability* (for i > 0);
- the *pseudoparadoxism* (for n_sup > 1 or n_inf < 0);
- the *tautologism* (for t_sup > 1).

Compared with all other types of classical probabilities, the neutrosophic probability introduces a percentage of "indeterminacy" - due to unexpected parameters hidden in some probability spaces, and let each component t, i, f be even boiling *over 1* (overflooded) or freezing *under 0* (underdried).

For example: an element in some tautological probability space may have t > 1, called "overprobable".  Similarly, an element in some paradoxist probability space may be "overindeterminate" (for i > 1), or "overunprobable" (for f > 1, in some unconditionally false appurtenances);  or "underprobable" (for t < 0, in some unconditionally false appurtenances), "underindeterminate" (for i < 0, in some unconditionally true or false appurtenances), "underunprobable" (for f < 0, in some unconditionally true appurtenances).

This is because we should make a distinction between unconditionally true (t > 1, and f < 0 or I < 0) and conditionally true appurtenances (t ≤ 1, and f ≤ 1 or I ≤ 1).

## 5 NEUTROSOPHIC STATISTICS:

Analysis of events, described by the neutrosophic probability, means *neutrosophic statistics*.
This is also a generalization of classical statistics.

In accordance with the development of neutrosophic probability the neutrosophic statistics could be better studied.  Here, above, it is only a definition in order to give scientists an impulse for research.

# Combination of paradoxical sources of information within the neutrosophic framework

Dr. Jean Dezert[*]

**Abstract -** *The recent emergence of Smarandache's logic as foundations for a new general unifying theory for uncertain reasoning is becoming both a new philosophical and mathematical research field and could modify deeply our perception and understanding of our outer and inner worlds in coming years. The ability for neutrosophy to include all existing logics as special cases is undoubtedly appealing. Beside of all potential advantages of neutrosophy to handle antinomies and uncertainties, the current mathematical neutrosophic logic, does not deal directly with the important problem of combination of evidences provided by different bodies of evidence. This paper is the first attempt to develop new foundations for the combination of sources of information in a very general framework where information can be both uncertain and paradoxical. We develop a new rule of combination close to the ad-hoc Dempster-Shafer rule of combination where both conjunctions and disjunctions of assertions are explicitly taking into account in the fusion process. Through several simple examples, we show the efficiency of this new theory of plausible and paradoxical reasoning to solve problems where the Dempster-Shafer theory usually fails. Finally a theoretical bridge between the neutrosophic logic and our new theory is presented, in order to solve the delicate problem of the combination of neutrosophic evidences. The neutrosophic logic seems to be an appealing general framework (prerequesite) for dealing with uncertain and paradoxical sources of information through this new theory.*

## 1 Introduction

The processing of uncertain information has always been a hot topic of research since mainly the 18th century. Up to middle of the 20th century, most theoretical advances have been devoted to the theory of probabilities through the works of eminent mathematicians like J. Bernoulli (1713), A. De Moivre (1718), T. Bayes (1763), P. Laplace (1774), K. Gauss (1823), S. Poisson (1837), I. Todhunter (1873), J. Bertrand (1889), E. Borel (1909), R. Fisher (1930), F. Ramsey (1931), A. Kolmogorov (1933), H. Jeffreys (1939), R. Cox (1946), I. Good (1950), R. Carnap (1950), G. Polya (1954), R. Jeffrey (1957), B. De Finetti (1958), M. Kendall (1963), L. Savage (1967),T. Fine (1973), E. Jaynes (1995) to name just few of them. With the development of the computers, the last half of the 20th century has became very prolific for the development new original theories dealing with uncertainty and imprecise information. Mainly three major theories are available now as alternative of the theory of probabilities for the automatic plausible reasoning in expert systems: the fuzzy set theory developed by L. Zadeh in sixties (1965), the Shafer's theory of evidence in seventies (1976) and the theory of possibilities by D. Dubois and H. Prade in eighties (1985). Only recently a new general and original theory, called neutrosophy, has been developed by F. Smarandache to unify all these existing theories in a common global framework. This paper is focused on the development of a new theory of plausible and paradoxical reasoning within the neutrosophical framework. After a brief presentation of probability and Dempster-Shafer theories in sections 2 and 3, we propose the foundations for a new theory in section 4 and discuss about the justification of our new rule of combination of uncertain and paradoxical sources of evidences. Several examples of our new inference will also been presented. In the last section of this paper, we will show how our theory can serve as theoretical tool for the problem of the combination of neutrosophical evidences.

## 2 The probability theory

Let $\Theta = \{\theta_i, i = 1, \ldots, n\}$ be a finite discrete set of *exhaustive* and *exclusive* hypotheses or outcomes of a random experiment. The probability $P(A)$ of $A \subseteq \Theta$ has been defined and interpreted differently (mainly through the geometrical approach, the subjective approach and the frequency approach [35]) since the seventieth century. The frequency approach $P\{A\}$ is defined as the ratio of the number of possible outcomes for event $A$ to the total number of possible outcomes for space $\Theta$. This is still now the easiest approach to introduce the notion of probability (chance) at a low mathematical level. The foundations of the probability theory with a new interpretation as the logic of science can be found through the works of E.T. Jaynes in [27, 25, 26].

---

[*]Dr. Jean Dezert is partially with Onera, 29 Av. Division Leclerc,92320 Châtillon, France. This personal research work is not supported by Onera and does not reflect the position or policy of Onera. Email: Jean.Dezert@onera.fr

## 2.1 Axiomatic approach of the theory of probabilities

Within the frequency approach of probability, one implicitly assumes that each elementary element of $\Theta$ is equally probable. Hence the definition of probabibility itself turns to fall actually in a vicious circle definition. Moreover the principle of sufficient reason (the hypothesis of equiprobable repartition for elementary components of $\Theta$ in case of no prior information), called also the principle of indifference, has been strongly critized especially for cases involving infinitely many possible outcomes because this can lead to confusing paradoxes. That is why, since the work of A. Kolmogorov in 1933, the axiomatic of the probability theory based on $\sigma$-algebras and measure theory has been definitely adopted. We remind now the four axioms of modern theory of probability:

A1: (Nonnegativity)
$$0 \leq P\{A\} \leq 1 \tag{1}$$

A2: (Unity) Any sure event (the sample space) has probability one
$$P\{\Theta\} = 1 \tag{2}$$

A3: (Finite additivity) If $A_1, \ldots, A_n$ are disjoint events, then
$$P\{A_1 \cup \ldots \cup A_n\} = P\{A_1\} + \ldots + P\{A_n\} \tag{3}$$

A4: (countable additivity) If $A_1, A_2, \ldots$ are disjoint events
$$P\{\bigcup_{i=1}^{\infty} A_i\} = \sum_{i=1}^{\infty} P\{A_i\} \tag{4}$$

## 2.2 Consequences of axioms and bayesian inference

From these axioms, all other probability laws (especially total probability theorem and Bayes's rule) can be derived. In particular,

$$P\{\emptyset\} = 0 \quad \text{and} \quad P\{A^c\} = 1 - P\{A\} \tag{5}$$

$$A \subset B \Rightarrow P\{A\} \leq P\{B\} \tag{6}$$

$$\forall A, B \subset \Theta, \quad P\{A \cup B\} = P\{A\} + P\{B\} - P\{A \cap B\} \tag{7}$$

$$\forall A_1, \ldots, A_n \subset \Theta, \quad P\{A_1 \cup \ldots \cup A_n\} \leq \sum_{i=1}^{n} P\{A_i\} \quad \text{(Boole's inequality)} \tag{8}$$

More precisely, in the general case, one has the Poincaré's equality

$$P\{A_1 \cup \ldots \cup A_n\} = \sum_{i=1}^{n} P\{A_i\} - \sum_{i<j} P\{A_i \cap A_j\} + \ldots$$
$$+ (-1)^{k-1} \sum_{i_1 < \ldots < i_k} P\{A_{i_1} \cap \ldots \cap A_{i_k}\} + \ldots + (-1)^n P\{\bigcap_{i=1,n} A_i\} \tag{9}$$

which can be also written under a more compact form as

$$P\{A_1 \cup \ldots \cup A_n\} = \sum_{\substack{I \subset \{1,\ldots,n\} \\ I \neq \emptyset}} (-1)^{|I|+1} P\{\bigcap_{i \in I} A_i\} \tag{10}$$

The probability of an event $A$ under the condition that event $B$ has occured (with probability $P\{B\} \neq 0$) is called the conditional (or a posteriori) probability of $A$ given $B$ and is defined as

$$P\{A \mid B\} = \frac{P\{A \cap B\}}{P\{B\}} \tag{11}$$

Events $A$ and $B$ are said to be independent if $P\{A \cap B\} = P\{A\}P\{B\}$ or equivalently $P\{A \mid B\} = P\{A\}$ and $P\{B \mid A\} = P\{B\}$.

The probability of any event $B$ can be recovered from any partition (i.e. a set of exhaustive and disjoint events) $A_1, \ldots, A_n$ of sample space $\Theta$ by the total probability theorem

$$P\{B\} = \sum_{i=1}^{n} P\{B \mid A_i\}P\{A_i\} \tag{12}$$

Since $P\{A \mid B\}P\{B\} = P\{A \cap B\} = P\{B \mid A\}P\{A\}$, one gets the famous Baye's formula, also called the Bayesian inference [4]

$$P\{A \mid B\} = \frac{P\{B \mid A\}P\{A\}}{P\{B\}} \tag{13}$$

## 2.3 Bayesian rule of combination

Suppose now that M independent sources of information (bodies of evidence) $\mathcal{B}_1, \ldots, \mathcal{B}_M$ provide $M$ *subjective* probability functions $P_1\{.\}, \ldots, P_M\{.\}$ over the same space $\Theta$, then the optimal bayesian fusion rule is obtained as follows (see [13] for a more general and theoretical justification).

$$P_{1,\ldots,M}\{\theta_i\} \triangleq [P_1 \oplus \ldots \oplus P_M]\{\theta_i\} = \frac{p_i^{1-M} \prod_{m=1,M} P_m\{\theta_i\}}{\sum_{i=1,n} p_i^{1-M} \prod_{m=1,M} P_m\{\theta_i\}} \tag{14}$$

where $p_i$ is the prior probability of $\theta_i$. This bayesian fusion rule of combination is however not defined when the sources are in full contradiction because in such case the normalization constant $\sum_{i=1,n} p_i^{1-M} \prod_{m=1,M} P_m\{\theta_i\}$ is zero. A source $\mathcal{B}_j$ is said to be in full conflict with a source $\mathcal{B}_k$ if, for all $\theta_i$, $P_j\{\theta_i\}P_k\{\theta_i\} = 0$. When the fusion is possible and when all the prior probabilities $p_i$ are unknown, one has then to use the principle of indifference by setting all $p_i = 1/n$ and the bayesian rule of combination reduces to

$$P_{1,\ldots,M}\{\theta_i\} = [P_1 \oplus \ldots \oplus P_M]\{\theta_i\} = \frac{\prod_{m=1,M} P_m\{\theta_i\}}{\sum_{i=1,n} \prod_{m=1,M} P_m\{\theta_i\}} \tag{15}$$

The Bayesian inference (13) can be interpreted as a special case of bayesian rule of combination (15) between two sources of information (prior and posterior information).

# 3 The Dempster-Shafer theory of evidence

We present now the basis of the Dempster-Shafer theory (DST) or the Mathematical Theory of Evidence (MTE) [48, 11] called also sometimes the theory of probable or evidential reasoning. The DST is usually considered as a generalization of the bayesian theory of subjective probability [52] and offers a simple and direct representation of ignorance. The DST has shown its compatibility with the classical probability theory, with boolean logic and has a feasible computational complexity [46] for problems of small dimension. The DST is a powerful theoretical tool which can be applied for the representation of incomplete knowledge, belief updating, and for combination of evidence [42, 17] through the Dempster-Shafer's rule of combination presented in the following. The Dempster-Shafer model of representation and processing of uncertainty has led to a huge number of practical applications in a wide range of domains (technical and medical diagnosis under unreliable measuring devices, information retrieval, integration of knowledge from heterogeneous sources for object identification and tracking, network reliability computation, multisensor image segmentation, autonomous navigation, safety control in large plants, map construction and maintenance, just to mention a few).

## 3.1 Basic probability masses

**Definition**

Let $\Theta = \{\theta_i, i = 1, \ldots, n\}$ be a finite discrete set of *exhaustive* and *exclusive* elements (hypotheses) called elementary elements. $\Theta$ has been called the frame of discernment of hypotheses or universe of discourse by G. Shafer. The cardinality (number of elementary elements) of $\Theta$ is denoted $|\Theta|$. The power set $\mathcal{P}(\Theta)$ of $\Theta$ which is the set of all subsets of $\Theta$ is usually noted $\mathcal{P}(\Theta) = 2^{\Theta}$ because its cardinality is exactly $2^{|\Theta|}$. Any element of $2^{\Theta}$ is then a composite event (disjunction) of the frame of discernment. The DST starts by defining a map associated to a body of evidence $\mathcal{B}$ (source of information), called basic assignment probability (bpa) or information granule $m(.) : 2^{\Theta} \to [0, 1]$ such that

$$m(\emptyset) = 0 \tag{16}$$

$$\sum_{A \in 2^\Theta} m(A) \equiv \sum_{A \subseteq \Theta} m(A) = 1 \tag{17}$$

$m(.)$ represents the strength of some evidence provided by the source of information under consideration. Condition (16) reflects the fact that no belief ought to be committed to $\emptyset$ and condition (17) reflects the convention that one's total belief has measure one [48]. The quantity $m(A)$ is called $A$'s basic probability number or sometimes $A$'s basic mass. $m(A)$ corresponds to the measure of the partial belief that is committed *exactly* to $A$ (degree of truth supported exactly by $A$) by the body of evidence $\mathcal{B}$ but not the total belief committed to $A$. All subsets $A$ for which $m(A) > 0$ are called focal elements of $m$. The set of all focal elements of $m(.)$ is called the core $\mathcal{K}(m)$ of $m$. Note that $m(A_1)$ and $m(A_2)$ can both be 0 even if $m(A_1 \cup A_2) \neq 0$. Even more peculiar, note that $A \subset B \not\Rightarrow m(A) < m(B)$. Hence, the bpa $m(.)$ is in general different from a probability distribution $p(.)$.

**Example**

Consider $\Theta = \{\theta_1, \theta_2, \theta_3\}$, then $2^\Theta = \{\emptyset, \theta_1, \theta_2, \theta_3, \theta_1 \cup \theta_2, \theta_1 \cup \theta_3, \theta_2 \cup \theta_3, \theta_1 \cup \theta_2 \cup \theta_3\}$. An information granule $m(.)$ on this frame of discernment $\Theta$ could be defined as

$$
\begin{aligned}
m(\emptyset) &\triangleq 0 & m(\theta_1 \cup \theta_2 \cup \theta_3) &= 0.05 \\
m(\theta_1) &= 0.40 & m(\theta_1 \cup \theta_2) &= 0.10 \\
m(\theta_2) &= 0.20 & m(\theta_2 \cup \theta_3) &= 0.10 \\
m(\theta_3) &= 0.05 & m(\theta_1 \cup \theta_3) &= 0.10
\end{aligned}
$$

In this particular example $\mathcal{K}(m) = \{\theta_1, \theta_2, \theta_3, \theta_1 \cup \theta_2, \theta_1 \cup \theta_3, \theta_2 \cup \theta_3, \theta_1 \cup \theta_2 \cup \theta_3\}$ and note that $\theta_1 \subset \{\theta_1 \cup \theta_2\}$ with $m(\theta_1) > m(\theta_1 \cup \theta_2)$.

## 3.2 Belief functions

To obtain the measure of the total belief committed to $A \in 2^\Theta$, one must add to $m(A)$ the masses $m(B)$ for all proper subsets $B \subset A$. G. Shafer has defined the belief (credibility) function $\mathrm{Bel}(.) : 2^\Theta \to [0, 1]$ associated with bpa $m(.)$ as

$$\mathrm{Bel}(A) = \sum_{B \subseteq A} m(B) \tag{18}$$

$\mathrm{Bel}(A)$ summarizes all our reasons to believe in $A$ (i.e. the lower probability to believe in $A$). More generally, a belief function $\mathrm{Bel}(.)$ can be characterized without reference to the information granule $m(.)$ if $\mathrm{Bel}(.)$ satisfies the following three conditions

$$\mathrm{Bel}(\Theta) = 1 \tag{19}$$

$$\mathrm{Bel}(\emptyset) = 0 \tag{20}$$

$$\forall n > 0, \forall A_1, \ldots, A_n \subset \Theta, \quad \mathrm{Bel}(A_1 \cup \ldots \cup A_n) \geq \sum_{\substack{I \subset \{1, \ldots, n\} \\ I \neq \emptyset}} (-1)^{|I|+1} \mathrm{Bel}\left(\bigcap_{i \in I} A_i\right) \tag{21}$$

For any given belief function $\mathrm{Bel}(.)$, one can always associate an unique information granule $m(.)$, called the Möbius inverse of belief function [44], and defined by [48]

$$\forall A \subseteq \theta, \quad m(A) = \sum_{B \subseteq A} (-1)^{|A-B|} \mathrm{Bel}(B) \tag{22}$$

The *vacuous belief function* having $\mathrm{Bel}(\Theta) = 1$ but $\mathrm{Bel}(A) = 0$ for all $A \neq \Theta$ describes the full ignorance on the frame of discernment $\Theta$. The corresponding bpa $m_v(.)$ is such that $m_v(\Theta) = 1$ and $m_v(A) = 0$ for all $A \neq \Theta$.

For any given belief function $\mathrm{Bel}(.)$ defined on frame $\Theta$, the following inequality holds

$$\forall A, B \subseteq \Theta, \quad \max(0, \mathrm{Bel}(A) + \mathrm{Bel}(B) - 1) \leq \mathrm{Bel}(A \cap B) \leq \min(\mathrm{Bel}(A), \mathrm{Bel}(B)) \tag{23}$$

**Bayesian belief functions**

Any belief function satisfying $\text{Bel}(\emptyset) = 0$, $\text{Bel}(\Theta) = 1$ and $\text{Bel}(A \cup B) = \text{Bel}(A) + \text{Bel}(B)$ whenever $A, B \subset \Theta$ and $A \cap B = \emptyset$ is called a *Bayesian belief function*. In such case, relation (21) coincides exactly with (10) and a probability function $P(.)$ is only a particular Dempster-Shafer's belief function. In this sense, the Dempster-Shafer theory can be considered as a generalization of the probability theory.

If $\text{Bel}(.)$ is a bayesian belief function, then all focal elements are only single points of $\mathcal{P}(\Theta)$. The basic probability mass assignement $m(.)$ commits a positive number $m(\theta_i)$ only to some elementary $\theta_i \in \Theta$ (possibly all $\theta_i$) and zero to all possible disjunctions of $\theta_1, \ldots, \theta_n$. In other words there exists a bayesian bpa $m(.) : \Theta \to [0, 1]$ such that

$$\sum_{\theta_i \in \Theta} m(\theta_i) = 1 \quad \text{and} \quad \forall A \subseteq \Theta, \quad \text{Bel}(A) = \sum_{\theta_i \in A} m(\theta_i) \tag{24}$$

## 3.3 Plausibility functions

Since the degree of belief $\text{Bel}(A)$ does not reveal to what extent one believes its negation $A^c$, Shafer has introduced the degree of doubt of $A$ as the total belief of $A^c$. The degree of doubt is less useful than the plausibility (credibility) $\text{Pl}(A)$ of $A$ which measures the total probability mass that can move into $A$. $\text{Pl}(A)$ can be interpreted as the *upper probability* of $A$. More precisely $\text{Pl}(A)$ is defined by

$$\text{Pl}(A) \triangleq 1 - \text{Dou}(A) = 1 - \text{Bel}(A^c) = \sum_{B \subseteq \Theta} m(B) - \sum_{B \subseteq A^c} m(B) = \sum_{B \cap A \neq \emptyset} m(B) \tag{25}$$

More generally, the dual of (21) implies

$$\forall n > 0, \forall A_1, \ldots, A_n \subset \Theta, \quad \text{Pl}(A_1 \cap \ldots \cap A_n) \leq \sum_{\substack{I \subset \{1, \ldots, n\} \\ I \neq \emptyset}} (-1)^{|I|+1} \text{Pl}(\bigcup_{i \in I} A_i) \tag{26}$$

The direct comparison of (18) with (25) indicates that

$$\forall A \subseteq \Theta, \quad \text{Bel}(A) \leq \text{Pl}(A) \tag{27}$$

For any given plausibility function $\text{Pl}(.)$ defined on frame $\Theta$, the following inequality holds

$$\forall A, B \subseteq \Theta, \quad \max(\text{Pl}(A), \text{Pl}(B)) \leq \text{Pl}(A \cup B) \leq \min(1, \text{Pl}(A) + \text{Pl}(B)) \tag{28}$$

Let $\Theta$ be a given frame of discernment and $m(.)$ a general bpa (neither a vacuous bpa, nor a bayesian bpa) provided by a body of evidence, then it is always possible to build the following *pignistic* probability [63] (bayesian belief function) by choosing $\forall \theta_i \in \Theta, P\{\theta_i\} = \sum_{B \subseteq \Theta | \theta_i \in B} \frac{1}{|B|} m(B)$. In such case, one always has

$$\forall A \subseteq \Theta, \quad \text{Bel}(A) \leq [P(A) = \sum_{\theta_i \in A} P\{\theta_i\}] \leq \text{Pl}(A) \tag{29}$$

Since $\text{Bel}(A)$ summarizes all our reasons to believe in $A$ and $\text{Pl}(A)$ expresses how much we should believe in $A$ if all currently unknown were to support $A$, the true belief in $A$ is somewhere in the interval $[\text{Bel}(A), \text{Pl}(A)]$. Now suppose that the true value of a parameter under consideration is known with some uncertainty $[\text{Bel}(A), \text{Pl}(A)] \subseteq [0, 1]$, then its corresponding bpa $m(A)$ can always be constructed by choosing

$$m(A) = \text{Bel}(A) \quad m(A \cup A^c) = \text{Pl}(A) - \text{Bel}(A) \quad m(A^c) = 1 - \text{Pl}(A) \tag{30}$$

## 3.4 The Dempster's rule of combination

Glenn Shafer has proposed the ad-hoc Dempster's rule of combination (orthogonal summation), symbolized by the operator $\oplus$, to combine two so-called distinct bodies of evidences $\mathcal{B}_1$ and $\mathcal{B}_2$ over the same frame of discernment $\Theta$. Let $\text{Bel}_1(.)$ and $\text{Bel}_2(.)$ be two belief functions over the same frame of discernment $\Theta$ and $m_1(.)$ and $m_2(.)$ their corresponding bpa masses. The combined global belief function $\text{Bel}(.) = \text{Bel}_1(.) \oplus \text{Bel}_2(.)$ is obtained from the combination of the information granules $m_1(.)$ and $m_2(.)$ as follows: $m(\emptyset) = 0$ and for any $C \neq \emptyset$ and $C \subseteq \Theta$,

$$m(C) \triangleq [m_1 \oplus m_2](C) = \frac{\sum_{A \cap B = C} m_1(A) m_2(B)}{\sum_{A \cap B \neq \emptyset} m_1(A) m_2(B)} = \frac{\sum_{A \cap B = C} m_1(A) m_2(B)}{1 - \sum_{A \cap B = \emptyset} m_1(A) m_2(B)} \tag{31}$$

The summation notation $\sum_{A \cap B = C}$ must be interpreted as the sum over all $A, B \subseteq \Theta$ such that $A \cap B = C$ (interpretation for other summation notations follows directly by analogy). The orthogonal sum $m(.)$ is a proper basic probability assignment if $K \triangleq 1 - \sum_{A \cap B = \emptyset} m_1(A)m_2(B) \neq 0$. If $K = 0$, which means $\sum_{A \cap B = \emptyset} m_1(A)m_2(B) = 1$ then orthogonal sum $m(.)$ *does not exist* and the bodies of evidences $\mathcal{B}_1$ and $\mathcal{B}_2$ are said to be totally (flatly) contradictory or in *full contradiction*. Such case arises whenever the cores of $\text{Bel}_1(.)$ and $\text{Bel}_2(.)$ are disjoint or equivalently when there exists $A \subset \Theta$ such that $\text{Bel}_1(A) = 1$ and $\text{Bel}_2(A^c) = 1$. The same problem of existence has already been pointed out previously in the presentation of the optimal Bayesian fusion rule.

The quantity $\log 1/K$ is called the *weight of conflict* between the the bodies of evidences $\mathcal{B}_1$ and $\mathcal{B}_2$. It is easy to show that the Dempster's rule of combination is commutative ($m_1 \oplus m_2 = m_2 \oplus m_1$) and associative ($[m_1 \oplus m_2] \oplus m_3 = m_1 \oplus [m_2 \oplus m_3]$). The vacuous belief function such that $m_v(\Theta) = 1$ and $m_v(A) = 0$ for $A \neq \Theta$ is the identity element for $\oplus$ binary operator, i.e. $m_v \oplus m = m \oplus m_v \equiv m$. If $\text{Bel}_1(.)$ and $\text{Bel}_2(.)$ are two combinable belief functions and if $\text{Bel}_1(.)$ is Bayesian, then $\text{Bel}_1 \oplus \text{Bel}_2$ is a bayesian belief function.

This rule of combination, initially proposed by G. Shafer without a strong theoretical justification (it's only an ''ad-hoc justification''), has been criticized in the past decades by many disparagers of this theory. Nowadays, this rule of combination has however been fully justified by the axiomatic of the transferable belief model developed by Ph. Smets in [59, 60, 61, 63]. We mention the fact that such theoretical justification had been already attempted by Cheng and Kashyap in [7]. Discussions on justifications and interpretations of the DST and the Dempster's rule of combination can be found in [16, 30, 31, 32, 42, 45, 68]. An interesting discussion on the justification of Dempster's rule of combination from the information entropy viewpoint based on the measurement projection and balance principles can be found in [66]. Connection of the DST with the Fuzzy Set Theory can be found in [5, 64]. The relationship between experimental observations and the DST belief functions is currently a hot topic of research. Several models have been developed for fitting belief functions with experimental data. A very recent detailed presentation and discussion on this problem can be found in [69].

We can see a very close similitude between the Dempster's rule and the optimal bayesian fusion rule (15). Actually these two rules coincides exactly when $m_1(.)$ and $m_2(.)$ become bayesian basic probability assignments and if we accept the principle of indifference within the optimal Bayesian fusion rule.

The complexity of DS rule of combination is important in general with large frames of discernment since the computational burden of finding all pairs $A$ and $B$ of subsets of $\Theta$ such that $A \cap B = C$ is $o(2^{|\Theta|-|C|} \times 2^{|\Theta|-|C|})$ which is a large number. For example, if $|\Theta| = 10$ and $|C| = 2$, we will have $o(2^{16}) = o(65536)$ tests to do to find $\{A \cap B | A \cap B = C\}$

**A simple example of the Dempster's rule of combination**

Consider the simple frame of discernment $\Theta = \{S(\text{unny}), R(\text{ainy})\}$ about the true nature of the weather at a given location $L$ for the next day and let consider two independent bodies of evidence $\mathcal{B}_1$ and $\mathcal{B}_2$ providing the following weather forecasts at $L$

$$m_1(S) = 0.80 \qquad m_1(R) = 0.12 \qquad m_1(S \cup R) = 0.08$$
$$m_2(S) = 0.90 \qquad m_2(R) = 0.02 \qquad m_2(S \cup R) = 0.08$$

Applying Dempster's rule of combination yields the following result

$$m(S) = (m_1 \oplus m_2)(S) = (0.72 + 0.072 + 0.064)/(1 - 0.108 - 0.016) \approx 0.9772$$
$$m(R) = (m_1 \oplus m_2)(R) = (0.0024 + 0.0096 + 0.0016)/(1 - 0.108 - 0.016) \approx 0.0155$$
$$m(S \cup R) = (m_1 \oplus m_2)(S \cup R) = 0.0064/(1 - 0.108 - 0.016) \approx 0.0073$$

Hence, in this example, the fusion of the two sources of evidence reinforces the belief that tomorrow will be a sunny day at location $L$ (assuming that both bodies of evidence are equally reliable).

**Another simple but disturbing example**

L. Zadeh has given the following example of a use of the Dempster's rule which shows an unexpected result. Two doctors examine a patient and agree that it suffers from either meningitis (M), concussion (C) or brain tumor (T). Thus $\Theta = \{M, C, T\}$. Assume that the doctors agree in their low expectation of a tumor, but disagree in likely cause and provide the following diagnosis

$$m_1(M) = 0.99 \qquad m_1(T) = 0.01$$
$$m_2(C) = 0.99 \qquad m_2(T) = 0.01$$

If we now combine beliefs using Dempster's rule of combination, one gets the unexpected final conclusion $m(T) = [m_1 \oplus m_2](T) = \frac{0.0001}{1 - 0.0099 - 0.0099 - 0.9801} = 1$ which means that the patient suffers with certainty from brain tumor !!!. This unexpected result arises from the fact that the two bodies of evidence (doctors) agree that patient does not suffer from tumor but are in almost full contradiction for the other causes of the disease. This very simple but practical example shows the limitations of practical use of the DS theory for automated reasoning. Some care must always be taken about the degree of conflict between sources before taking final decision from the result of the Dempster's rule of combination. A justification of non effectiveness of the Dempster's rule in such kind of example based on an information entropy argument has already been presented in [66].

**Conditional Belief functions**

Let $m_B(A) = 1$ if $B \subseteq A$ and $m_B(A) = 0$ for if $B \not\subset A$ (the subset $B$ is the only focal element of $Bel_B$ and its basic probability number is one). Then $Bel_B$ is a belief function that focuses all of the belief on $B$ (note that $Bel_B$ is not in general a Bayesian belief function unless $|B|=1$). If we now consider another belief function Bel over $\Theta$ combinable with $Bel_B$, then the orthogonal sum of Bel with $Bel_B$ denoted as $Bel(. | B) = Bel \oplus Bel_B$ is defined for all $A \subset \Theta$ by [48]

$$Bel(A \mid B) = \frac{Bel(A \cup B^c) - Bel(B^c)}{1 - Bel(B^c)} \quad (32)$$

and

$$Pl(A \mid B) = \frac{Pl(A \cap B)}{Pl(B)} \quad (33)$$

If Bel is a Bayesian belief function, then

$$Bel(A \mid B) = \frac{Bel(A \cap B)}{Bel(B)} = Pl(A \mid B) \quad (34)$$

which coincides exactly with the classical conditional probability P(A|B) defined in (11).

## 4  A new theory for plausible and paradoxical reasoning [DSmT]*

### 4.1 Introduction

As seen in the previous Zadeh's troubling example, the use of the DST must be done only with extreme caution if one has to take a final and important decision with the result of the Dempter's rule of combination. In most (if not all) of practical applications based on the DST, some *ad-hoc* or heuristic recipes are added to the fusion process to correctly manage or reduce the possibility of high degree of conflict between sources. Otherwise, the fusion results lead to unreliable/dangerous conclusion or cannot provide a result at all when the degree of conflict becomes high. Even if the DST provides fruitful results in many applications (mainly in artificial intelligence and systems expert areas) in past decades, we argue that this theory is still too limited because it is based on the following very restrictive constraints :

C1- The DST considers a discrete and finite frame of discernment based on a set of exhaustive and exclusive elementary elements.

C2- The bodies of evidence are assumed independent (each source of information does not take into account the knowledge of other sources) and provide a belief function on the power set $2^\Theta$.

These two constraints are very strong in many practical problems involving uncertain and probable reasoning and dealing with fusion of uncertain and imprecise information. A discussion about this important remark had already been discussed earlier in [33, 34, 47]. In [47], the author proposed a new partitioning management technique to overcome mainly the C2 constraint. The first constraint is very severe actually since it does not allow paradoxes on elements of the frame of discernment $\Theta$. The DST accepts as foundation the commonly adopted principle of the third exclude. Even if at first glance, it makes sense in the traditional classical thought, we can develop a new theory that does not accept the principle of the third exclude and accepts and deals with paradoxes. This is the main purpose of this paper.

The constraint C1 assumes that each elementary hypothesis of the frame of discernment $\Theta$ is finely and precisely well defined and we are able to discriminate between all elementary hypotheses without ambiguity and difficulty. We argue that this constraint is too limited and that it is not always possible in practice to model a frame of discernment satisfying C1 even for some very simple problems where each elementary hypothesis corresponds to a fuzzy concept or attribute. In such cases, the elementary elements of the "frame of discernment" cannot be precisely separated without ambiguity such that no refinement of the frame of discernment satisfying the first constraint is possible. As a simple example, consider an armed robbery situation having a witness and the frame of discernment (associated to the possible size of the thief) having only two elementary imprecise classes $\Theta = \{\theta_1 = \text{small}, \theta_2 = \text{tall}\}$. An investigator asks the witness about the size of the thief and the witness declares that the thief was tall with bpa number $m(\theta_1) = 0.80$, small with bpa number $m(\theta_2) = 0.15$ and is uncertain (either tall or small) with $m(\theta_1 \cup \theta_2) = 0.05$. The investigator will have to deal only with this information although the smallness and the tallness have not been precisely defined. The use of this testimony by the investigator (having in other side extra-information about the thief from other sources of information) to infer on the true size of the thief is delicate especially with the important missing information about the size of the witness (who could be either a basketball player, a dwarf or most probably has a size on the average as you and me - assuming you are neither a dwarf or a basketball player. These both hypotheses are not incompatible actually since some dwarfs really enjoy to play basketball).

---

\* **This has been called Dezert-Smarandache Theory for Plausible and Paradoxical Reasoning** [ref.].

In many situations, we argue that the frame of discernment $\Theta$ can only be described in terms of imprecise elements which cannot be clearly separated and which cannot be considered as fully disjoint and that the refinement of initial frame into a new one satisfying C1 is like a graal quest and cannot be accomplished. Our last remark about C1 constraint concerns the universal nature of the frame of discernment. As shown in our previous simple example, it is clear that, in general, the "same" frame of discernment is interpreted differently by the bodies of evidence. Some subjectivity, or at least some fortuitous biases, on the information provided by a source of information is almost unavoidable, otherwise this would assume (as foundation for the DST) that all bodies of evidence have an objective/universal (possibly uncertain) interpretation or measure of the phenomena under consideration. This vision seems to be too excessive because usually independent bodies of evidence provide their beliefs about some hypotheses only with respect to their own worlds of knowledge and experience. We don't go deeper here in the techniques of refinements and coarsenings of compatible frame of discernments which is a prerequesite to the Dempster's rule of combination. This has already been presented in details in chapter 6 of [48]. We just want to emphaze here that these nice appealing techniques cannot be used at all in all cases where C1 cannot be satisfied and we have more generally to accept the idea to deal with paradoxical information. To convince the reader to accept our radical new way of thought, just think about the true nature of a photon? For experts working in particle physics, photons look like particles, for physicists working in electromagnetic field theory, photons are considered as electromagnetic waves. Both interpretations are true, there is no unicity on the true nature of the photon and actually a photon holds both aspects which appears as a paradoxe for most of human minds. This notion has been accepted in modern physics only with great difficulty and many vigourous discussions about this fundamental question have held at the beginning of the 20th century between all eminent physicists at that time.

The constraint C2 hides a strong difficulty already discussed in the previous paragraph. In order to apply the Dempter's rule of combination of two independent bodies of evidence $\mathcal{B}_1$ and $\mathcal{B}_2$), it is necessary that both frames of discernment $\Theta_1$ and $\Theta_2$ (related to each source $\mathcal{B}_1$ and $\mathcal{B}_2$) have to be compatible and to correspond to the same "universal vision" of the possibilities of the answer of the question under consideration. This constraint is itself very difficult to satisfy actually since each source of information has usually only a personal (and maybe biaised) interpretation of elements of frame of discernment. The belief provided by each local source of information mainly depends on the own knowledge frame of the source without reference to the (inaccessible) absolute truth of the space of possibilities. Therefore, C2 is in many cases also a too strong hypothesis to accept for foundations of a general theory of probable and paradoxical reasoning. A general theory should include the possibility to deal with evidences arising from different sources of information which have no access to absolute interpretation of the elements of the frame of discernment $\Theta$ under consideration. This yields to accept paradoxical information as basis for a new general theory of probable reasoning. Actually we will show in the first example on section 4.3 that paradoxical information resulting of fusion of several bodies of evidence is very informative and can be used to help us to take legitimous final decision.

In other words, this new theory can be interpreted as a general and direct extension of probability theory and the Dempster-Shafer theory in the following sense. Let $\Theta = \{\theta_1, \theta_2\}$ be the simplest frame of discernment involving only two elementary hypotheses (with no more additional assumptions on $\theta_1$ and $\theta_2$), then

- the probability theory deals with basic assignment masses $m(.) \in [0, 1]$ such that
$$m(\theta_1) + m(\theta_2) = 1$$

- the Dempster-Shafer theory extends the probability theory by dealing with basic assignment masses $m(.) \in [0, 1]$ such that
$$m(\theta_1) + m(\theta_2) + m(\theta_1 \cup \theta_2) = 1$$

- our general theory extends the two previous theories by accepting the possibility of paradoxical information and deals with new basic assignment masses $m(.) \in [0, 1]$ such that
$$m(\theta_1) + m(\theta_2) + m(\theta_1 \cup \theta_2) + m(\theta_1 \cap \theta_2) = 1$$

## 4.2 Hyper-power set and general basic probability assignment $m(.)$

### 4.2.1 Hyper-power set definition

Let $\Theta = \{\theta_1, \ldots, \theta_n\}$ be a set of $n$ elementary elements considered as exhaustive which cannot be precisely defined and separated so that no refinement of $\Theta$ in a new larger set $\Theta_{ref}$ of disjoint elementary hypotheses is possible and let's consider the classical set operators $\cup$ (disjunction) and $\cap$ (conjunction). The exhaustive hypothesis about $\Theta$ is not a strong constraint since when $\theta_i, i = 1, n$ does not constitute an exhaustive set of elementary possibilities, we can always add an extra element $\theta_0$ such that $\theta_i, i = 0, n$ describes now an exhaustive set. We will assume therefore, from now on and in the following, that $\Theta$ characterizes an exhaustive frame of discernment. $\Theta$ will be called a *general* frame of discernment in the sequel to emphaze the fact that $\Theta$ does not satisfy the Dempster-Shafer C1 constraint.

The classical power set $\mathcal{P}(\Theta) = 2^\Theta$ has been defined as the set of all proper subsets of $\Theta$ when all elements $\theta_i$ are disjoint. We extend here this notion and define now the "hyper-power" set $D^\Theta$ as the set of all composite possibilities build from $\Theta$ with $\cup$ and $\cap$ operators such that $\forall A \in D^\Theta, B \in D^\Theta, (A \cup B) \in D^\Theta$ and $(A \cap B) \in D^\Theta$. Obviously, one would always have $D^\Theta \subset 2^{\Theta_{ref}}$ if the refined power set $2^{\Theta_{ref}}$ could be defined and accessible which is unfortunately not possible in general as already argued. The cardinality of hyper-power set $D^\Theta$ is majored by $2^{2^n}$ when $\text{Card}(\Theta) = |\Theta| = n$. The generation of hyper-power set $D^\Theta$ corresponds to the famous Dedekind's problem on enumerating the set of monotone Boolean functions (i.e., functions expressible using only AND and OR set operators) [10]. This problem is also related with the Sperner systems [65, 37] based on finite poset (called also as antichains in literature) [8]. The number of antichains on the $n$-set $\Theta$ are equal to the number of monotonic increasing Boolean functions of $n$ variables, and also the number of free distributive lattices with $n$ generators [18, 20, 28, 29, 38, 54]. Determining these numbers is exactly the Dedekind's problem. The choice of letter $D$ in our notation $D^\Theta$ to represent the hyper-power set of $\Theta$ is in honour of the great mathematician R. Dedekind. The general solution of the Dedekind's problem (for $n > 10$) has not been found yet. We just know that the cardinality numbers of $D^\Theta$ follow the integers of the Dedekind's sequence minus one when $\text{Card}(\Theta) = n$ increases.

**Examples**

1. If we consider $\Theta = \{\}$ (empty set) then $D^\Theta = \{\emptyset\}$ and $|D^\Theta| = 1$
2. If we consider $\Theta = \{\theta_1\}$ then $D^\Theta = \{\emptyset, \theta_1\}$ and $|D^\Theta| = 2$
3. If we consider $\Theta = \{\theta_1, \theta_2\}$ then $D^\Theta = \{\emptyset, \theta_1, \theta_2, \theta_1 \cup \theta_2, \theta_1 \cap \theta_2\}$ and $|D^\Theta| = 5$
4. If we consider $\Theta = \{\theta_1, \theta_2, \theta_3\}$ then

$$D^\Theta = \{\emptyset, \theta_1, \theta_2, \theta_3, \theta_1 \cup \theta_2, \theta_1 \cup \theta_3, \theta_2 \cup \theta_3, \theta_1 \cap \theta_2, \theta_1 \cap \theta_3, \theta_2 \cap \theta_3, \theta_1 \cup \theta_2 \cup \theta_3, \theta_1 \cap \theta_2 \cap \theta_3,$$
$$(\theta_1 \cup \theta_2) \cap \theta_3, (\theta_1 \cup \theta_3) \cap \theta_2, (\theta_2 \cup \theta_3) \cap \theta_1, (\theta_1 \cap \theta_2) \cup \theta_3, (\theta_1 \cap \theta_3) \cup \theta_2, (\theta_2 \cap \theta_3) \cup \theta_1,$$
$$(\theta_1 \cup \theta_2) \cap (\theta_1 \cup \theta_3) \cap (\theta_2 \cup \theta_3)\}$$

and $|D^\Theta| = 19$

It is not difficult, although tedious, to check that $\forall A \in D^\Theta, B \in D^\Theta, (A \cup B) \in D^\Theta$ and $(A \cap B) \in D^\Theta$. The extension to larger frame of discernment is possible but requires a higher computational burden. The general and direct analytic computation of $|D^\Theta|$ for a $n$-set $\Theta$ with $n > 10$ is not known and is still under investigations in the mathematical community. Cardinality numbers $|D^\Theta|$ follow the Dedekind's sequence (minus one), i.e. $1, 2, 5, 19, 167, 7580, 7828353, \ldots$ when $\text{Card}(\Theta) = n = 0, 1, 2, 3, 4, 5, 6, \ldots$.

## 4.3 General basic assignment numbers

### 4.3.1 Definition

Let $\Theta$ be a *general* frame of discernment of the problem under consideration. We define a map $m(.) : D^\Theta \to [0, 1]$ associated to a given body of evidence $\mathcal{B}$ which can support paradoxical information, as follows

$$m(\emptyset) = 0 \tag{35}$$

$$\sum_{A \in D^\Theta} m(A) = 1 \tag{36}$$

The quantity $m(A)$ is called $A$'s general basic probability number. As in the DST, all subsets $A \in D^\Theta$ for which $m(A) > 0$ are called focal elements of $m(.)$ and the set of all focal elements of $m(.)$ is also called the core $\mathcal{K}(m)$ of $m$. The belief and plausibility functions are defined in the same way as in the DST, i.e.

$$\text{Bel}(A) = \sum_{B \in D^\Theta, B \subseteq A} m(B) \tag{37}$$

$$\text{Pl}(A) = \sum_{B \in D^{\Theta}, B \cap A \neq \emptyset} m(B) \tag{38}$$

Note that, we don't define here explicitly the complementary $A^c$ of a proposition $A$ since $m(A^c)$ cannot be precisely evaluated from $\cup$ and $\cap$ operators on $D^{\Theta}$ since we include the possibility to deal with a complete paradoxical source of information such that $\forall A \in D^{\Theta}, \forall B \in D^{\Theta}, m(A \cap B) > 0$. These definitions are compatible with the DST definitions when the sources of information become uncertain but rational (they do not support paradoxical information). We still have $\forall A \in D^{\Theta}, \text{Bel}(A) \leq \text{Pl}(A)$.

### 4.3.2 Construction of pignistic probabilities from general basic assignment numbers

The construction of a pignistic probability measure from the general basic numbers $m(.)$ over $D^{\Theta}$ with $|\Theta| = n$ is still possible and is given by the general expression of the form

$$\forall i = 1, \ldots, n \qquad P\{\theta_i\} = \sum_{A \in D^{\Theta}} \alpha_{\theta_i}(A) m(A) \tag{39}$$

where $\alpha_{\theta_i}(A) \in [0, 1]$ are weighting coefficients which depend on the inclusion or non-inclusion of $\theta_i$ with respect to proposition $A$. No general analytic expression for $\alpha_{\theta_i}(A)$ has been derived yet even if $\alpha_{\theta_i}(A)$ can be obtained explicitly for simple examples. When general bpa $m(.)$ reduces to classical bpa (i.e. the DS bpa without paradoxe), then $\alpha_{\theta_i}(A) = \frac{1}{|A|}$ if $\theta_i \subseteq A$ and therefore one gets

$$\forall i = 1, \ldots, n \qquad P\{\theta_i\} = \sum_{A \subseteq \Theta | \theta_i \in A} \frac{1}{|A|} m(A) \tag{40}$$

We present now two examples of pignistic probabilities reconstruction from a general and non degenerated bpa $m(.)$ (i.e. $\nexists A \in D^{\Theta}$ with $A \neq \emptyset$ such that $m(A) = 0$) over $D^{\Theta}$.

- Example 1 : If $\Theta = \{\theta_1, \theta_2\}$ then

$$P\{\theta_1\} = m(\theta_1) + \frac{1}{2} m(\theta_1 \cup \theta_2) + \frac{1}{2} m(\theta_1 \cap \theta_2)$$

$$P\{\theta_2\} = m(\theta_2) + \frac{1}{2} m(\theta_1 \cup \theta_2) + \frac{1}{2} m(\theta_1 \cap \theta_2)$$

- Example 2 : If $\Theta = \{\theta_1, \theta_2, \theta_3\}$ then

$$\begin{aligned}
P\{\theta_1\} =& m(\theta_1) + \frac{1}{2} m(\theta_1 \cup \theta_2) + \frac{1}{2} m(\theta_1 \cup \theta_3) + \frac{1}{2} m(\theta_1 \cap \theta_2) + \frac{1}{2} m(\theta_1 \cap \theta_3) \\
&+ \frac{1}{3} m(\theta_1 \cup \theta_2 \cup \theta_3) + \frac{1}{3} m(\theta_1 \cap \theta_2 \cap \theta_3) \\
&+ \frac{1/2 + 1/3}{3} m((\theta_1 \cup \theta_2) \cap \theta_3) + \frac{1/2 + 1/3}{3} m((\theta_1 \cup \theta_3) \cap \theta_2) + \frac{1/2 + 1/2 + 1/3}{3} m((\theta_2 \cup \theta_3) \cap \theta_1) \\
&+ \frac{1/2 + 1/2 + 1/3}{5} m((\theta_1 \cap \theta_2) \cup \theta_3) + \frac{1/2 + 1/2 + 1/3}{5} m((\theta_1 \cap \theta_3) \cup \theta_2) \\
&+ \frac{1 + 1/2 + 1/2 + 1/3}{5} m((\theta_2 \cap \theta_3) \cup \theta_1) + \frac{1/2 + 1/2 + 1/3}{4} m((\theta_1 \cup \theta_2) \cap (\theta_1 \cup \theta_2) \cap (\theta_2 \cup \theta_3))
\end{aligned}$$

$$\begin{aligned}
P\{\theta_2\} =& m(\theta_2) + \frac{1}{2} m(\theta_1 \cup \theta_2) + \frac{1}{2} m(\theta_2 \cup \theta_3) + \frac{1}{2} m(\theta_1 \cap \theta_2) + \frac{1}{2} m(\theta_2 \cap \theta_3) \\
&+ \frac{1}{3} m(\theta_1 \cup \theta_2 \cup \theta_3) + \frac{1}{3} m(\theta_1 \cap \theta_2 \cap \theta_3) \\
&+ \frac{1/2 + 1/3}{3} m((\theta_1 \cup \theta_2) \cap \theta_3) + \frac{1/2 + 1/2 + 1/3}{3} m((\theta_1 \cup \theta_3) \cap \theta_2) + \frac{1/2 + 1/3}{3} m((\theta_2 \cup \theta_3) \cap \theta_1) \\
&+ \frac{1/2 + 1/2 + 1/3}{5} m((\theta_1 \cap \theta_2) \cup \theta_3) + \frac{1 + 1/2 + 1/2 + 1/3}{5} m((\theta_1 \cap \theta_3) \cup \theta_2) \\
&+ \frac{1/2 + 1/2 + 1/3}{5} m((\theta_2 \cap \theta_3) \cup \theta_1) + \frac{1/2 + 1/2 + 1/3}{4} m((\theta_1 \cup \theta_2) \cap (\theta_1 \cup \theta_2) \cap (\theta_2 \cup \theta_3))
\end{aligned}$$

$$P\{\theta_3\} = m(\theta_3) + \frac{1}{2}m(\theta_1 \cup \theta_3) + \frac{1}{2}m(\theta_2 \cup \theta_3) + \frac{1}{2}m(\theta_1 \cap \theta_3) + \frac{1}{2}m(\theta_2 \cap \theta_3)$$

$$+ \frac{1}{3}m(\theta_1 \cup \theta_2 \cup \theta_3) + \frac{1}{3}m(\theta_1 \cup \theta_2 \cup \theta_3)$$

$$+ \frac{1/2+1/2+1/3}{3}m(\theta_1 \cup \theta_2 \cup \theta_3) + \frac{1/2+1/3}{3}m(\theta_1 \cup \theta_3 \cup \theta_2) + \frac{1/2+1/3}{3}m(\theta_2 \cup \theta_3 \cup \theta_1)$$

$$+ \frac{1+1/2+1/2+1/3}{5}m((\theta_1 \cap \theta_2) \cup \theta_3) + \frac{1/2+1/2+1/3}{5}m((\theta_1 \cap \theta_3) \cup \theta_2)$$

$$+ \frac{1/2+1/2+1/3}{5}m(\theta_2 \cup \theta_3) \cup \theta_1 + \frac{1/2+1/2+1/3}{4}m((\theta_1 \cup \theta_2) \cap (\theta_1 \cup \theta_2) \cap (\theta_2 \cup \theta_3))$$

The evaluation of weighting coefficients $\alpha\theta_i(A)$ has been obtained from the geometrical interpretation of the relative contribution of the distinct parts of A with proposition $\theta_i$ under consideration. For example, consider $A = (\theta_1 \cap \theta_2) \cup \theta_3$ which corresponds to the area $a_1 \cup a_2 \cup a_3 \cup a_4 \cup a_5$ on the following Venn diagram.

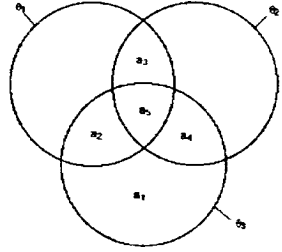

Figure 1: Representation of $A = (\theta_1 \cap \theta_2) \cup \theta_3 \equiv a_1 \cup a_2 \cup a_3 \cup a_4 \cup a_5$

$a_1$ which is shared only by $\theta_3$ will contribute to $\theta_3$ with weight 1; $a_2$ which is shared by $\theta_1$ and $\theta_3$ will contribute to $\theta_3$ with weight 1/2; $a_3$ which is not shared by $\theta_3$ will contribute to $\theta_3$ with weight 0; $a_4$ which is shared by $\theta_2$ and $\theta_3$ will contribute to $\theta_3$ with weight 1/2; $a_5$ which is shared by both $\theta_1$, $\theta_2$ and $\theta_3$ will contribute to $\theta_3$ with weight 1/3. Since moreover, one must have

$\forall A \in D^\Theta$ with $m(A) \neq 0$, $m(A) \neq 0, \sum_{i=1}^{n} \alpha\theta_i(A)m(A) = m(A)$, it is necessary to normalize $\alpha\theta_i(A)$. Therefore $\alpha\theta_i(A)$ is given by

$$\alpha\theta_3(A) = \frac{1+1/2+1/2+1/3}{5} \text{ and similarly } \alpha\theta_2(A) = \frac{1/2+1/2+1/3}{5}, \alpha\theta_1(A) = \frac{1/2+1/2+1/3}{5}$$

All $\alpha\theta_i(A)$, $\forall A \in D^\Theta$ entering in derivation of $P\{\theta_i\}$ can be obtained using similar process.

### 4.4 General rule of combination of paradoxical sources of evidence*

#### 4.4.1 The rule of combination
Let's consider now two distinct (but potentially paradoxical) bodies of evidences $B_1$ and $B_2$ over the same frame of discernment $\theta$ with belief functions $Bel_1(.)$ and $Bel_2(.)$ associated with information granules $m_1(.)$ and $m_2(.)$. The combined global belief function $Bel(.) = Bel_1(.) \oplus Bel_2(.)$ is obtained through the combination of the granules $m_1(.)$ and $m_2(.)$ by the simple rule

$$\forall C \varepsilon D^\theta, m(C) \underline{\underline{\Delta}} [m_1 \oplus m_2](C) = \sum_{A,B\varepsilon D^\theta, A \cap B = C} m_1(A)m_2(B) \quad (41)$$

Since $D^\Theta$ is closed under $\cup$ and $\cap$ operators, this new rule of combination guarantees that $m(.) : D^\Theta \rightarrow [0, 1]$ is a proper general information granule satisfying (35) and (36). The global belief function $Bel(.)$ is then obtained from the granule $m(.)$ through (37). This rule of combination is commutative and associative and can always be used for fusion of paradoxical and/or rational sources of information (bodies of evidence). Obviously, the decision process will have to be made with more caution to take final decision based on the general granule $m()$ when internal paradoxical conflicts arise.

It is important to note that any fusion of sources of information generates either uncertainties, paradoxes or in general both. This is intrinsic to the general fusion process itself. For instance, let's consider the frame of discernment $\Theta = \{\theta_1, \theta_2\}$ and the following very simple examples:

---

*This has been called **Dezert-Smarandache Rule of Combination of Paradoxical Sources of Evidence** [ref.].

- If we consider the two rational information granules

$$m_1(\theta_1) = 0.80 \qquad m_1(\theta_2) = 0.20 \qquad m_1(\theta_1 \cup \theta_2) = 0 \qquad m_1(\theta_1 \cap \theta_2) = 0$$

$$m_2(\theta_1) = 0.90 \qquad m_2(\theta_2) = 0.10 \qquad m_2(\theta_1 \cup \theta_2) = 0 \qquad m_2(\theta_1 \cap \theta_2) = 0$$

then

$$m(\theta_1) = 0.72 \qquad m(\theta_2) = 0.02 \qquad m(\theta_1 \cup \theta_2) = 0 \qquad m(\theta_1 \cap \theta_2) = 0.26$$

- If we consider the two uncertain information granules

$$m_1(\theta_1) = 0.80 \qquad m_1(\theta_2) = 0.15 \qquad m_1(\theta_1 \cup \theta_2) = 0.05 \qquad m_1(\theta_1 \cap \theta_2) = 0$$

$$m_2(\theta_1) = 0.90 \qquad m_2(\theta_2) = 0.05 \qquad m_2(\theta_1 \cup \theta_2) = 0.05 \qquad m_2(\theta_1 \cap \theta_2) = 0$$

then

$$m(\theta_1) = 0.805 \qquad m(\theta_2) = 0.0175 \qquad m(\theta_1 \cup \theta_2) = 0.0025 \qquad m(\theta_1 \cap \theta_2) = 0.175$$

- If we consider the two paradoxical information granules

$$m_1(\theta_1) = 0.80 \qquad m_1(\theta_2) = 0.15 \qquad m_1(\theta_1 \cup \theta_2) = 0 \qquad m_1(\theta_1 \cap \theta_2) = 0.05$$

$$m_2(\theta_1) = 0.90 \qquad m_2(\theta_2) = 0.05 \qquad m_2(\theta_1 \cup \theta_2) = 0 \qquad m_2(\theta_1 \cap \theta_2) = 0.05$$

then

$$m(\theta_1) = 0.72 \qquad m(\theta_2) = 0.0075 \qquad m(\theta_1 \cup \theta_2) = 0 \qquad m(\theta_1 \cap \theta_2) = 0.2725$$

- If we consider the two uncertain and paradoxical information granules

$$m_1(\theta_1) = 0.80 \qquad m_1(\theta_2) = 0.10 \qquad m_1(\theta_1 \cup \theta_2) = 0.05 \qquad m_1(\theta_1 \cap \theta_2) = 0.05$$

$$m_2(\theta_1) = 0.90 \qquad m_2(\theta_2) = 0.05 \qquad m_2(\theta_1 \cup \theta_2) = 0.03 \qquad m_2(\theta_1 \cap \theta_2) = 0.02$$

then

$$m(\theta_1) = 0.789 \qquad m(\theta_2) = 0.0105 \qquad m(\theta_1 \cup \theta_2) = 0.0015 \qquad m(\theta_1 \cap \theta_2) = 0.199$$

Note that this general fusion rule can also be used with intuitionist logic in which the sum of bpa is allowed to be less than one ($\sum m(A) < 1$) and with the paraconsistent logic in which the sum of bpa is allowed to be greater than one ($\sum m(A) > 1$) as well. In such cases, the fusion result does not provide in general $\sum m(A) = 1$. By example, let's consider the fusion of the paraconsistent source $\mathcal{B}_1$ with $m_1(\theta_1) = 0.60$, $m_1(\theta_2) = 0.30$, $m_1(\theta_1 \cup \theta_2) = 0.20$, $m_1(\theta_1 \cap \theta_2) = 0.10$ with the intuitionist source $\mathcal{B}_2$ with $m_2(\theta_1) = 0.50$, $m_2(\theta_2) = 0.20$, $m_2(\theta_1 \cup \theta_2) = 0.10$, $m_2(\theta_1 \cap \theta_2) = 0.10$. In such case, the fusion result of these two sources of information yields the following global paraconsistent bpa $m(.)$

$$m(\theta_1) = 0.46 \qquad m(\theta_2) = 0.13 \qquad m(\theta_1 \cup \theta_2) = 0.02 \qquad m(\theta_1 \cap \theta_2) = 0.47 \qquad \Rightarrow \sum m = 1.08 > 1$$

In practice, for the sake of fair comparison between several alternatives or choices, it is better and simpler to deal with normalized bpa to take a final important decision for the problem under consideration. A nice property of the new rule of combination of non-normalized bpa is its invariance to the pre- or post-normalization process as we will show right now. In the previous example, the post-normalization of bpa $m(.)$ will yield the new bpa $m'(.)$

$$m'(\theta_1) = \frac{0.46}{1.08} \approx 0.426 \qquad m'(\theta_2) = \frac{0.13}{1.08} \approx 0.12 \qquad m'(\theta_1 \cup \theta_2) = \frac{0.02}{1.08} \approx 0.019 \qquad m'(\theta_1 \cap \theta_2) = \frac{0.47}{1.08} \approx 0.435$$

The fusion of pre-normalization of bpa $m_1(.)$ and $m_2(.)$ will yield the same normalized bpa $m'(.)$ since

$$m'_1(\theta_1) = \frac{0.6}{1.2} = 0.50 \qquad m'_1(\theta_2) = \frac{0.3}{1.2} = 0.25 \qquad m'_1(\theta_1 \cup \theta_2) = \frac{0.2}{1.2} \approx 0.17 \qquad m'_1(\theta_1 \cap \theta_2) = \frac{0.1}{1.2} \approx 0.08$$

$$m'_2(\theta_1) = \frac{0.5}{0.9} \approx 0.56 \qquad m'_2(\theta_2) = \frac{0.2}{0.9} \approx 0.22 \qquad m'_2(\theta_1 \cup \theta_2) = \frac{0.1}{0.9} \approx 0.11 \qquad m'_2(\theta_1 \cap \theta_2) = \frac{0.1}{0.9} \approx 0.11$$

$$m'(\theta_1) \approx 0.426 \qquad m'(\theta_2) \approx 0.12 \qquad m'(\theta_1 \cup \theta_2) \approx 0.019 \qquad m'(\theta_1 \cap \theta_2) \approx 0.435$$

It is easy to verify from the general fusion table that the pre or post normalization step yields always the same global normalized bpa even for the general case (when $|\Theta| = n$) because the post-normalization constant $\sum m(A)$ is always equal to the product of the two pre-normalization constants $\sum m_1(A)$ and $\sum m_2(A)$.

### 4.4.2 Justification of the new rule of combination from information entropy

Let's consider two bodies of evidence $\mathcal{B}_1$ and $\mathcal{B}_2$ characterized respectively by their bpa $m_1(.), m_2(.)$ and their cores $\mathcal{K}_1 = \mathcal{K}(m_1), \mathcal{K}_2 = \mathcal{K}(m_2)$. Following Sun's notation [66], each source of information will be denoted

$$\mathcal{B}_1 = \begin{bmatrix} \mathcal{K}_1 \\ m_1 \end{bmatrix} = \begin{bmatrix} f_1^{(1)} & f_2^{(1)} & \dots & f_k^{(1)} \\ m_1(f_1^{(1)}) & m_1(f_2^{(1)}) & \dots & m_1(f_k^{(1)}) \end{bmatrix} \quad (42)$$

$$\mathcal{B}_2 = \begin{bmatrix} \mathcal{K}_2 \\ m_2 \end{bmatrix} = \begin{bmatrix} f_1^{(2)} & f_2^{(2)} & \dots & f_l^{(2)} \\ m_2(f_1^{(2)}) & m_2(f_2^{(2)}) & \dots & m_2(f_l^{(2)}) \end{bmatrix} \quad (43)$$

where $f_i^{(1)}, i = 1, k$ are focal elements of $\mathcal{B}_1$ and $f_j^{(2)}, j = 1, l$ are focal elements of $\mathcal{B}_2$.

Let's consider now the combined information associated with a new body of evidence $\mathcal{B}$ resulting from the fusion of $\mathcal{B}_1$ and $\mathcal{B}_2$ having bpa $m(.)$ with core $\mathcal{K}$. We denote $\mathcal{B}$ as

$$\mathcal{B} \triangleq \mathcal{B}_1 \oplus \mathcal{B}_2 = \begin{bmatrix} \mathcal{K} \\ m \end{bmatrix} = \begin{bmatrix} f_1^{(1)} \cap f_1^{(2)} & f_1^{(1)} \cap f_2^{(2)} & \dots & f_k^{(1)} \cap f_l^{(2)} \\ m(f_1^{(1)} \cap f_1^{(2)}) & m(f_1^{(1)} \cap f_2^{(2)}) & \dots & m(f_k^{(1)} \cap f_l^{(2)}) \end{bmatrix} \quad (44)$$

The fusion of 2 informations granules can be represented with the general table fusion as follows

| $\oplus$ | $m_1(f_1^{(1)})$ | $m_1(f_2^{(1)})$ | ... | $m_1(f_i^{(1)})$ | ... | $m_1(f_k^{(1)})$ |
|---|---|---|---|---|---|---|
| $m_2(f_1^{(2)})$ | $m(f_1^{(1)} \cap f_1^{(2)})$ | $m(f_2^{(1)} \cap f_1^{(2)})$ | ... | $m(f_i^{(1)} \cap f_1^{(2)})$ | ... | $m(f_k^{(1)} \cap f_1^{(2)})$ |
| $m_2(f_2^{(2)})$ | $m(f_1^{(1)} \cap f_2^{(2)})$ | $m(f_2^{(1)} \cap f_2^{(2)})$ | ... | $m(f_i^{(1)} \cap f_2^{(2)})$ | ... | $m(f_k^{(1)} \cap f_2^{(2)})$ |
| ... | ... | ... | ... | ... | ... | ... |
| $m_2(f_j^{(2)})$ | $m(f_1^{(1)} \cap f_j^{(2)})$ | $m(f_2^{(1)} \cap f_j^{(2)})$ | ... | $m(f_i^{(1)} \cap f_j^{(2)})$ | ... | $m(f_k^{(1)} \cap f_j^{(2)})$ |
| ... | ... | ... | ... | ... | ... | ... |
| $m_2(f_l^{(2)})$ | $m(f_1^{(1)} \cap f_l^{(2)})$ | $m(f_2^{(1)} \cap f_l^{(2)})$ | ... | $m(f_i^{(1)} \cap f_l^{(2)})$ | ... | $m(f_k^{(1)} \cap f_l^{(2)})$ |

We look for the optimal rule of combination, i.e. the bpa $m(.) = m_1(.) \oplus m_2(.)$ which maximizes the joint entropy of the two information sources. The justification for the Maxent criteria is discussed in [24, 27]. Thus, one has to find $m(.)$ such that [66, 67].

$$\max_m [H(m)] \equiv \max_m \left[ -\sum_{i=1}^{k}\sum_{j=1}^{l} m(f_i^{(1)} \cap f_j^{(2)}) \log[m(f_i^{(1)} \cap f_j^{(2)})] \right] \equiv -\min_m [-H(m)] \quad (45)$$

satisfying both

- the measurement projection principle (marginal bpa), i.e. $\forall i = 1, \dots, k$ and $\forall j = 1, \dots, l$

$$m_1(f_i^{(1)}) = \sum_{j=1}^{l} m(f_i^{(1)} \cap f_j^{(2)}) \quad \text{and} \quad m_2(f_j^{(2)}) = \sum_{i=1}^{k} m(f_i^{(1)} \cap f_j^{(2)}) \quad (46)$$

These constraints state that the marginal bpa $m_1(.)$ is obtained by the summation over each column of the fusion table and the marginal bpa $m_2(.)$ is obtained by the summation over each row of the fusion table.

- the measurement balance principle (the sum of all cells of the fusion table must be unity)

$$\sum_{i=1}^{k}\sum_{j=1}^{l} m(f_i^{(1)} \cap f_j^{(2)}) = 1 \quad (47)$$

Using the concise notation $m_{ij} \triangleq m(f_i^{(1)} \cap f_j^{(2)})$, the Lagrangian associated with this optimization problem under equality constraints is given by (we consider here the minimization of $-J(m)$ appearing in r.h.s of (45))

$$\mathcal{L}(m, \lambda) = \sum_{i=1}^{k}\sum_{j=1}^{l} m_{ij} \ln[m_{ij}]$$

$$+ \sum_{i=1}^{k} \lambda_i [m_1(f_i^{(1)}) - \sum_{j=1}^{l} m_{ij}] + \sum_{j=1}^{l} \gamma_i [m_2(f_j^{(2)}) - \sum_{i=1}^{k} m_{ij}] + \eta [\sum_{i=1}^{k}\sum_{j=1}^{l} m_{ij} - 1] \quad (48)$$

which can be written more concisely as

$$\mathcal{L}(m,\lambda) = -H(m) + \lambda' \mathbf{g}(m) \tag{49}$$

where $m = [m_{11}\ m_{12}\ \ldots\ m_{kl}]'$ and

$$\lambda = \begin{bmatrix} \lambda_1 \\ \vdots \\ \lambda_k \\ \gamma_1 \\ \vdots \\ \gamma_l \\ \eta \end{bmatrix} \quad \text{and} \quad \mathbf{g}(m) = \begin{bmatrix} m_1(f_1^{(1)}) - \sum_{j=1}^{l} m_{1j} \\ \vdots \\ m_1(f_k^{(1)}) - \sum_{j=1}^{l} m_{kj} \\ m_2(f_1^{(2)}) - \sum_{i=1}^{k} m_{i1} \\ \vdots \\ m_2(f_l^{(2)}) - \sum_{i=1}^{k} m_{il} \\ \sum_{i=1}^{k} \sum_{j=1}^{l} m_{ij} - 1 \end{bmatrix} \tag{50}$$

Following the classical method of Lagrange multipliers, one has to find optimal solution $(m^*, \lambda^*)$ such that

$$\frac{\partial \mathcal{L}}{\partial m}(m^*, \lambda^*) = \mathbf{0} \quad \text{and} \quad \frac{\partial \mathcal{L}}{\partial \lambda}(m^*, \lambda^*) = \mathbf{0} \tag{51}$$

The first $k \times l$ equations express the general solution $m[\lambda]$ and the $k+l+1$ last equations determine $\lambda^*$ and therefore by substitution into $m[\lambda]$, the optimal solution $m^* = m[\lambda^*]$. One has to solve

$$\frac{\partial \mathcal{L}}{\partial m} = \begin{bmatrix} \frac{\partial \mathcal{L}}{\partial m_{11}} \\ \vdots \\ \frac{\partial \mathcal{L}}{\partial m_{ij}} \\ \vdots \\ \frac{\partial \mathcal{L}}{\partial m_{kl}} \end{bmatrix} = \begin{bmatrix} \ln(m_{11}) + 1 + \eta - \lambda_1 - \gamma_1 \\ \vdots \\ \ln(m_{ij}) + 1 + \eta - \lambda_i - \gamma_j \\ \vdots \\ \ln(m_{kl}) + 1 + \eta - \lambda_k - \gamma_l \end{bmatrix} = \begin{bmatrix} 0 \\ \vdots \\ \vdots \\ \vdots \\ 0 \end{bmatrix} = \mathbf{0} \tag{52}$$

which yields $\forall i, j$,

$$m_{ij} = e^{-\eta - 1} e^{\lambda_i} e^{\gamma_j} \tag{53}$$

and

$$\frac{\partial \mathcal{L}}{\partial \lambda} = \begin{bmatrix} \frac{\partial \mathcal{L}}{\partial \lambda_1} \\ \vdots \\ \frac{\partial \mathcal{L}}{\partial \lambda_k} \\ \frac{\partial \mathcal{L}}{\partial \gamma_1} \\ \vdots \\ \frac{\partial \mathcal{L}}{\partial \gamma_l} \\ \frac{\partial \mathcal{L}}{\partial \eta} \end{bmatrix} = \begin{bmatrix} 0 \\ \vdots \\ 0 \\ 0 \\ \vdots \\ 0 \\ 0 \end{bmatrix} = \mathbf{0} \Leftrightarrow \begin{bmatrix} e^{-\eta-1} \sum_{j=1}^{l} e^{\lambda_1} e^{\gamma_j} \\ \vdots \\ e^{-\eta-1} \sum_{j=1}^{l} e^{\lambda_k} e^{\gamma_j} \\ e^{-\eta-1} \sum_{i=1}^{k} e^{\gamma_1} e^{\lambda_i} \\ \vdots \\ e^{-\eta-1} \sum_{i=1}^{k} e^{\gamma_l} e^{\lambda_i} \\ e^{-\eta-1} \sum_{i=1}^{k} \sum_{j=1}^{l} e^{\gamma_l} e^{\lambda_i} \end{bmatrix} = \begin{bmatrix} m_1(f_1^{(1)}) \\ \vdots \\ m_1(f_k^{(1)}) \\ m_2(f_1^{(2)}) \\ \vdots \\ m_2(f_l^{(2)}) \\ 1 \end{bmatrix} \tag{54}$$

The last constraint in (54) can also be written as

$$e^{-\eta-1} \sum_{i=1}^{k} \sum_{j=1}^{l} e^{\gamma_l} e^{\lambda_i} = e^{-\eta-1} \left(\sum_{i=1}^{k} e^{\lambda_i}\right)\left(\sum_{j=1}^{l} e^{\gamma_l}\right) = 1 \tag{55}$$

Now with basic algebraic manipulation, the optimal global bpa $m_{ij}\ \forall i, j$ we are searching for, can be expressed as

$$m_{ij} = e^{-\eta-1} e^{\lambda_i} e^{\gamma_j}$$

$$= e^{-\eta-1} e^{\lambda_i} e^{\gamma_j} \times \overbrace{e^{-\eta-1} \left(\sum_{i=1}^{k} e^{\lambda_i}\right)\left(\sum_{j=1}^{l} e^{\gamma_l}\right)}^{1}$$

$$= \underbrace{\left(e^{-\eta-1} e^{\lambda_i} \sum_{j=1}^{l} e^{\gamma_l}\right)}_{m_1(f_i^{(1)})} \underbrace{\left(e^{-\eta-1} e^{\gamma_j} \sum_{i=1}^{k} e^{\lambda_i}\right)}_{m_2(f_j^{(2)})}$$

Thus, the solution of the maximisation of the joint entropy is obtained by choosing $\forall i, j$

$$m_{ij} = m(f_i^{(1)} \cap f_j^{(2)}) = m_1(f_i^{(1)}) m_2(f_j^{(2)}) \tag{56}$$

Since it may exist several combinations yielding to the same focal element, the bpa of all focal elements equal to $f_i^{(1)} \cap f_j^{(2)}$ over the fusion space is

$$m(f_i^{(1)} \cap f_j^{(2)}) = \sum_{i,j} m_1(f_i^{(1)}) m_2(f_j^{(2)}) \tag{57}$$

which coincides exactly with the new rule of combination expressed previously.

### 4.4.3 Numerical example of entropy calculation

We present here a very simple numerical example of the derivation of entropies of individual sources of informations and the combined (joint) entropy of combined sources. Let's consider the simple frame of discernment $\Theta = \{\theta_1, \theta_2\}$ and the two following (uncertain and paradoxical) information granules

$$m_1(\theta_1) = 0.60 \quad m_1(\theta_2) = 0.20 \quad m_1(\theta_1 \cup \theta_2) = 0.10 \quad m_1(\theta_1 \cap \theta_2) = 0.10$$

$$m_2(\theta_1) = 0.50 \quad m_2(\theta_2) = 0.20 \quad m_2(\theta_1 \cup \theta_2) = 0.10 \quad m_2(\theta_1 \cap \theta_2) = 0.20$$

The fusion rule can be described through the following fusion table

| $\oplus$ | $m_1(\theta_1) = 0.60$ | $m_1(\theta_2) = 0.20$ | $m_1(\theta_1 \cup \theta_2) = 0.10$ | $m_1(\theta_1 \cap \theta_2) = 0.10$ |
|---|---|---|---|---|
| $m_2(\theta_1) = 0.50$ | 0.30 $(\theta_1)$ | 0.10 $(\theta_1 \cap \theta_2)$ | 0.05 $(\theta_1)$ | 0.05 $(\theta_1 \cap \theta_2)$ |
| $m_2(\theta_2) = 0.20$ | 0.12 $(\theta_1 \cap \theta_2)$ | 0.04 $(\theta_2)$ | 0.02 $(\theta_2)$ | 0.02 $(\theta_1 \cap \theta_2)$ |
| $m_2(\theta_1 \cup \theta_2) = 0.10$ | 0.06 $(\theta_1)$ | 0.02 $(\theta_2)$ | 0.01 $(\theta_1 \cup \theta_2)$ | 0.01 $(\theta_1 \cap \theta_2)$ |
| $m_2(\theta_1 \cap \theta_2) = 0.20$ | 0.12 $(\theta_1 \cap \theta_2)$ | 0.04 $(\theta_1 \cap \theta_2)$ | 0.02 $(\theta_1 \cap \theta_2)$ | 0.02 $(\theta_1 \cap \theta_2)$ |

(58)

Each cell of the table provides a part of the global bpa $m(.)$ contribution for the corresponding proposition $M$ indicated between parentheses. The entropies of individual sources are given by

$$H(M_1) = -0.60 \ln(0.60) - 0.20 \ln(0.20) - 0.10 \ln(0.10) - 0.10 \ln(0.10) = 1.0889 \; nats$$
$$H(M_2) = -0.50 \ln(0.50) - 0.20 \ln(0.20) - 0.10 \ln(0.10) - 0.20 \ln(0.20) = 1.2206 \; nats$$

The conditional entropies $H(M_1|M_2)$ and $H(M_2|M_1)$ are given by [9]

$$\begin{aligned}
H(M_1|M_2) =& m_2(M_2 = \theta_1) H(M_1|M_2 = \theta_1) + m_2(M_2 = \theta_2) H(M_1|M_2 = \theta_2) \\
&+ m_2(M_2 = \theta_1 \cap \theta_2) H(M_1|M_2 = \theta_1 \cap \theta_2) + m_2(M_2 = \theta_1 \cup \theta_2) H(M_1|M_2 = \theta_1 \cup \theta_2) \\
=& 0.5 H\left[\left(\frac{0.30}{0.50}, \frac{0.10}{0.50}, \frac{0.05}{0.50}, \frac{0.05}{0.50}\right)\right] + 0.2 H\left[\left(\frac{0.12}{0.20}, \frac{0.04}{0.20}, \frac{0.02}{0.20}, \frac{0.02}{0.20}\right)\right] \\
&+ 0.1 H\left[\left(\frac{0.06}{0.10}, \frac{0.02}{0.10}, \frac{0.01}{0.10}, \frac{0.01}{0.10}\right)\right] + 0.2 H\left[\left(\frac{0.12}{0.20}, \frac{0.04}{0.20}, \frac{0.02}{0.20}, \frac{0.02}{0.20}\right)\right] \\
=& (0.5 \times 1.0889) + (0.2 \times 1.0889) + (0.1 \times 1.0889) + (0.2 \times 1.0889) \\
=& 1.0889 \; nats
\end{aligned}$$

$$\begin{aligned}
H(M_2|M_1) =& m_1(M_1 = \theta_1) H(M_2|M_1 = \theta_1) + m_1(M_1 = \theta_2) H(M_2|M_1 = \theta_2) \\
&+ m_1(M_1 = \theta_1 \cap \theta_2) H(M_2|M_1 = \theta_1 \cap \theta_2) + m_1(M_1 = \theta_1 \cup \theta_2) H(M_2|M_1 = \theta_1 \cup \theta_2) \\
=& 0.6 H\left[\left(\frac{0.30}{0.60}, \frac{0.12}{0.60}, \frac{0.06}{0.60}, \frac{0.12}{0.60}\right)\right] + 0.2 H\left[\left(\frac{0.10}{0.20}, \frac{0.04}{0.20}, \frac{0.02}{0.20}, \frac{0.04}{0.20}\right)\right] \\
&+ 0.1 H\left[\left(\frac{0.05}{0.10}, \frac{0.02}{0.10}, \frac{0.01}{0.10}, \frac{0.02}{0.10}\right)\right] + 0.1 H\left[\left(\frac{0.05}{0.10}, \frac{0.02}{0.10}, \frac{0.01}{0.10}, \frac{0.02}{0.10}\right)\right] \\
=& (0.6 \times 1.2206) + (0.2 \times 1.2206) + (0.1 \times 1.2206) + (0.1 \times 1.2206) \\
=& 1.2206 \; nats
\end{aligned}$$

Therefore, one has

$$H(M_1) = H(M_1|M_2) = 1.0889 \; nats \quad \text{and} \quad H(M_2) = H(M_2|M_1) = 1.2206 \; nats$$

The joint entropy $H(M) = H(M_1, M_2)$ is directly obtained from the cells of the fusion table and one gets

$$H(M) = H[(0.3, 0.1, 0.05, 0.05, 0.12, 0.04, 0.02, 0.02, 0.06, 0.02, 0.01, 0.01, 0.12, 0.04, 0.02, 0.02)] = 2.3095 \ nats$$

Hence, one has verified the classical result (chain rule) of the information theory, i.e.

$$H(M) = H(M_1) + H(M_2|M_1) = H(M_2) + H(M_1|M_2)$$

or more specially because of the independence of the two sources of information

$$H(M) = H(M_1) + H(M_2)$$

Note that $H(M)$ must be evaluated using the full description of the fusion table (from all the cells of the table) and not from the global bpa $m(.)$. Otherwise a smaller value for $H(M)$ is deduced, as it can be easily shown. From the fusion table, one gets the final bpa $m(.)$ with

$$m(\theta_1) = 0.41 \quad m(\theta_2) = 0.08 \quad m(\theta_1 \cup \theta_2) = 0.01 \quad m(\theta_1 \cap \theta_2) = 0.50$$

The evaluation of $H(M)$ from bpa $m(.)$ yields the value

$$\tilde{H}(M) = H[(0.41, 0.08, 0.01, 0.50)] = 0.96023 \ nats < H(M)$$

*Remark*

Note that in this example, the combination of the two sources reduces the uncertainty of judgment of each local information sources since $\tilde{H}(M) < H(M_1)$ and $\tilde{H}(M) < H(M_2)$. This is unfortunalely not a valid conclusion in general as many people (wrongly) think. We argue that the fusion of independent sources of information does not necessarily reduces the uncertainty of judgment. To convince the reader, just take the similar example with the following new information granules

$$m_1(\theta_1) = 0.900 \quad m_1(\theta_2) = 0.090 \quad m_1(\theta_1 \cup \theta_2) = 0.009 \quad m_1(\theta_1 \cap \theta_2) = 0.001$$
$$m_2(\theta_1) = 0.090 \quad m_2(\theta_2) = 0.900 \quad m_2(\theta_1 \cup \theta_2) = 0.009 \quad m_2(\theta_1 \cap \theta_2) = 0.001$$

It is not too difficult to check that global bpa $m(.)$ is

$$m(\theta_1) = 0.08991 \quad m(\theta_2) = 0.08991 \quad m(\theta_1 \cup \theta_2) = 0.000081 \quad m(\theta_1 \cap \theta_2) = 0.820099$$

with corresponding entropies

$$H(M_1) = H(M_2) = 0.36084 \ nats \quad \text{and} \quad (\tilde{H}(M) = 0.59659) < (H(M) = 0.72168)$$

but $\tilde{H}(M) > H(M_1)$ and $\tilde{H}(M) > H(M_2)$. Thus in this case, the fusion increases actually the uncertainty of the final judgment.

### 4.4.4 Definition for the generalized entropy of a source

The evaluation of the entropy $H(m)$ of a given source from the direct extension of its classical definition, with convention (see [9]) $0 \ln(0) = 0$ and with bpa $m(.)$, i.e.

$$H(m) = - \sum_{A \in D^\Theta} m(A) \ln(m(A))$$

seems to not be the best measure for the self-information of a general (uncertain and paradoxical) source of information because it does not catch the intrinsic informational strength (i.i.s. for short) $s(A)$ of the propositions $A$. An extension of the classical entropy in the DST framework had already been proposed in 1983 by R. Yager based on the weight of conflict between the belief function Bel and the certain support function $\text{Bel}_A$ focused on each proposition $A$ (see [70] for details). In the classical definition (based only on probability measure), one always has $s(A) \equiv |A| = 1$. This does not hold in our general theory of plausible and paradoxical reasoning and we propose to generalize the notion of entropy in the following manner to measure correctly the self-information of a general source :

$$H_g(m) = - \sum_{A \in D^\Theta} \frac{1}{s(A)} m(A) \ln(\frac{1}{s(A)} m(A)) \tag{59}$$

$H_g(m)$ will be called the *generalized entropy* of the source associated with bpa $m(.)$. This general definition introduces the cardinality of a general (irreducible) proposition $A$ which can be derived from the two following important rules

$$s\left(\bigcup_{i=1,n} B_i\right) = s\left(B_1 \cup \ldots \cup B_n\right) = \frac{\sum_{i=1,n} 1/s\left(B_i\right)}{\prod_{i=1,n} 1/s\left(B_i\right)} \tag{60}$$

$$s\left(\bigcap_{i=1,n} B_i\right) = s\left(B_1 \cap \ldots \cap B_n\right) = \frac{\prod_{i=1,n} s\left(B_i\right)}{\sum_{i=1,n} s\left(B_i\right)} \tag{61}$$

It is very important to note that these rules apply only on irreductible propositions (logical atoms) $A$. A proposition $A$ is said to be irreductible (or equivalently has a compact form) if and only if it does not admit other equivalent form with a smaller number of operands and operators. For example $(\theta_1 \cup \theta_3) \cap (\theta_2 \cup \theta_3)$ is not an irreductible proposition since it can be reduced to its equivalent logical atom $(\theta_1 \cap \theta_2) \cup \theta_3$. To compute the i.i.s. $s(A)$ of any proposition $A$ using the rules (60) and (61), the proposition has first to be reduced to its minimal representation (irreductible form).

*Examples*

Here are few examples of the value of the cardinality for some elementary and composite irreductible propositions $A$. We recall that $\theta_i$ involved in $A$ are singletons such that $|\theta_i| = 1$.

$$A = \theta_1 \cup \theta_2 \Rightarrow s(A) = 2$$
$$A = \theta_1 \cap \theta_2 \Rightarrow s(A) = 1/2$$
$$A = \theta_1 \cup \theta_2 \cup \theta_3 = (\theta_1 \cup \theta_2) \cup \theta_3 = \theta_1 \cup (\theta_2 \cup \theta_3) = \theta_2 \cup (\theta_1 \cup \theta_3) \Rightarrow s(A) = 3$$
$$A = \theta_1 \cap \theta_2 \cap \theta_3 = (\theta_1 \cap \theta_2) \cap \theta_3 = \theta_1 \cap (\theta_2 \cap \theta_3) = \theta_2 \cap (\theta_1 \cap \theta_3) \Rightarrow s(A) = 1/3$$
$$A = (\theta_1 \cap \theta_2) \cup \theta_3 \Rightarrow s(A) = 3/2$$
$$A = (\theta_1 \cup \theta_2) \cap \theta_3 \Rightarrow s(A) = 2/3$$
$$A = (\theta_1 \cap \theta_2) \cup (\theta_3 \cap \theta_4) \Rightarrow s(A) = 1$$
$$A = (\theta_1 \cup \theta_2) \cap (\theta_3 \cup \theta_4) \Rightarrow s(A) = 1$$
$$A = (\theta_1 \cap \theta_2) \cup (\theta_3 \cap \theta_4 \cap \theta_5) \Rightarrow s(A) = 5/6$$
$$A = (\theta_1 \cup \theta_2) \cap (\theta_3 \cup \theta_4 \cup \theta_5) \Rightarrow s(A) = 6/5$$

Thus the evaluation of $s(A)$ for any general irreductible proposition $A$ can always be obtained from the two basic rules (60) and (61). This generalized definition makes sense with the notion of entropy and is coherent with classical definition (i.e. $H_g(m) \equiv H(m)$ when $m(.)$ becomes a bayesian bpa $p(.)$). Let $\Theta = \{\theta_1, \ldots, \theta_n\}$ be a *general* frame of discernment of the problem under consideration and a general body of evidence with information granule $m(.)$ on $D^\Theta$, then the generalized entropy $H_g(m)$ takes its minimal value $-n \ln(n)$ when the source provides the maximum of paradoxe which is obtained when $m(\theta_1 \cap \ldots \cap \theta_n) = 1$. It is important to note that the maximum of uncertainty is not obtained when $m(\theta_1 \cup \ldots \cup \theta_n) = 1$ but rather for a specific $m()$ which distributes some weight of evidence assignment to each proposition $A \in D^\Theta$ because there is less information (from the information theory viewpoint) when there exists several propositions with non nul bpa rather than one. One has also to take into account the intrinsic self-information of the propositions to get a good measure of global information provided by a source. The generalized entropy includes both aspects of the information (the intrinsic and the classical aspect). The uniform distribution for $m(.)$ does not generate the maximum generalized-entropy because of the different intrinsic self-information of each proposition (see next example). We argue that the generalized entropy of any source defined with respect to a frame $\Theta$ appears to be a very useful tool to measure the degree of uncertainty and paradox of any given source of information.

*Example*

We give here some values of $H_g(m)$ for different sources of information over the same frame $\Theta = \{\theta_1, \theta_2\}$. The sources have been classified from the most informative one $\mathcal{B}_1$ up to the less informative one $\mathcal{B}_{16}$. $\mathcal{B}_{16}$ corresponds to the source containing minimal information on the hyper-power set of the frame $\Theta$ (thus $\mathcal{B}_{16}$ has the minimal discrimination power between all possible propositions). There does not exist a source $\mathcal{B}_k$ such that $H_g^{\mathcal{B}_k}(m) > H_g^{\mathcal{B}_{16}}(m)$ for this simpliest example. Finding $m^*(.)$ such that $H_g(m^*)$ takes its maximal value for a general frame $\Theta$ with $|\Theta| = n$ is called the general whitening source problem. No solution for this problem has been obtained so far.

| | $m(\theta_1)$ | $m(\theta_2)$ | $m(\theta_1 \cup \theta_2)$ | $m(\theta_1 \cap \theta_2)$ | $H_g(m)$ |
|---|---|---|---|---|---|
| $\mathcal{B}_1$ | 0 | 0 | 0 | 1 | $-1.386$ |
| $\mathcal{B}_2$ | 0 | 0 | 0.3 | 0.7 | $-0.186$ |
| $\mathcal{B}_3$ | 1 | 0 | 0 | 0 | 0 |
| $\mathcal{B}_4$ | 0 | 1 | 0 | 0 | 0 |
| $\mathcal{B}_5$ | 0.1 | 0.2 | 0 | 0.7 | 0.081 |
| $\mathcal{B}_6$ | 0 | 0 | 1 | 0 | 0.346 |
| $\mathcal{B}_7$ | 0.8 | 0.2 | 0 | 0 | 0.500 |
| $\mathcal{B}_8$ | 0 | 0 | 0.7 | 0.3 | 0.673 |
| $\mathcal{B}_9$ | 0.5 | 0.5 | 0 | 0 | 0.693 |
| $\mathcal{B}_{10}$ | 0.7 | 0.2 | 0.1 | 0 | 0.721 |
| $\mathcal{B}_{11}$ | 0.7 | 0.2 | 0 | 0.1 | 0.893 |
| $\mathcal{B}_{12}$ | 0.1 | 0.2 | 0.7 | 1 | 0.919 |
| $\mathcal{B}_{13}$ | 0.1 | 0.2 | 0.3 | 0.4 | 1.015 |
| $\mathcal{B}_{14}$ | 0.1 | 0.2 | 0.4 | 0.3 | 1.180 |
| $\mathcal{B}_{15}$ | 0.25 | 0.25 | 0.25 | 0.25 | 1.299 |
| $\mathcal{B}_{16}$ | 0.25 | 0.25 | 0.35 | 0.15 | 1.359 |

$\mathcal{B}_1$ is the most informative source because all the weights of evidence about the truth are focused only on the smaller element $\theta_1 \cap \theta_2$ of hyper-powerset $D^\Theta$. $\mathcal{B}_2$ is less informative than $\mathcal{B}_1$ because there exists an ambiguity between the two propositions $\theta_1 \cup \theta_2$ and $\theta_1 \cap \theta_2$. $\mathcal{B}_3$ and $\mathcal{B}_4$ are less informative than $\mathcal{B}_1$ because the weights of evidence about the truth are focused on larger elements ($\theta_1$ or $\theta_2$ respectively) of $D^\Theta$. $\mathcal{B}_6$ is less informative than $\mathcal{B}_3$ or $\mathcal{B}_4$ because the weight of evidence about the truth is focused on a bigger element $\theta_1 \cup \theta_2$ of $D^\Theta$. $\mathcal{B}_7$ is less informative than previous sources since there exists an ambiguity between the two propositions $\theta_1$ and $\theta_2$ but it is more informative than $\mathcal{B}_9$ since the discrimination power (our easiness to decide which proposition supports the truth) is higher with $\mathcal{B}_7$ than with $\mathcal{B}_9$. Note that even if in this very simple example, it is not obvious to see that $\mathcal{B}_{16}$ is the less informative (white) source of information. Most of readers would have probably thought to choose either $\mathcal{B}_6$ or $\mathcal{B}_{15}$. This comes from the confusion between the intrinsic information supported by the proposition itself and the information supported by the whole bpa $m(.)$.

### 4.4.5 Zadeh's example

Let's take back the disturbing Zadeh's example given in section 3.4. Two doctors examine a patient and agree that it suffers from either meningitis (M), concussion (C) or brain tumor (T). Thus $\Theta = \{M, C, T\}$. Assume that the doctors agree in their low expectation of a tumor, but disagree in likely cause and provide the following diagnosis

$$m_1(M) = 0.99 \qquad m_1(T) = 0.01 \qquad \text{and} \quad \forall A \in D^\Theta, A \neq T, A \neq M, \quad m_1(A) = 0$$

$$m_2(C) = 0.99 \qquad m_2(T) = 0.01 \qquad \text{and} \quad \forall A \in D^\Theta, A \neq T, A \neq C, \quad m_2(A) = 0$$

The new general rule of combination (41), yields the following combined information granule

$$m(M \cap C) = 0.9801 \qquad m(M \cap T) = 0.0099 \qquad m(C \cap T) = 0.0099 \qquad m(T) = 0.0001$$

From this granule, one gets

$$\text{Bel}(M) = m(M \cap C) + m(M \cap T) = 0.99$$
$$\text{Bel}(C) = m(M \cap C) + m(T \cap C) = 0.99$$
$$\text{Bel}(T) = m(T) + m(M \cap T) + m(C \cap T) = 0.0199$$

If both doctors can be considered as equally reliable, the combined information granule $m(.)$ mainly focuses weight of evidence on the paradoxical proposition $M \cap C$ which means that patient suffers both meningitis and concussion but almost surely not from brain tumor. This conclusion is coherent with the common sense actually. Then, no therapy for brain tumor (like heavy and ever risky brain surgical intervention) will be chosen in such case. This really helps to take important decision to save the life of the patient in this example. A deeper medical examination adapted to both meningitis and concussion will almost surely be done before applying the best therapy for the patient. Just remember that in this case, the DST had concluded that the patient had brain tumor with certainty . . . .

### 4.4.6 Mahler's example revisited

Let's consider now the following example excerpt from the R. Mahler's paper [36]. We consider that our classification knowledge base consists of the three (imaginary) new and rare diseases corresponding to following frame of discernment

$$\Theta = \{\theta_1 = kotosis, \theta_2 = phlegaria, \theta_3 = pinpox\}$$

We assume that the three diseases are equally likely to occur in the patient population but there is some evidence that *phlegaria* and *pinpox* are the same disease and there is also a small possibility that *kotosis* and *phlegaria* might be the same disease. Finally, there is a small possibility that all three diseases are the same. This information can be expressed by assigning a priori bpa as follows

$$m_0(\theta_1) = 0.2 \qquad m_0(\theta_2) = 0.2 \qquad m_0(\theta_3) = 0.2$$
$$m_0(\theta_2 \cap \theta_3) = 0.2 \quad m_0(\theta_1 \cap \theta_2) = 0.1 \quad m_0(\theta_1 \cap \theta_2 \cap \theta_3) = 0.1$$

Let $Bel(.)$ the prior belief measure corresponding to this prior bpa $m(.)$. Now assume that Doctor $D_1$ and Doctor $D_2$ examine a patient and deliver diagnoses with following reports:

- Report for $D_1$:     $m_1(\theta_1 \cup \theta_2 \cup \theta_3) = 0.05 \qquad m_1(\theta_2 \cup \theta_3) = 0.95$

- Report for $D_2$:     $m_2(\theta_1 \cup \theta_2 \cup \theta_3) = 0.20 \qquad m_2(\theta_2) = 0.80$

The combination of the evidences provided by the two doctors $m' = m_1 \oplus m_2$ obtained by the general rule of combination (41) yields the following bpa $m'(.)$

$$m'(\theta_2) = 0.8 \qquad m'(\theta_2 \cup \theta_3) = 0.19 \qquad m'(\theta_1 \cup \theta_2 \cup \theta_3) = 0.01$$

The combination of bpa $m'(.)$ with prior evidence $m_0(.)$ yields the final bpa $m = m_0 \oplus m' = m_0 \oplus [m_1 \oplus m_2]$ with

$$m(\theta_1) = 0.002 \qquad m(\theta_2) = 0.200 \qquad m(\theta_3) = 0.040$$
$$m(\theta_1 \cap \theta_2) = 0.260 \qquad m(\theta_2 \cap \theta_3) = 0.360 \quad m(\theta_1 \cap \theta_2 \cap \theta_3) = 0.100$$
$$m(\theta_1 \cap (\theta_2 \cup \theta_3)) = 0.038$$

Therefore the final belief function given by (37) is

$$Bel(\theta_1) = 0.002 + 0.260 + 0.100 + 0.038 = 0.400$$
$$Bel(\theta_2) = 0.200 + 0.260 + 0.360 + 0.100 = 0.920$$
$$Bel(\theta_3) = 0.040 + 0.360 + 0.100 = 0.500$$
$$Bel(\theta_1 \cap \theta_2) = 0.260 + 0.100 = 0.360$$
$$Bel(\theta_2 \cap \theta_3) = 0.360 + 0.100 = 0.460$$
$$Bel(\theta_1 \cap (\theta_2 \cup \theta_3)) = 0.038 + 0.100 = 0.138$$
$$Bel(\theta_1 \cap \theta_2 \cap \theta_3) = 0.100$$

Thus, on the basis of all the evidences one has, we are able to conclude with high a degree of belief that the patient has phlegaria which is coherent with the Mahler's conclusion based on his Conditioned Dempster-Shafer theory developed from his conditional event algebra although a totally new and simpliest approach has been adopted here.

### 4.4.7 A thief identification example

Let's revisit a very simple thief identification example. Assume that a 75 years old grandfather is taking a walk with his 9 years old grandson in a park. They saw at 50 meters away, a 45 years old pickpocket robbing the bag of an old lady. A policeman looking for some witnesses of this event asks separately the grandfather and his grandchild if they have seen the thief (they both answer yes) and how was the thief (a young or an old man). The grandfather (source of information $B_1$ reports that the thief was a young man with high confidence 0.99 and with only a low uncertainty 0.01. His grandson reports that the thief was a old man with high confidence 0.99 and with only a low uncertainty 0.01. These two witnesses provide fair reports (with respect to their own world of knowledge) even if apparently they appear as paradoxical. The policeman then send the two reports with only minimal information about witnesses (saying only their names and that they were a priori fully trustable) to an investigator. The investigator has no possibility to meet or to call back the witnesses in order to get more details.

Under such condition, what would be the best decision to be taken by the investigator about the age of the thief to eventually help to catch him? Such kind of simple examples occur quite frequently in witnesses problems actually. A rational investigator will almost surely suspect a mistake or an error in one or both reports since they appear apparently in full contradiction. The investigator will then try to take his final decision with some other better information (if any).

If the investigator uses our new plausible and paradoxical reasoning, he will defined the following bpa with respect to the frame of discernment $\Theta = \{\theta_1 = \text{young}, \theta_2 = \text{old}\}$ and the available reports $\mathcal{B}_1$ and $\mathcal{B}_2$ with following bpa

$$m_1(\theta_1) = 0.99 \quad m_1(\theta_2) = 0 \quad m_1(\theta_1 \cup \theta_2) = 0.01 \quad m_1(\theta_1 \cap \theta_2) = 0$$

$$m_2(\theta_1) = 0 \quad m_2(\theta_2) = 0.99 \quad m_2(\theta_1 \cup \theta_2) = 0.01 \quad m_2(\theta_1 \cap \theta_2) = 0$$

The fusion of these two sources of information yields the global bpa $m(.)$ with

$$m(\theta_1) = 0.0099 \quad m(\theta_2) = 0.0099 \quad m(\theta_1 \cup \theta_2) = 0.0001 \quad m(\theta_1 \cap \theta_2) = 0.9801$$

Thus, from this global information, the investigator has no better choice but to consider with almost certainty that the thief was both a young and old man. By assuming that the expected life duration is around 80 years, the inspector will deduce that the true age of the thief is around 40 years old which is not too far from the truth. At least, this conclusion could be helpful to interrogate some suspicious individuals.

### 4.4.8 A model to generate information granules $m(.)$ from intervals

We present here a model to generate information granules $m(.)$ from information represented by intervals. It is very common in practice that uncertain sources of information provide evidence on a given proposition in term of basic intervals $[\epsilon_*, \epsilon^*] \subset [0, 1]$ rather than a direct bpa $m(.)$. In such cases, some preprocessing must be done before applying the general rule of combination between such sources to take the final decision.

In the DST framework, we recall that the simplest and easiest transformation to convert $[\epsilon_*, \epsilon^*]$ into bpa has already been proposed by A. Appriou in [3]. The basic idea was to interpret $\epsilon_*$ as the minimal credibility committed to $A$ and $\epsilon^*$ as the plausibility committed to $A$. In other words, the Appriou's transformation model within the DST is the following one

$$\epsilon_* = m(A)$$
$$\epsilon^* = 1 - m(A^c)$$
$$\epsilon^* - \epsilon_* = m(A \cup A^c)$$

This model can be directly extended within our new theory of plausible and paradoxical reasoning by setting now.

$$\epsilon_* = m(A) + \frac{1}{2}m(A \cap A^c)$$
$$\epsilon^* = 1 - m(A^c) - \frac{1}{2}m(A \cap A^c)$$
$$\epsilon^* - \epsilon_* = m(A \cup A^c)$$

or equivalently

$$m(A) + \frac{1}{2}m(A \cap A^c) = \epsilon_* \tag{62}$$

$$m(A^c) + \frac{1}{2}m(A \cap A^c) = 1 - \epsilon^* \tag{63}$$

$$m(A \cup A^c) = \epsilon^* - \epsilon_* \tag{64}$$

This appealing model presents nice properties specially when $\epsilon^* = \epsilon_* = 0$ or when $\epsilon^* = \epsilon_* = 1$. This model is moreover coherent with the previous Appriou's model whenever the source becomes rational (i.e $m(A \cap A^c) = 0$). This new model presents however a degree of freedom since one has only two constraints (62) and (63) for three unknowns $m(A)$, $m(A^c)$ and $m(A \cap A^c)$. Thus in general, without an additional constraint, there exists many possible choices for $m(A)$, $m(A^c)$ and $m(A \cap A^c)$ and therefore there exists several bpa $m(.)$ satisfying this transformation model. Without extra prior information, it becomes difficult to justify the choice of a specific bpa versus all other admissible possibilities for $m(.)$.

To solve this important drawback, we propose to add the constraint on the maximization of the generalized-entropy $H_g(m)$. This will allow us to obtain from $[\epsilon_*, \epsilon^*]$ the unique bpa $m(.)$ having the minimum of specificity and admissible with our transformation model. From definition of $H_g(m)$ and previous equations (62)-(64), one gets

$$H_g(m) = -(\epsilon_* - m(A \cap A^c)/2)\ln(\epsilon_* - m(A \cap A^c)/2) - (1 - \epsilon^* - m(A \cap A^c)/2)\ln(1 - \epsilon^* - m(A \cap A^c)/2)$$
$$- \frac{1}{2}(\epsilon^* - \epsilon_*)\ln(\frac{1}{2}(\epsilon^* - \epsilon_*)) - 2m(A \cap A^c)\ln(2m(A \cap A^c))$$

The maximization of $H_g(m)$ is obtained for the optimal value $m^\star(A \cap A^c)$ such that $\frac{\partial H_g}{\partial m(A \cap A^c)}(m^\star(A \cap A^c)) = 0$ and $\frac{\partial^2 H_g}{\partial m(A \cap A^c)^2}(m^\star(A \cap A^c)) < 0$. The annulation of the first derivative is obtained by the solution of the equation

$$\frac{1}{2}\ln(\epsilon_* - m^\star/2) + \frac{1}{2}\ln(1 - \epsilon^* - m^\star/2) - 2m^\star \ln(2m^\star) - 1 = 0$$

or equivalently after basic algebraic manipulations

$$64e^2(m^\star)^4 - (m^\star)^2 + 2(1 - \epsilon^* + \epsilon_*)m^\star - 4(1 - \epsilon^*)\epsilon_* = 0 \qquad (65)$$

The solution of this equation does not admit a simple analytic expression but can be easily found using classical numerical methods. It is also easy to check that the second derivative is always negative and therefore $H_g(m)$ reaches its maximal value when

$$m(A) + \frac{1}{2}m^\star(A \cap A^c) = \epsilon_* \qquad (66)$$

$$m(A^c) + \frac{1}{2}m^\star(A \cap A^c) = 1 - \epsilon^* \qquad (67)$$

$$m(A \cup A^c) = \epsilon^* - \epsilon_* \qquad (68)$$

This completes the definition of our new transformation model. Note that $[\epsilon_*, \epsilon^*]$ can also be generated from bpa $m(.)$ through (62)-(64).

*Numerical examples*

- $[\epsilon_*, \epsilon^*] = [0.0, 0.0]$    $m(A \cap A^c) = 0.000$    $m(A) = 0.000$    $m(A^c) = 1.000$    $m(A \cup A^c) = 0.000$
- $[\epsilon_*, \epsilon^*] = [0.2, 0.2]$    $m(A \cap A^c) \approx 0.164$    $m(A) \approx 0.118$    $m(A^c) \approx 0.718$    $m(A \cup A^c) = 0.000$
- $[\epsilon_*, \epsilon^*] = [0.5, 0.5]$    $m(A \cap A^c) \approx 0.192$    $m(A) \approx 0.404$    $m(A^c) \approx 0.404$    $m(A \cup A^c) = 0.000$
- $[\epsilon_*, \epsilon^*] = [0.8, 0.8]$    $m(A \cap A^c) \approx 0.164$    $m(A) \approx 0.718$    $m(A^c) \approx 0.118$    $m(A \cup A^c) = 0.000$
- $[\epsilon_*, \epsilon^*] = [1.0, 1.0]$    $m(A \cap A^c) = 0.000$    $m(A) = 1.000$    $m(A^c) = 0.000$    $m(A \cup A^c) = 0.000$
- $[\epsilon_*, \epsilon^*] = [0.2, 0.4]$    $m(A \cap A^c) \approx 0.152$    $m(A) \approx 0.124$    $m(A^c) \approx 0.524$    $m(A \cup A^c) = 0.200$
- $[\epsilon_*, \epsilon^*] = [0.6, 0.8]$    $m(A \cap A^c) \approx 0.152$    $m(A) \approx 0.524$    $m(A^c) \approx 0.124$    $m(A \cup A^c) = 0.200$
- $[\epsilon_*, \epsilon^*] = [0.4, 0.6]$    $m(A \cap A^c) \approx 0.170$    $m(A) \approx 0.315$    $m(A^c) \approx 0.315$    $m(A \cup A^c) = 0.200$
- $[\epsilon_*, \epsilon^*] = [0.3, 0.9]$    $m(A \cap A^c) \approx 0.100$    $m(A) \approx 0.250$    $m(A^c) \approx 0.050$    $m(A \cup A^c) = 0.600$
- $[\epsilon_*, \epsilon^*] = [0.0, 1.0]$    $m(A \cap A^c) = 0.000$    $m(A) = 0.000$    $m(A^c) = 0.000$    $m(A \cup A^c) = 1.000$

## 5 Plausible and paradoxical reasoning in the neutrosophy framework

### 5.1 Neutrosophy and the neutrosophic logic

The neutrosophy is a new branch of philosophy, introduced by Florentin Smarandache in 1980, which studies the origin, nature and scope of neutralities, as well as their interactions with different ideational spectra. Neutrosophy considers a proposition, theory, event, concept, or entity, $A$ in relation to its opposite, *anti-A* and that which is not $A$, *non-A*, and that which is neither $A$ nor *anti-A*, denoted by *neut-A*. Neutrosophy serves as the basis for the neutrosophic logic [22].

The Neutrosophic Logic (NL) or Smarandache's logic is a general framework for the unification of all existing logics [56, 57, 58]). The main idea of NL is to characterize each logical statement in a 3D neutrosophic space where each dimension of the space represents respectively the truth (T), the falsehood (F) and the indeterminacy (I) of the statement under consideration where T, I and F are standard or non-sandard real subsets of $]^-0; 1^+[$. Moreover in NL, each statement is allowed to be over or under true, over or under false and over or under inderterminate by using hyper real numbers developed in the non-standard analysis theory [43, 14]. The neutrosophical value $\mathfrak{N}(A) = (T(A), I(A), F(A))$ in a frame of discernment (world of discourse) $\Theta$ of a statement A is then defined as a subset (a volume not necessary connexe; i.e. a set of disjoint volumes) of the neutrosophic space. Any statement $A$ represented by a triplet $\mathfrak{N}(A)$ is called a *neutrosophic event* or $\mathfrak{N}-$ event. The subset $\mathfrak{N}_t \triangleq T(A)$ characterizes the truth part of statement A. $\mathfrak{N}_i \triangleq$ (A) and $\mathfrak{N}_f \triangleq F(A)$ represent the inderterminacy and the falsehood of A. This Smarandache's representation is close to the human reasoning. It characterizes and catches the imprecision of knowledge or linguistic inexactitude received by various observers, uncertainty due to incomplete knowledge of acquisition errors or stochasticity, and vagueness due to lack of clear contours or limits. This approach allows theoretically to consider any kinds of logical statements. For example, the fuzzy set logic or the classical modal logic (which works with statements verifying $T(A)$, $I(A) \equiv 0$, $F(A) = 1 - T(A)$, where $T$ is a real number belonging to $[0;1]$) are included in NL. The neutrosophic logic can easily handle also paradoxes. We emphaze the fact that in general the neutrosohic value $\mathfrak{N}(A)$ of a proposition $A$ can also depend on dynamical parameters which can evolve with time, space, etc. For seak of concise notation, we omit to introduce this dependence in our notations in the sequel.

**Basic operations on sets**

Beside this modelling, F. Smarandache has introduced the following operations on sets. Consider $S_1$ and $S_2$ be two (unidimensional) standard or non-standard real subsets. The addition, substraction, multiplication and division (by a non null finite number) of these sets are defined as follows :

- **Addition**

$$S_1 \oplus S_2 = S_2 \oplus S_1 \triangleq \{x \mid x = s_1 + s_2, \forall s_1 \in S_1, \forall s_2 \in S_2\} \tag{69}$$

- **Substraction**

$$S_1 \ominus S_2 = -(S_2 \ominus S_1) \triangleq \{x \mid x = s_1 - s_2, \forall s_1 \in S_1, \forall s_2 \in S_2\} \tag{70}$$

  For *real positive* subsets, the Inf and Sup values of $S_1 \ominus S_2$ are given by

$$\text{Inf}[S_1 \ominus S_2] = \text{Inf}[S_1] - \text{Sup}[S_2] \quad \text{and} \quad \text{Sup}[S_1 \ominus S_2] = \text{Sup}[S_1] - \text{Inf}[S_2]$$

- **Multiplication**

$$S_1 \odot S_2 = S_2 \odot S_1 \triangleq \{x \mid x = s_1 \cdot s_2, \forall s_1 \in S_1, \forall s_2 \in S_2\} \tag{71}$$

  For *real positive* subsets, one gets

$$\text{Inf}[S_1 \odot S_2] = \text{Inf}[S_1] \cdot \text{Inf}[S_2] \quad \text{and} \quad \text{Sup}[S_1 \odot S_2] = \text{Sup}[S_1] \cdot \text{Sup}[S_2]$$

- **Division of a set by a non null standard number**
  Let $k \in \mathbb{R}^*$, then

$$S_1 \oslash k \triangleq \{x \mid x = s_1/k, \forall s_1 \in S_1\} \tag{72}$$

**Neutrosophic topology**

Let's construct now a neutrosophic topology (NT) [56] on interval $]^-0; 1^+[$, by considering the associated family of standard or non-standard subsets included in $]^-0; 1^+[$ and the empty set $\emptyset$, which is closed under set union and finite intersection. The union and intersection of two any propositions $A$ and $B$ (corresponding to either the part of truth, indeterminacy or falsehood of a given assertion defined on $]^-0; 1^+[$) are defined as follows

$$A \cup B = (A \oplus B) \ominus (A \odot B) \quad \text{and} \quad A \cap B = A \odot B \tag{73}$$

The neutrosophic complement of $A$ is defined as $\bar{A} = \{1^+\} \ominus A$ and the $\mathfrak{N}$ − value of an assertion A is characterized by a mapping function $\mathfrak{N}(.)$ such that

$$\mathfrak{N} : A \quad \mapsto \quad \mathfrak{N}(A) = (T(A), I(A), F(A)) \subset ]^-0; 1^+[^3 \tag{74}$$

The interval $]^-0; 1^+[$, endowed with this topology, forms a *neutrosophic topological space.*

Consider now two statements $A_1$ and $A_2$, then one defines the following basic neutrosophic operators:

$$\mathfrak{N}(A_1) \boxplus \mathfrak{N}(A_2) = (T_1 \oplus T_2, I_1 \oplus I_2, F_1 \oplus F_2) \tag{75}$$

$$\mathfrak{N}(A_1) \boxminus \mathfrak{N}(A_2) = (T_1 \ominus T_2, I_1 \ominus I_2, F_1 \ominus F_2) \tag{76}$$

$$\mathfrak{N}(A_1) \boxdot \mathfrak{N}(A_2) = (T_1 \odot T_2, I_1 \odot I_2, F_1 \odot F_2) \tag{77}$$

where $T_i = T(A_i)$, $I_i = I(A_i)$, $F_i = F(A_i)$ for $i = 1, 2$.

Since the truth, falsehood and indeterminacy of any statement must belong to $]^-0; 1^+[$, the result of each previous operator $\boxplus$, $\boxminus$ and $\boxdot$ must be in $]^-0; 1^+[^3$. Therefore upper and lower bounds of $T_1 \oplus T_2$ must be set respectively to $^-0$ and $1^+$ whenever $inf(T_1 \oplus T_2) < 0$ or $sup(T_1 \oplus T_2) > 1$. The same remark applies for $\boxminus$ and $\boxdot$ operators and for falsehood and inderterminacy part of compounded statement.

All classical logical operators and connectors can be extented in the $\mathfrak{N}$– Logic. For notation convenience, we will identify logical operators with their classical counterpart in set theory as pointed out in [35] (hence the following equivalences will be used $\neg A \equiv \bar{A}$, $A_1 \wedge A_2 \equiv A_1 \cap A_2$ and $A_1 \vee A_2 \equiv A_1 \cup A_2$ throughout this paper). We recall here only important operators used in the sequel. Additional neutrosophic logical operators like (strong disjunction, implication, equivalence, Sheffer's and Pierce's connectors) and general, physics and philosophical examples of application of neutrosophic operators can be found in [56, 58].

- **Negation**

$$\mathfrak{N}(\bar{A}) = (\{1\} \ominus T(A), \{1\} \ominus I(A), \{1\} \ominus F(A)) \tag{78}$$

- **Conjunction**

$$\mathfrak{N}(A_1 \cap A_2) = \mathfrak{N}(A_1) \boxdot \mathfrak{N}(A_2) = (T_1 \odot T_2, I_1 \odot I_2, F_1 \odot F_2) \tag{79}$$

- **Disjunction**

$$\begin{aligned}\mathfrak{N}(A_1 \cup A_2) &= (T_1 \cup T_2, I_1 \cup I_2, F_1 \cup F_2) \\ &= ((T_1 \oplus T_2) \ominus (T_1 \odot T_2), (I_1 \oplus I_2) \ominus (I_1 \odot I_2), (I_1 \oplus I_2) \ominus (I_1 \odot I_2)) \\ &= [\mathfrak{N}(A_1) \boxplus \mathfrak{N}(A_2)] \boxminus [\mathfrak{N}(A_1) \boxdot \mathfrak{N}(A_2)]\end{aligned} \tag{80}$$

$\mathfrak{N}$ – **Membership function over a neutrosophic set**

Let $\Theta$ be a world of discourse (called frame of discernment in the DST). Each $\mathfrak{N}$−element $x$ of $\Theta$ is characterized by its own neutrosophical basic assignment ($\mathfrak{N}$– value) $\mathfrak{N}(x) \triangleq (T(x), I(x), F(x))$ with $T(x), I(x)$ and $F(x) \subset ]^{-}0; 1^{+}[$. The $\mathfrak{N}$– membership function of any neutrosophical element $x$ with any subset $M \subset \Theta$ is defined in similar way by

$$\mathfrak{N}(x \mid M) \triangleq (T_M(x), I_M(x), F_M(x)) \tag{81}$$

with $T_M(x), I_M(x)$ and $F_M(x) \subset ]^{-}0; 1^{+}[$. The $\mathfrak{N}$– value of $x$ over $M$ can be interpreted, by abuse of language, as its membership function to $M$ in the following sense: $x$ is $t\%$ true in the set $M$, $i\%$ indeterminate (unknown if it is) in $M$, and $f\%$ false in $M$, where $t$ varies in $T$, $i$ varies in $I$, $f$ varies in $F$. The standard notation $x \in M$ will be used in the sequel to denote the neutrosophical membership of $x$ to $M$. One can say actually that any element $x$ of a given frame of discernment supported by a body of evidence neutrosophically belongs to any set, due to the percentages of truth/indeterminacy/falsity involved, which varies between 0 and 1 or even less than 0 or greater than 1. From this definition and previous neutrosophic rules, one gets directly the following basic neutrosophical set operations :

- **Complement of** $M$

If $x \in M$ with $\mathfrak{N}(x \mid M) \triangleq (T_M(x), I_M(x), F_M(x))$, then $x \notin M$ with

$$\mathfrak{N}(x \mid \bar{M}) = (\{1\} \ominus T_M(x), \{1\} \ominus I_M(x), \{1\} \ominus F_M(x)) \tag{82}$$

- **Intersection** $M \cap N$

If $x \in M$ with $\mathfrak{N}(x \mid M) \triangleq (T_M(x), I_M(x), F_M(x))$ and $x \in N$ with $\mathfrak{N}_{|W}(x \mid N) \triangleq (T_N(x), I_N(x), F_N(x))$, then $x \in M \cap N$ with

$$\mathfrak{N}(x \mid M \cap N) = (T_M(x) \odot T_N(x), I_M(x) \odot I_N(x), F_M(x) \odot F_N(x)) \tag{83}$$

- **Union** $M \cup N$

If $x \in M$ with $\mathfrak{N}(x \mid M) \triangleq (T_M(x), I_M(x), F_M(x))$ and $x \in N$ with $\mathfrak{N}(x \mid N) \triangleq (T_N(x), I_N(x), F_N(x))$, then $x \in M \cup N$ with

$$\mathfrak{N}(x \mid M \cup N) = (T_{M \cup N}(x), I_{M \cup N}(x), F_{M \cup N}(x)) \tag{84}$$

where

$$T_{M \cup N}(x) \triangleq [T_M(x) \oplus T_N(x)] \ominus [T_M(x) \odot T_N(x)] \tag{85}$$

$$I_{M \cup N}(x) \triangleq [I_M(x) \oplus I_N(x)] \ominus [I_M(x) \odot I_N(x)] \tag{86}$$

$$F_{M \cup N}(x) \triangleq [F_M(x) \oplus F_N(x)] \ominus [F_M(x) \odot F_N(x)] \tag{87}$$

- **Difference** $M - N$

    Since $M - N \triangleq M - \bar{N}$, if $x \in M$ with $\mathfrak{N}(x \mid M) \triangleq (T_M(x), I_M(x), F_M(x))$ and $x \in N$ with $\mathfrak{N}(x \mid N) \triangleq (T_N(x), I_N(x), F_N(x))$, then $x \in M - N$ with

    $$\mathfrak{N}(x \mid M - N) = (T_{M-N}(x), I_{M-N}(x), F_{M-N}(x)) \tag{88}$$

    where

    $$T_{M-N}(x) \triangleq T_M(x) \ominus [T_M(x) \odot T_N(x)] \tag{89}$$
    $$I_{M-N}(x) \triangleq I_M(x) \ominus [I_M(x) \odot I_N(x)] \tag{90}$$
    $$F_{M-N}(x) \triangleq F_M(x) \ominus [F_M(x) \odot F_N(x)] \tag{91}$$

- **Inclusion** $M \subset N$

    We will said that $M \subset N$ if for all $x \in M$ with $\mathfrak{N}(x \mid M) \triangleq (T_M(x), I_M(x), F_M(x))$ and $x \in N$ with $\mathfrak{N}(x \mid N) \triangleq (T_N(x), I_N(x), F_N(x))$, one has jointly $T_M(x) \subset T_N(x)$, $I_M(x) \subset I_N(x)$ and $F_M(x) \subset F_N(x)$.

## 5.2 Combination of neutrosophic evidences

Let's consider a general finite frame of discernment $\Theta = \{\theta_1, \ldots, \theta_n\}$ and two bodies of (neutrosophic) evidence $\mathcal{B}_1$ and $\mathcal{B}_2$. In the neutrosophic framework, we assume that each body of evidence provides some report of evidence (i.e. $\mathfrak{N}$ – value) committed to some elements of the hyper-power set $D^\Theta$. In other words, the information one has to deal with is the reports:

- Report for $\mathcal{B}_1$ : $R_1 = \{\mathfrak{N}_1(A_1), \ldots, \mathfrak{N}_1(A_m)\}$ for $A_1, \ldots, A_m \in D^\Theta$
- Report for $\mathcal{B}_2$ : $R_2 = \{\mathfrak{N}_2(B_1), \ldots, \mathfrak{N}_2(B_n)\}$ for $B_1, \ldots, B_n \in D^\Theta$

where each neutrosophic value for a proposition corresponds actually to a given triplet $(T(.), I(.), F(.)) \subset ]^-0; 1^+[^3$. Within the neutrosophic logic, one has the full degree of freedom between the $\mathfrak{N}$ – values for a report.

Our major concern now is to solve the difficult question on how to combine such kind of information to get the global and most pertinent information about the problem under consideration. So, is it possible to construct a new global report (and hopefully more informative) $R$ from $R_1$ and $R_2$ ? Unfortunately, the neutrosophic logic which is a new appealing and modelling tool to deal with uncertainties on propositions of same universe of discourse does not provide a clear and direct mathematical mechanism for dealing with combination of such kind of evidences. We propose in this section a possible issue for this important question based on our new generalization of the DST.

The main idea for combining such kind of evidences is to convert the reports into two proper general bpa $m_{R_1}(.)$ and $m_{R_2}(.)$ and then combine them using the general rule of combination (41). The combination of neutrosophic evidences is a two-level process.

**Level 1 : the general bpa transformation**

The major difficulty is the mapping of the set of neutrosophic values $\{\mathfrak{N}(.)\}$ into a set of corresponding elementary bpa $m(.)$. Several cases are now examined.

- Case 1 (simpliest case) : We assume that each neutrosophic evidence corresponds only to a triplet of real positive or null numbers belonging to $[0; 1]$ (i.e. $T(.)$, $I(.)$ and $F(.)$ are restricted to real numbers $\in [0; 1]$).
    Since in the neutrosophic logic, $T(.)$, $I(.)$ and $F(.)$ have no strong mathematical relationships, the easiest solution within the classical DST would be to use the following transformation

    $$m_{(}A) = T(A)/c \qquad m(A^c) = F(A)/c \qquad m(A \cup A^c) = I(A)/c$$

    where $c$ is a normalization constant such that $m_{(}A) + m(A^c) + m(A \cap A^c) = 1$.

    In our general theory of plausible and paradoxical reasoning, it seems more judicious to use the following mapping based on our general modelling of information granule described in section 4.4.8. Thus, we are now able to construct

from $T(.)$, $I(.)$ and $F(.)$ the three corresponding elementary bpa as follows for any proposition $C \in D^\Theta$ involved in a given report :

$$
\begin{array}{lll}
m_1(A) = T(A) - \frac{1}{2}m_1^\star & m_2(A) = 1 - F(A) - \frac{1}{2}m_2^\star & m_3(A) = 0 \\
m_1(A^c) = 1 - T(A) - \frac{1}{2}m_1^\star & m_2(A^c) = F(A) - \frac{1}{2}m_2^\star & m_3(A^c) = 0 \\
m_1(A \cap A^c) = m_1^\star & m_2(A \cap A^c) = m_2^\star & m_3(A \cap A^c) = 1 - I(A) \\
m_1(A \cup A^c) = 0 & m_2(A \cup A^c) = 0 & m_3(A \cup A^c) = I(A)
\end{array}
$$

where $m_1^\star$ is given by the solution of equation $64e^2(m_1^\star)^4 - (m_1^\star)^2 + 2m_1^\star - 4(1-T(A))T(A) = 0$ and $m_2^\star$ by the solution of equation $64e^2(m_2^\star)^4 - (m_2^\star)^2 + 2m_2^\star - 4(1-F(A))F(A) = 0$. The mapping $m_3(.)$ comes from the necessity to not assign a prior preference to $A$ rather than to $A^c$ when only indeterminacy is available.

- <u>Case 2</u> : We assume now that each neutrosophic evidence corresponds only to a triplet of real intervals belonging to $[0;1]$. In this case, the more general mapping is proposed.

$$
\begin{array}{lll}
m_1(A) = m_T - \frac{1}{2}m_1^\star & m_2(A) = 1 - M_F - \frac{1}{2}m_2^\star & m_3(A) = (M_I - m_I)/2 \\
m_1(A^c) = 1 - M_T - \frac{1}{2}m_1^\star & m_2(A^c) = m_F - \frac{1}{2}m_2^\star & m_3(A^c) = (M_I - m_I)/2 \\
m_1(A \cap A^c) = m_1^\star & m_2(A \cap A^c) = m_2^\star & m_3(A \cap A^c) = 1 - M_I \\
m_1(A \cup A^c) = M_T - m_T & m_2(A \cup A^c) = M_F - m_F & m_3(A \cup A^c) = m_I
\end{array}
$$

where $m_T \triangleq \text{Inf}(T(A))$, $M_T \triangleq \text{Sup}(T(A))$, $m_F \triangleq \text{Inf}(F(A))$, $M_F \triangleq \text{Sup}(F(A))$ and $m_I \triangleq \text{Inf}(I(A))$, $M_I \triangleq \text{Sup}(I(A))$. $m_1^\star$ is given by the solution of equation $64e^2(m_1^\star)^4 - (m_1^\star)^2 + 2(1 - M_T + m_T)m_1^\star - 4(1 - M_T)m_T = 0$ and $m_2^\star$ by the solution of equation $64e^2(m_2^\star)^4 - (m_2^\star)^2 + 2(1 - M_F + m_F)m_2^\star - 4(1 - M_I)m_I = 0$.

- <u>Case 3</u> (general case) : We assume now that each component of neutrosophic value $(T(A) = \bigcup_i T_i(A), I(A) = \bigcup_j I_j(A), F(A) = \bigcup_k F_k(A))$ is actually the union of subintervals of $[0;1]$. In such general case, we propose to construct for each possible combinations of $(T_i(A), I_j(A), F_k(A))$ a corresponding general bpa as for case 2 then combine all bpa using the general rule of combination to get the global bpa relative to the proposition under consideration.

**Level 2 : the combination of evidences**

We have just shown how general elementary bpa can be evaluated from each neutrosophic values of a report. For the report $R_1$, we have now in hands a set of bpa $m_1(.), \ldots m_m(.)$ associated to every proposition in this report. Similarly, we get also another set of bpa $m_1'(.), \ldots m_n'(.)$ for report $R_2$. For each set of bpa, we are now able to compute the global general bpa $m_{R_1}(.)$ and $m_{R_2}(.)$ from the general rule of combination (41) by

$$m_{R_1} = m_1 \oplus m_2 \oplus \ldots \oplus m_m$$

$$m_{R_2} = m_1' \oplus m_2' \oplus \ldots \oplus m_n'$$

The next step of the combination is then to combine the bpa $m_{R_1}$ with $m_{R_2}$ by applying for the last time the general rule of combination (41) to finally get the global result we are looking for; i.e.

$$m(.) = m_{R_1} \oplus m_{R_2}$$

From the global bpa $m(.)$ defined on the hyper-power set $D^\Theta$, we will then be able to evaluate the degree of belief of each proposition of $D^\Theta$ which will help us to take the most pertinent decision for the problem under consideration.

# 6  Conclusion

In this paper, the foundations for a new theory of paradoxical and plausible reasoning has been developed which takes into account in the combination process itself the possibility for uncertain and paradoxical information. The basis for the development of this theory is to work with the hyper-power set of the frame of discernment relative to the problem under consideration rather than its classical power set since, in general, the frame of discernment cannot be fully described in terms of an exhaustive and exclusive list of disjoint elementary hypotheses. In such general case, no refinement is possible to apply directly the Dempster-Shafer theory (DST) of evidence. In our new theory, the rule of combination is justified from the maximum entropy principle and there is no mathematical impossibility to combine sources of evidence even if they appear at first glance in contradiction (in the Shafer's sense) since the paradox between sources is fully taken into account in our formalism. We have also shown that in general, the combination of evidence yields unavoidable paradoxes. This theory has shown, through many illustrated examples, that conclusions drawn from it, provides results which agree perfectly with the human reasoning and is useful to take a decision on complex problems where DST usually fails. The last part of this work has been devoted to the development of a theoretical bridge between the neutrosophic logic and this new theory, in order to solve the delicate problem of the combination of neutrosophic evidences. The neutrosophic logic serves here as the most general framework (prerequesite) for dealing with uncertain and paradoxical sources of information through this new theory.

# Dynamic Fuzzy Sets and Neutrosophic Concepts

Andrzej Buller

ATR Human Information Science Laboratories
2-2-2 Hikaridai, Seika-cho, Soraku-gun, Kyoto 619-0288, Japan
buller@atr.co.jp

*Abstract:* This paper presents the concept of *Dynamic Fuzzy Set* (DFS) as an instance of *Neutrosophic Set* (NS) useful in the modeling of dynamics of mental processes. For each element of a NS there is a triple (T, F, I) understood as functions/operators dependent on time, space, and other not necessarily known parameters. T determines a degree of truth that a given element is in the set, while F and I are for a degree of false and indeterminacy, respectively. DFS has been defined as such a set that for each of its elements there is a function that for a given time returns a value of membership of a given element to the set. Based on DFS one can perform an extraordinary fuzzy inferencing as battle between populations of copies of contradictory statements. It has been observed that resulting membership values may change in time even when the processed data remain constant. Psychological justification of DFS came from empirically confirmed fact that that subject's feelings about a perceived person or social situation can oscillate from a highly positive value to a highly negative value and back even in absence of new data about related objects.

## 1 INTRODUCTION

The concept of the fuzzy set, invented in the 60's [8] has proved to reflect the way humans categorize [4]. However, while categorization is an important determinant of human behavior, it is not the only determinant. In order to model mental mechanisms underlying human interactions with surrounding reality a number of logics has been invented, however each of them deals rather with a selected isolated facet of mental activity. Neutrosophic approach proposed by Florentin Smarandache tries to bring all the logics together towards a unitary, formally analyzable, general model of reality [5]. A part of the neutrosophic model called *Neutrosophic Set* (NS) intended to generalize fuzzy set. From another attempt to generalize the concept of fuzzy set the idea of *Dynamic Fuzzy Set* (DFS) emerged. Let as, therefore, compare the two definitions:

*Definition 1* (Smarandache [5: 77]):
    **Neutrosophic Set** is a set such that an element belongs to the set with a neutrosophic probability, i.e. t% is true that the element is in the set, f% false, and i% indeterminate.

*Definition 2* (based on Buller [1]):
    **Dynamic Fuzzy Set** is a fuzzy set such that an element belongs to the set with a membership value that changes in time.

In [5: 78] some examples of NS application suggesting that NS much better than classic fuzzy set reflect complex truth about reality one can find can be found. Nevertheless, the early definition of NS cited above says nothing about a space to which t%, f% and i% belong. If t%, f% and i% were interpreted as real numbers, DFS could be treated as a different sort of generalization of fuzzy set, not included in the concept of NS. But the concept of NS evolved. Let's present another definition:

*Definition 3* (based on Smarandache [8: 111] and Smarandache [9]):

**Neutrosophic Set** included in a universe of discourse is a set such that an element from the universe belongs to the set in such a way that t% is true that the element is in the set, f% false, and i% indeterminate, where t varies in T, f varies in F, i varies in I, where T, I, F are functions/operators depending on many known and unknown parameters; here T(…), I(…), F(…) are standard or non-standard subsets included in the non-standard unit interval $]^-0, 1^+[$.

Compared to the Fuzzy Set, the Neutrosophic Set can distinguish between 'absolute membership' (appurtenance) of an element to a set ($T=1^+$), and 'relative membership' ($T=1$), whereas the 'partial membership' is represented by $0 < T < 1$. Also, the sum of neutrosophic membership components (truth, indeterminacy, falsehood) are not required to be 1 as in fuzzy membership components, but may be any number between 0 and 3.

Since *Definition 3* allows T, I and F to be dependent on time (as one of the known parameters), DFS may be recognized as an instance of NS such that, for example, t% : **T** $\rightarrow$ [0, 1], i = 0 = *const.*, f% : **T** $\rightarrow$ [0, 1], $\forall_{t \in \mathbf{T}}$ f% = 1 – t%, where **T** is a space of integers representing moments in time.

The importance of DFS comes from increasing interest in dynamic processes in psychology. The development of social cognition caused in the 90's an increasing interest in intrinsic dynamics of human categorization. Experiments confirmed that when judging a perceived person or social situation people sometimes oscillate from highly positive feelings to highly negative feelings even in absence of new data [3]. Similar oscillations were demonstrated by computational models of human working memory that treated mental process as a "debate" in a "society of memes" in a cellular working memory [2]. The idea of *Dynamic Fuzzy Calculus* (DFC) based on DFS emerged from the research on the models.

## 2 Dynamic Fuzzy Calculus

Let **E** = {$\emptyset$, $e_1$, $\neg e_1$, $e_2$, $\neg e_2$, …} be the space of particular *notions*, where $\emptyset$ denotes a *notion of reference* while $\neg$ is the operator of *negation*. Let **T** be the space of integers representing *time*, while **A** be a space of functions such that $\forall_{\mu \in \mathbf{A}} \mu : \mathbf{E}^2 \rightarrow [0, 1]^K$, where K is a positive integer.

Let the entity of interest be a space **B°** such that
$$\forall_{a \in \mathbf{B°}} a : \mathbf{T} \rightarrow \mathbf{A}$$
Any element of **B°** may be considered as a system of relationships between notions. Each of the relationships applies to a pair of notions and expresses itself as a K-element vector of real numbers not lesser than 0 not greater than 1. Elements of the vector can be interpreted as strengths of several sorts of memberships of the second notion in a pair to a DFS represented by the first notion in the pair. Each of the strengths may, by definition, change over time.

When a system consists of a set of related notions and the relationships between the notions change in time we say could say that in the system a kind of mental activity takes palace. Buller [1] argues that thinking is nothing but changing relationships between notions, and concludes that any element of the space **B°** could be considered to be a model of a *mind* or even to be a mind itself. The thing that is the source of inconstancy of the relationships within **B°** is a *brain* that can be understood as any machine dedicated to a given $a \in$ **B°** that processes every $a_t$ onto $a_{t+1}$.

The space **B′** contains minds that use brains understood as external devices to determine relationships between specific notions for consecutive moments of time. Formally:

$$a \in \mathbf{B'} \Leftrightarrow (\ a \in \mathbf{B^\circ}, \exists_{\langle x, F \rangle} x = (x_0, x_1, x_2, \ldots, x_N),$$
$$x : \mathbf{T} \to \mathbf{M}^N, F : \mathbf{A} \times \mathbf{M}^N \to \mathbf{A} \times \mathbf{M}^N,$$
$$(a_{t+1}, x_{t+1}) = F(a_t, x_t)\ ),$$

where $\mathbf{M} = \mathbf{E}^2$ is the space of *memes*, while N is a positive integer. The brain may be, therefore, considered as the couple $\langle x, F \rangle$.

Appendix contains a description of an instance of $a \in$ **B′** reflecting an isolated process of social judgment. Some simulation results appeared to fit the psychological evidence that subjects may oscillate from highly positive feelings to highly negative feelings without new data about objects of reference, i.e. based exclusively on intrinsic dynamics of mental process. It has been suggested that the membership of reference tended to a 2-state limit cycle attractor.

### Neurotrophic dynamics

Let us consider a situation wherein somebody becomes a victim of a sexual assault. Among a number of possible plots of the victim's feelings towards the perpetrator let us consider: (a) start from extreme hatred in $t_0$ and then a linear decay of the hatred towards indifference in $t_k$, and (b) constant extreme hatred that in 0.75 $t_k$ turns suddenly into perverse love. In terms of DFS we can describe each of the cases via providing of appropriate formula for function $\mu_{H,t}$, where $\mu_{H,t}=1$ means that the perpetrator in time t fully belongs to the dynamic fuzzy set of people extremely hated by the victim, $\mu_{H,t}=0$ may mean that the perpetrator in time t fully belongs to the dynamic fuzzy set of people extremely loved by the victim, while $\mu_{H,t}=0.5$ may mean that the perpetrator in time t is perfectly indifferent to the victim. If one wanted to describe the victim's feelings for t∈ [ $t_0$, $t_k$] in terms of classic fuzzy sets, the simplest solution would be to provide $\mu_H$ as the average value of $\mu_{H,t}$ in the period [ $t_0$, $t_k$]. The unavoidable price for this simplicity is loss of the knowledge ablot the essence of the victim's mental states. Indeed, in both considered cases the $\mu_H$ would take the same value 0.75. As for full neutrosophic generalization, much more knowledge about subjects' mental states could be added to the model if the neutrosophic component I were included and allowed to change in time. However, the question how to measure subjects' indeterminacy, especially in reference to social judgment, remains open.

### Conclusions

Psychological evidence provides justification for cognitive models based on Dynamic Fuzzy Sets (DFS) and Dynamic Fuzzy Calculus that deal with membership values changing over time. Subjects' oscillation from highly positive feelings to highly negative ones and back when knowledge about an object of interest remains constant can be interpreted in terms of membership to a dynamic fuzzy set tending to a 2-state limit cycle attractor. DFS-based modeling also seems to lead to the philosophical thesis that mind is nothing but relationships that occur between certain notions, and that change over time. Neutrosophic approach intends to generalize existing logics to the most possible extent. DFS can be assumed to be an instance of Neutrsophic Set. A possibility generalization of DFS in the course of adding of more neutrosophic components, for example, indeterminacy factor, encourages to search for new ways of experimental exploration of dynamics of social judgment and for new ways of mathematical description of mental processes.

**Appendix.** *Dynamic Fuzzy Calculus in a model of social judgment* (based on [1])

Let us consider an instance of $a \in \mathbf{B}'$ such that:

1. K=2;

2. $\mathbf{E} = \{\emptyset, \mathbf{N}, \mathbf{n}, \mathbf{R}, \mathbf{r}, \mathbf{A}, \mathbf{a}\}$, where $\mathbf{n} = \neg \mathbf{N}$, $\mathbf{r} = \neg \mathbf{R}$, $\mathbf{a} = \neg \mathbf{A}$;

3. $\forall_{t \in \mathbf{T}} \mathbf{m} \notin \{\langle s|m \rangle \in \mathbf{E}^2 \mid m \neq \emptyset \Rightarrow s \notin \{\emptyset, m\}\} \Rightarrow \mu_{s,t}(m) = 0$;

4. $\forall_{t \in \mathbf{T}} \mu_{s,t}(\neg m) = 1 - \mu_{s,t}(m)$;

5. $\forall_{t \in \mathbf{T}} x_{0,t} = \langle \emptyset | \emptyset \rangle$,

6. if $u_{i,t} = \langle \varnothing | \varnothing \rangle$ then $x_{i,t} = v_t$, while if $u_{i,t} \neq \langle \varnothing | \varnothing \rangle$ then $x_{i,t} = u_{i,t}$, where

7. $\forall_{\mathbf{m'} \in \mathbf{M}} P(v_t = \mathbf{m'}) = \mu_1(\mathbf{m'}, t)\mu_2(\mathbf{m'}, t) / \sum_{\mathbf{m} \in \mathbf{M}} \mu_1(\mathbf{m}, t)\mu_2(\mathbf{m}, t)$, where $\forall_Z P(Z)$ is probability of Z;

8. $u_{i,t} : \mathbf{Z}_N \times \mathbf{T} \to \mathbf{M}$, $\forall$ i,j,k$\in \mathbf{Z}_N$, $u_{i,t} = \psi(x_{j,t}, x_{k,t})$, where $j = L_i(\mathbf{x}_t)$, $k = R_i(\mathbf{x}_t)$, where $\mathbf{Z}_N = \{\iota \in \mathbf{I} \mid 0 \leq \iota \leq N\}$, where $\mathbf{Z}$ is the space of integers, while

9. $\psi : \mathbf{M}^2 \to \mathbf{M}$, $\forall_{m \in \mathbf{M}} \forall_{\alpha, \beta \in \mathbf{E} \mid \alpha \neq \beta, \neg \alpha \neq \beta}$
$m \in \{ \langle \beta | \alpha \rangle, \langle \neg \alpha | \varnothing \rangle \} \Rightarrow \psi(\langle \alpha | \varnothing \rangle, m) = \langle \varnothing | \varnothing \rangle$;
$m \notin \{ \langle \beta | \alpha \rangle, \langle \neg \alpha | \varnothing \rangle \} \Rightarrow \psi(\langle \alpha | \varnothing \rangle, m) = \langle \alpha | \varnothing \rangle$;
$m \notin \{ \langle \alpha | \varnothing \rangle, \langle \neg \alpha | \varnothing \rangle \} \Rightarrow \psi(\langle \beta | \alpha \rangle, m) = \langle \beta | \alpha \rangle$;
$\psi(\langle \beta | \alpha \rangle, \langle \alpha | \varnothing \rangle) = \langle \beta | \varnothing \rangle$; $\psi(\langle \beta | \alpha \rangle, \langle \neg \alpha | \varnothing \rangle) = \langle \neg \beta | \varnothing \rangle$;

10. $L, R : \mathbf{Z}_N \times \mathbf{M}^N \to \mathbf{Z}_N$,

| | | | | | | | | |
|---|---|---|---|---|---|---|---|---|
| $x_{p,t} \neq \langle \varnothing \mid \varnothing \rangle$ | 0 | 0 | 0 | 0 | 1 | 1 | 1 | 1 |
| $x_{q,t} \neq \langle \varnothing \mid \varnothing \rangle$ | 0 | 1 | 0 | 1 | 0 | 1 | 0 | 1 |
| $x_{r,t} \neq \langle \varnothing \mid \varnothing \rangle$ | 0 | 0 | 1 | 1 | 0 | 0 | 1 | 1 |
| $L_{w(\varphi,\lambda,0)}(\mathbf{x}_t)$ | r | p | r | p | r | p | p | p |
| $R_{w(\varphi,\lambda,0)}(\mathbf{x}_t)$ | q | q | p | q | p | q | r | r |
| $L_{w(\varphi,\lambda,1)}(\mathbf{x}_t)$ | q | q | q | q | q | q | q | q |
| $R_{w(\varphi,\lambda,1)}(\mathbf{x}_t)$ | p | p | p | p | p | p | p | 0 |
| $L_{w(\varphi,\lambda,2)}(\mathbf{x}_t)$ | r | r | r | r | p | r | r | r |
| $R_{w(\varphi,\lambda,2)}(\mathbf{x}_t)$ | p | p | q | q | q | q | p | p |

where $p = w(\varphi_p, \lambda_p, 2)$, $q = w(\varphi_q, \lambda_q, 1)$, $r = w(\varphi_r, \lambda_r, 0)$, $\varphi_p = (\varphi + \varphi_{max}) \bmod \varphi_{max} + 1$, $\varphi_q = \varphi$, $\varphi_r = (\varphi + 1) \bmod \varphi_{max} + 1$, $\lambda_q = (\lambda + \lambda_{max}) \bmod \lambda_{max} + 1$, if $\varphi \bmod 2 \neq 0$ then $\lambda_p = \lambda_r = \lambda$, if $\varphi \bmod 2 = 0$ then $\lambda_p = \lambda_r = (\lambda + \lambda_{max}) \bmod \lambda_{max} + 1$, where

11. $w$ is any one-by-one function such that $w : \mathbf{G} \to \mathbf{Z}_N$, where $\mathbf{G} = \{(\varphi, \lambda, \delta) \in \mathbf{Z} \mid 0 \leq \varphi \leq \varphi_{max}, 0 \leq \lambda \leq \lambda_{max}, 0 \leq \delta \leq 2\} \cup \{e\}$, $(\varphi_{max}, + 1)(\lambda_{max} + 1) = N$, e is any element such that $e \neq \varnothing$;

12. $\forall_{\mathbf{m} \notin \{\langle \mathbf{A} | \varnothing \rangle, \langle \mathbf{a} | \varnothing \rangle\}} \mu_1(\mathbf{m}, t) = $ const., while $\mu_1(\langle \mathbf{A} | \varnothing \rangle, t) = N_\mathbf{A} / (N_\mathbf{A} + N_\mathbf{a})$, where for a given notion e, $N_e = C_{e, N}$, $C_{e, 0} = 0$, $C_{e, i+1} = C_{e, i} + j$, $j = 1$ if $x_{i, t} = e$, while $j = 0$ if $x_{i, t} \neq e$; where $C_{J, i}$ is an auxiliary counter.

13. $\forall_{\mathbf{m} \in \mathbf{M}} \mu_2(\mathbf{m}, t) = $ const.

Assuming that **N** and **R** represent 'nicety' and 'richness' of a date proponent, respectively, while **A**|∅ is a meme suggesting subject's readiness to having the date, let us simulate the plot of $\mu_1(\langle \mathbf{A} | \varnothing \rangle, t)$. For given constant values $\mu_1(\langle \mathbf{N} | \varnothing \rangle, t) = .6$, $\mu_1(\langle \mathbf{R} | \varnothing \rangle, t) = .4$, $\mu_1(\langle \mathbf{A} | \mathbf{N} \rangle, t) = \mu_1(\langle \mathbf{A} | \mathbf{R} \rangle, t) = 1.0$, $\mu_2(\langle \mathbf{N} | \varnothing \rangle, t) = \mu_2(\langle \mathbf{R} | \varnothing \rangle, t) = .67$, $\mu_2(\langle \mathbf{A} | \mathbf{N} \rangle, t) = \mu_2(\langle \mathbf{A} | \mathbf{R} \rangle, t) = .2$, F produces a plot of $\mu_1(\langle \mathbf{A} | \varnothing \rangle, t)$ such as, for examle:

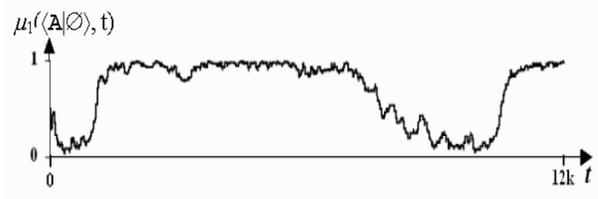

# NeuroFuzzy and Neutrosophic Approach to Compute the Rate of Change in New Economies

*Dr. Mohammad Khoshnevisan*

School of Accounting and Finance, Griffith University (Gold Coast)

PMB 50 Gold Coast Mail Center, Queensland 9726, Australia

*Dr. Sarjinder Singh*

Department of Mathematics and Statistics, University of Saskatchewan, Canada

ABSTRACT:

Probabilistic methods have been used to model uncertainty in many optimization problems, an alternative is to use NeuroFuzzy and Neutrosophic control methods to model uncertainties with new economies. New economies have a solid impact on the growth of financial markets. In this day of age, one needs to be able to monitor these rapid dynamical changes with as much precision as possible. The NeuroFuzzy and Neutrosophic control based optimization minimizes the system possibility of failure.

Therefore, we are able to formulate our optimal control model with significant reliability. NeuroFuzzy systems combine the advantages of fuzzy systems - the transparent representation of knowledge and the ability to cope with uncertainties with the advantages of neural nets. In this paper, we have proposed and designed a NeuroFuzzy and Neutrosophic and control system for analysing the impact of new economies growth and decline.

1. **Introduction:**

In recent years, substantial progress has been made in our ability to model and deal with uncertainty through computationally - orientated techniques. In finance, there is a long-standing tradition that probability theory is the only available tool for dealing with a stochastic nature. Recently, however the validity has been called into question by researchers, the development of other tools - centered on possibility theory and neurofuzzy control for dealing with uncertainty. Sengupta created a model of Shumpeterian Dynamics. In that model, two key variables were introduced through a two staged formulation using Schumpeterian dynamics. Sengupta formulated a model for innovation input in the form of knowledge capitol, the stochastic process followed as a birth and death process mechanism, where $B_t$ represents new ideas and $D_t$ represents the relative obsolescence or destruction of the old. In this paper we will illustrate a comparison of modeling using Schumpeterian dynamics versus Neurofuzzy and Neutrosophic control modeling solely emphasizing this application to new economies. As we all are aware of the hundreds of new companies entering the market - $B_t$ - birth at the same instance there are companies exiting the market $D_t$ - death. we will demonstrate the entrance and exiting phenomena using Neurofuzzy control.

## 2. The stochastic model proposed by Sengupta

Sengupta[11,12] has given the mean the variance of the innovation input as follows:

$E(x_{t+1}) = x_t \, e^{b_t - d_t} = e^{r_1} x_t$, $r_1 = b_t - d_t$ (2.1)

and $\text{Var}(x_{t+1}) = b_t + d_t / b_t - d_t \, [\, e^{r_1} - 1] \, e^{r_1} x_t$ if $b_t \neq d_t$ (2.2)

$= 2 b_t d_t$. One can observe that equation (2.2) does not hold when $b_t = d_t$. That is, if $b_t = d_t$, the variance $(x_{t+1})$ is not defined and is not equal to $2 b_t d_t$. The elimination of (2.2) will lead to some jump discontinuity of the growth model

Given the stochastic innovation process $E(x_{t+1})=x_t+B_t-D_t$ and using the assumption that the birth and death flows are endogenously determined. That is,

$B_t=b_0+b_1x_t+b_2x_{t-1}$ (2.3)

$D_t=d_0+d_1x_t+d_2x_{t-1}$ (2.4)

It is noted that $x_t= -\beta_0 +\beta_1 x_{t+1} +\beta_2 x_{t+1}$.

The generalized formula of the past history and future expectations given in his paper will only consider the case when $b_t \neq d_t$. In his model $x_t$ has been expressed as values of t +1 with different coefficients $\beta_1$ and $\beta_2$ and in the general form which is given by

$x_t= -\beta_0+ \sum_{i=1}^{m} \beta_{11i} x_{t-i} + \sum_{i=1}^{n} \beta_{2\,2i} x_{t+i}$ (2.5)

$x_t$ has been expressed as functions of t+1 and t-1.

Equation (2.5) will lead to some designated path of optimization for only a special case when $b_t \neq d_t$ are and it does not include all the possibilities such as when $\beta_1=\beta_2$.

We have shown when $b_t=d_t$ the variance is not defined and the $E(x_{t+1})=x_t$.

It is claimed that when $b_2>d_2$ with $b_1 \geq d_1$ the net result is positive growth and when $b_2<d_2$ with $b_1 \geq d_1$ the future impact through demand pull is much stronger than the past trend in innovation. This is true only if $b_t \neq d_t$. In this model, it is assumed that birth and death of new economies take two absolute forms. i.e. being true or false. We have considered a Neutrosophic design for three states being T,I,F. We have an alternative approach for dealing with this model. In section, 3 we will build our NeuroFuzzy and Neutrosophic system.

### 3. Fuzzy and Neutrosophic rules and three state Design for the birth and death rate of new economies

The fuzzy rules in our model is given as follows:

IF $b_t$ is Positive-Big AND $d_t$ is Positive-Small THEN the control output Growth=1
(T)

IF $b_t$ is Positive-Small AND $d_t$ is Positive-Big THEN the control output Growth=2
(I)

IF is Positive-Small AND $d_t$ is Positive-Small THEN the control output Growth=3 (T)

IF is Positive-Big AND $d_t$ is Positive-Big THEN the control output Growth=3 (F)

These functions results in a relationship matrix with a membership function defined as $\mu_{3R}(b_t,d_t,growth)=\mu_{A \times B}(b_t,d_t)\text{-------}> \mu_E(growth)===>(T,I,F)$

$\mu_{3R}(b_t,d_t,growth)$ is the membership value of error element ($b_t$), change in error element ($d_t$) and the resulting control action growth in the three dimensional implication function ( A x B---->E).

The first step is to classify our neural model. We have used 90 samples and applied Neutrosophic classification ( T,I,F) in three states ( T.I,F) , the Neutrosophic approach is to cast out undesirable clusters in our training data. The linguistic variables, terms, membership functions, rule blocks, rules and interfaces have been incorporated, in our computation, the winner neuron is 4 meaning that in every training iteration 4 rules are applied. This was the optimum for our model after 40 runs were executed. In our model, the average error of all samples is computed and Neutrosophic error of the worse sample is determined. These values indicate the dynamical changes of the birth and death of new economies can be measured in a set of optimum target shooting in three distinct forms ( T,I,F) . We have not included the long computational results in this short paper.

### 4. Computational Method:

Now we need to build our NeuroFuzzy system. We have run 10000 iteration using Fuzzytech software and as a result the parameters of our system has been chosen as follows:

Maximum step=200

Maximum Deviation=3%

Average Deviation=.1

Neural network learning method:

Batch Random

Step width =.03

**Winner Neuron=4**

Isodata cluster configuration

1 b(t) .03 .90 2 input

2 d(t) .01 .90 2 input

3 Growth 1 3 2 output

1 b(t) .03 .90 2 output

2 d(t) .01 .90 2 output

3 Growth 1 3 2 output

We have defined the defuzzification box and the best compromise is used in most control applications. The rule block is 1 with random Dos value generator.

## 5.Conclusions:

Qualitatively, these examples demonstrate some of the potential power of NeuroFuzzy and Netrosophic control as applied to new economies problems. Rapid changes by the introduction of ever changing technology have made the use of traditional statistical methods questionable. In financial markets one would need to know the rate of growth of new economies at any given time to be able to form a dynamical optimization pattern. The fine tuning of NeroFuzzy and Neutrosophic controllers are intuitive and would likely further benefit from automatic processing by neural networks or genetic algorithms. Future research will focus on finding a systematic approach to tuning the NeuroFuzzy and Neutrosophic controllers, which would be of great use in financial engineering and risk optimizations and analyzing stock price behavior.

# Definitions Derived from Neutrosophics (Addenda)

Florentin Smarandache
Department of Mathematics
University of New Mexico
Gallup, NM 87301, USA

*Abstract*: As an addenda to the previous papers on "Neutrosophy" and "Neutrosophic Logic" from these proceedings, thirty-two new definitions are presented, derived from neutrosophic set, neutrosophic probability, neutrosophic statistics, and neutrosophic logic. Each one is independent, short, and cross-referenced like in a dictionary. They constitute the author's future exploration.

*Keywords*: Fuzzy set, fuzzy logic; neutrosophic logic;
Neutrosophic set, intuitionistic set, paraconsistent set, faillibilist set, paradoxist set, pseudo-paradoxist set, tautological set, nihilist set, dialetheist set, trivialist set;
Classical probability and statistics, imprecise probability;
Neutrosophic probability and statistics, intuitionistic probability and statistics, paraconsistent probability and statistics, faillibilist probability and statistics, paradoxist probability and statistics, pseudo-paradoxist probability and statistics, tautological probability and statistics, nihilist probability and statistics, dialetheist probability and statistics, trivialist probability and statistics;
Neutrosophic logic, paradoxist logic (or paradoxism), pseudo-paradoxist logic (or pseudo-paradoxism), tautological logic (or tautologism).

*2000 MSC*:  03E99, 03-99, 03B99, 60A99, 62A01, 62-99.

*Introduction*:
As a consequence to [1], [3], [4-7] we display the below unusual extensions of definitions resulted/deviated from neutrosophics in the Set Theory, Probability, and Logic. Some of them are listed in the Dictionary of Computing [2]. Further development of these definitions (including properties, applications, etc.) is in our research plan.

*1. Definitions of New Sets*

====================================================

**1.1. Neutrosophic Set**:

<logic, mathematics> A set which generalizes many existing classes of sets, especially the fuzzy set.

Let U be a universe of discourse, and M a set included in U. An element x from U is noted, with respect to the set M, as x(T,I,F), and belongs to M in the following way: it is T% in the set (membership appurtenance), I% indeterminate (unknown if it is in the set), and F% not in the set (non-membership);
here T,I,F are real standard or non-standard subsets, included in the non-standard unit interval $]^-0, 1^+[$, representing truth, indeterminacy, and falsity percentages respectively.

Therefore: $^-0 \leq \inf(T) + \inf(I) + \inf(F) \leq \sup(T) + \sup(I) + \sup(F) \leq 3^+$.

Generalization of {classical set}, {fuzzy set}, {intuitionistic set}, {paraconsistent set}, {faillibilist set}, {paradoxist set}, {tautological set}, {nihilist set}, {dialetheist set}, {trivialist}.

Related to {neutrosophic logic}.

=====================================================

**1.2. Intuitionistic Set**:

<logic, mathematics> A set which provides incomplete information on its elements.

A class of {neutrosophic set} in which every element x is incompletely known, i.e. x(T,I,F) such that sup(T)+sup(I)+sup(F)<1;
here T,I,F are real standard or non-standard subsets, included in the non-standard unit interval $]^-0, 1^+[$, representing truth, indeterminacy, and falsity percentages respectively.

Contrast with {paraconsistent set}.

Related to {intuitionistic logic}.

=====================================================

**1.3. Paraconsistent Set**:

<logic, mathematics> A set which provides paraconsistent information on its elements.

A class of {neutrosophic set} in which every element x(T,I,F) has the
property that sup(T)+sup(I)+sup(F)>1;
here T,I,F are real standard or non-standard subsets, included in the
non-standard unit interval ]⁻0, 1⁺[, representing truth, indeterminacy,
and falsity percentages respectively.

Contrast with {intuitionistic set}.

Related to {paraconsistent logic}.

======================================================

**1.4. Faillibilist Set**:

<logic, mathematics> A set whose elements are uncertain.

A class of {neutrosophic set} in which every element x has a
percentage of indeterminacy, i.e. x(T,I,F) such that inf(I)>0;
here T,I,F are real standard or non-standard subsets, included
in the non-standard unit interval ]⁻0, 1⁺[, representing truth,
indeterminacy, and falsity percentages respectively.

Related to {faillibilism}.

======================================================

**1.5. Paradoxist Set**:

<logic, mathematics> A set which contains and doesn't contain
itself at the same time.

A class of {neutrosophic set} in which every element x(T,I,F) has
the form x(1,I,1), i.e. belongs 100% to the set and doesn't
belong 100% to the set simultaneously;
here T,I,F are real standard or non-standard subsets, included in
the non-standard unit interval ]⁻0, 1⁺[ , representing truth,
indeterminacy, and falsity percentages respectively.

Related to {paradoxism}.

======================================================

**1.6. Pseudo-Paradoxist Set**:

<logic, mathematics> A set which totally contains and partially doesn't contain
itself at the same time,

or partially contains and totally doesn't contain itself at the same time.

A class of {neutrosophic set} in which every element x(T,I,F) has
the form x(1,I,F) with $0<\inf(F)\leq\sup(F)<1$ or x(T,I,1) with $0<\inf(T)\leq\sup(T)<1$,
i.e. belongs 100% to the set and doesn't belong F% to the set simultaneously, with
$0<\inf(F)\leq\sup(F)<1$,
or belongs T% to the set and doesn't belong 100% to the set simultaneously, with
$0<\inf(T)\leq\sup(T)<1$;
here T,I,F are real standard or non-standard subsets, included in
the non-standard unit interval $]^-0, 1^+[$, representing truth,
indeterminacy, and falsity percentages respectively.

Related to {pseudo-paradoxism}.

========================================================

**1.7. Tautological Set**:

<logic, mathematics> A set whose elements are absolutely
determined in all possible worlds.

A class of {neutrosophic set} in which every element x has the
form $x(1^+, ^-0, ^-0)$, i.e. absolutely belongs to the set;
here T,I,F are real standard or non-standard subsets, included
in the non-standard unit interval $]^-0, 1^+[$, representing truth,
indeterminacy, and falsity percentages respectively.

Contrast with {nihilist set} and {nihilism}.

Related to {tautologism}.

========================================================

**1.8. Nihilist Set**:

<logic, mathematics> A set whose elements absolutely
don't belong to the set in all possible worlds.

A class of {neutrosophic set} in which every element x has the
form $x(^-0, ^-0, 1^+)$, i.e. absolutely doesn't belongs to the set;
here T,I,F are real standard or non-standard subsets, included
in the non-standard unit interval $]^-0, 1^+[$, representing truth,
indeterminacy, and falsity percentages respectively.

The empty set is a particular set of {nihilist set}.

Contrast with {tautological set}.

Related to {nihilism}.

====================================================

**1.9. Dialetheist Set**:

 <logic, mathematics> /di:-al-u-theist/ A set which contains at least one element which also belongs to its complement.

 A class of {neutrosophic set} which models a situation where the intersection of some disjoint sets is not empty.

 There is at least one element x(T,I,F) of the dialetheist set M which belongs at the same time to M and to the set C(M), which is the complement of M;
 here T,I,F are real standard or non-standard subsets, included in the non-standard unit interval $]^-0, 1^+[$, representing truth, indeterminacy, and falsity percentages respectively.

 Contrast with {trivialist set}.

 Related to {dialetheism}.

 ====================================================

**1.10. Trivialist Set**:

 <logic, mathematics> A set all of whose elements also belong to its complement.

 A class of {neutrosophic set} which models a situation where the intersection of any disjoint sets is not empty.

 Every element x(T,I,F) of the trivialist set M belongs at the same time to M and to the set C(M), which is the complement of M;
 here T,I,F are real standard or non-standard subsets, included in the non-standard unit interval $]^-0, 1^+[$, representing truth, indeterminacy, and falsity percentages respectively.

 Contrast with {dialetheist set}.

 Related to {trivialism}.

 ====================================================

## 2. Definitions of New Probabilities and Statistics

=======================================================

**2.1. Neutrosophic Probability**:

<probability> The probability that an event occurs is (T, I, F), where T,I,F are real standard or non-standard subsets, included in the non-standard unit interval $]^-0, 1^+[$, representing truth, indeterminacy, and falsity percentages respectively.

Therefore: $^-0 \leq \inf(T) + \inf(I) + \inf(F) \leq \sup(T) + \sup(I) + \sup(F) \leq 3^+$.

Generalization of {classical probability} and {imprecise probability}, {intuitionistic probability}, {paraconsistent probability}, {faillibilist probability}, {paradoxist probability}, {tautological probability}, {nihilistic probability}, {dialetheist probability}, {trivialist probability}.

Related with {neutrosophic set} and {neutrosophic logic}.

The analysis of neutrosophic events is called **Neutrosophic Statistics**.

=======================================================

**2.2. Intuitionistic Probability**:

<probability> The probability that an event occurs is (T, I, F), where T,I,F are real standard or non-standard subsets, included in the non-standard unit interval $]^{-0, 1+}[$, representing truth, indeterminacy, and falsity percentages respectively, and $n\_sup = \sup(T)+\sup(I)+\sup(F) < 1$, i.e. the probability is incompletely calculated.

Contrast with {paraconsistent probability}.

Related to {intuitionistic set} and {intuitionistic logic}.

The analysis of intuitionistic events is called **Intuitionistic Statistics**.

=======================================================

**2.3. Paraconsistent Probability**:

<probability> The probability that an event occurs is (T, I, F),
where T,I,F are real standard or non-standard subsets, included in the
non-standard unit interval $]^-0, 1^+[$, representing truth,
indeterminacy, and falsity percentages respectively,
and n_sup = sup(T)+sup(I)+sup(F) > 1,
i.e. contradictory information from various sources.

Contrast with {intuitionistic probability}.

Related to {paraconsistent set} and {paraconsistent logic}.

The analysis of paraconsistent events is called
**Paraconsistent Statistics**.

========================================================

**2.4. Faillibilist Probability**:

<probability> The probability that an event occurs is (T, I, F),
where T,I,F are real standard or non-standard subsets, included in the
non-standard unit interval $]^-0, 1^+[$, representing truth,
indeterminacy, and falsity percentages respectively,
and inf(I) > 0,
i.e. there is some percentage of indeterminacy in calculation.

Related to {faillibilist set} and {faillibilism}.

The analysis of faillibilist events is called **Faillibilist Statistics**.

========================================================

**2.5. Paradoxist Probability**:

<probability> The probability that an event occurs is (1, I, 1),
where I is a standard or non-standard subset, included in the
non-standard unit interval $]^-0, 1^+[$, representing indeterminacy.

Paradoxist probability is used for paradoxal events (i.e. which
may occur and may not occur simultaneously).

Related to {paradoxist set} and {paradoxism}.

The analysis of paradoxist events is called **Paradoxist Statistics**.

========================================================

**2.6. Pseudo-Paradoxist Probability**:

<probability> The probability that an event occurs is either (1, I, F) with
0<inf(F)≤sup(F)<1, or (T, I, 1) with 0<inf(T)≤sup(T)<1,
where T,I,F are standard or non-standard subset, included in the
non-standard unit interval $]^{-}0, 1^{+}[$, representing the truth, indeterminacy, and
falsity percentages respectively.

Pseudo-Paradoxist probability is used for pseudo-paradoxal events (i.e. which
may certainly occur and may not partially occur simultaneously,
or may partially occur and may not certainly occur simultaneously).

Related to {pseudo-paradoxist set} and {pseudo-paradoxism}.

The analysis of pseudo-paradoxist events is called **Pseudo-Paradoxist Statistics**.

========================================================

**2.7. Tautological Probability**:

<probability> The probability that an event occurs is more than one,
i.e. $(1^{+}, {}^{-}0, {}^{-}0)$.

Tautological probability is used for universally sure events (in all
possible worlds, i.e. do not depend on time, space, subjectivity, etc.).

Contrast with {nihilistic probability} and {nihilism}.

Related to {tautological set} and {tautologism}.

The analysis of tautological events is called **Tautological Statistics**.

========================================================

**2.8. Nihilist Probability**:

<probability> The probability that an event occurs is less than zero,
i.e. $({}^{-}0, {}^{-}0, 1^{+})$.

Nihilist probability is used for universally impossible events (in all
possible worlds, i.e. do not depend on time, space, subjectivity, etc.).

Contrast with {tautological probability} and {tautologism}.

Related to {nihilist set} and {nihilism}.

The analysis of nihilist events is called **Nihilist Statistics**.

========================================================

**2.9. Dialetheist Probability**:

<probability> /di:-al-u-theist/ A probability space where at least one event and its complement are not disjoint.

A class of {neutrosophic probability} that models a situation where the intersection of some disjoint events is not empty.

Here, similarly, the probability of an event to occur is (T, I, F), where T,I,F are real standard or non-standard subsets, included in the non-standard unit interval $]^-0, 1^+[$, representing truth, indeterminacy, and falsity percentages respectively.

Contrast with {trivialist probability}.

Related to {dialetheist set} and {dialetheism}.

The analysis of dialetheist events is called **Dialetheist Statistics**.

========================================================

**2.10. Trivialist Probability**:

<probability> A probability space where every event and its complement are not disjoint.

A class of {neutrosophic probability} which models a situation where the intersection of any disjoint events is not empty.

Here, similarly, the probability of an event to occur is (T, I, F), where T,I,F are real standard or non-standard subsets, included in the non-standard unit interval $]^-0, 1^+[$, representing truth, indeterminacy, and falsity percentages respectively.

Contrast with {dialetheist probability}.

Related to {trivialist set} and {trivialism}.

The analysis of trivialist events is called **Trivialist Statistics**.

===========================================================

### 3. Definitions of New Logics

===========================================================

This definition is not quite new but because the next ones are connected to it we recall it:

**3.1. Neutrosophic Logic**:

<logic, mathematics> A logic which generalizes many existing classes of logics, especially the fuzzy logic.

In this logic each proposition is estimated to have the percentage of truth in a subset T, the percentage of indeterminacy in a subset I, and the percentage of falsity in a subset F;
here T,I,F are real standard or non-standard subsets, included in the non-standard unit interval $]^-0, 1^+[$, representing truth, indeterminacy, and falsity percentages respectively.

Therefore: $^-0 \leq \inf(T) + \inf(I) + \inf(F) \leq \sup(T) + \sup(I) + \sup(F) \leq 3^+$.

Generalization of {classical or Boolean logic}, {fuzzy logic}, {multiple-valued logic}, {intuitionistic logic}, {paraconsistent logic}, {faillibilist logic, or failibilism}, {paradoxist logic, or paradoxism}, {pseudo-paradoxist logic, or pseudo-paradoxism}, {tautological logic, or tautologism}, {nihilist logic, or nihilism}, {dialetheist logic, or dialetheism}, {trivialist logic, or trivialism}.

Related to {neutrosophic set}.

========================================================

**3.2. Paradoxist Logic (or Paradoxism)**:

<logic, mathematics> A logic devoted to paradoxes, in which each proposition has the logical vector value (1, I, 1);
here I is a real standard or non-standard subset, included in the non-standard unit interval $]^-0, 1^+[$, representing the indeterminacy.

As seen, each paradoxist (paradoxal) proposition is true and false simultaneously.

Related to {paradoxist set}.

===========================================================

### 3.3. Pseudo-Paradoxist Logic (or Pseudo-Paradoxism):

<logic, mathematics> A logic devoted to pseudo-paradoxes,
in which each proposition has the logical vector value:
either (1, I, F), with $0<\inf(F)\leq\sup(F)<1$,
or (T, I, 1), with $0<\inf(T)\leq\sup(T)<1$;
here I is a real standard or non-standard subset, included in the
non-standard unit interval $]^-0, 1^+[$, representing the indeterminacy.

As seen, each pseudo-paradoxist (pseudo-paradoxal) proposition is:
either totally true and partially false simultaneously,
or partially true and totally false simultaneously.

Related to {pseudo-paradoxist set}.

=====================================================

### 3.4. Tautological Logic (or Tautologism):

<logic, mathematics> A logic devoted to tautologies, in which each
proposition has the logical vector value $(1^+, {}^-0, {}^-0)$.

As seen, each tautological proposition is absolutely true (i. e, true in all
possible worlds).

Related to {tautological set}.

=====================================================

General References:

[1]  Jean Dezert, *Open Questions on Neutrosophic Inference*, Multiple-Valued Logic Journal, 2001 (to appear).
[2]  Denis Howe, *On-Line Dictionary of Computing*, http://foldoc.doc.ic.ac.uk/foldoc/
[3]  Charles Le, *Preamble to Neutrosophy and Neutrosophic Logic*, Multiple-Valued Logic Journal, 2001 (to appear).
[4]  Florentin Smarandache, *A Unifying Field in Logics: Neutrosophic Logic*, Multiple-Valued Logic Journal, 2002 (to appear).
[5]  Florentin Smarandache, organizer, *First International Conference on Neutrosophy, Neutrosophic Probability, Set, and Logic*, University of New Mexico, 1-3 December 2001; http://www.gallup.unm.edu/~smarandache/FirstNeutConf.htm
[6]  Florentin Smarandache, *Neutrosophy, a New Branch of Philosophy*, Multiple-Valued Logic Journal, 2002 (to appear).

# A Short Note on Financial Data Set Detection using Neutrosophic Probability

Dr. Mohammad Khoshnevisan
School of Accounting and Finance
Griffith University, PMB 50 Gold Coast, QLD 9726 Australia
and
Sukanto Bhattacharya
School of Information Technology
Bond University, Gold Coast, QLD 4229, Australia

**Abstract:**
This study actually draws from and builds on "Benford's law and its application in financial misrepresented of a data " (Kumar and Bhattacharya, 2001). here we have simply added a *neutrosophic dimension* to the problem of determining the *conditional probability* that a financial misrepresentation of the data set, has been actually committed, given that no Type I error occurred while rejecting the null hypothesis $H_0$: the observed first-digit frequencies approximate a benford distribution; and accepting the alternative hypothesis $H_1$: the observed first-digit frequencies do not approximate a benford distribution.

**Keywords**: Financial misrepresented data set, Benford's law, probability distributions, neutrosophic probability

**Introduction:**

**1. Testing for manipulation in a set of accounting data.**

Kumar and Bhattacharya, 2001, proposed a Monte Carlo adaptation of Benford's law. There has been some research already on the application of Benford's law. However, most of the practical work in this regard has been concentrated in detecting the first digit frequencies from the account balances selected on basis of some known audit sampling method and then directly comparing the result with the expected Benford frequencies. We have voiced slight reservations about this technique in so far as that the Benford frequencies are necessarily **steady state frequencies** and may not therefore be truly reflected in the sample frequencies. As samples are always of finite sizes, it is therefore perhaps not entirely fair to arrive at any conclusion on the basis of such a direct comparison, as the sample frequencies won't be steady state frequencies.

However, if we draw digits randomly using the **inverse transformation technique** from within random number ranges derived from a cumulative probability distribution function based on the Benford frequencies; then the problem boils down to running a *goodness of fit* kind of test to identify any significant difference between observed and simulated first-digit frequencies. This test may be conducted using a known sampling distribution like for example the **Pearson's $\chi^2$ distribution**. The random number ranges for the Monte Carlo simulation are to be drawn from a cumulative probability distribution function based on the following Benford probabilities given in Table I.

Table I

| First Significant Digit | 1 | 2 | 3 | 4 | 5 | 6 | 7 | 8 | 9 |
|---|---|---|---|---|---|---|---|---|---|
| Benford Probability | 0.301 | 0.176 | 0.125 | 0.097 | 0.079 | 0.067 | 0.058 | 0.051 | 0.046 |

The first-digit probabilities can be best approximated mathematically by the log-based formula as was derived by Benford: P (First significant digit = d) = $\log_{10} [1 + (1/d)]$.

## 2. Computational Algorithm:

Define a finite sample size n and draw a sample from the relevant account balances using a suitable audit sampling procedure

1. Perform a continuous Monte Carlo run of length $\lambda^* \approx (1/2\varepsilon)^{2/3}$ grouped in epochs of size n using a customized MS-Excel spreadsheet.

2. Test for significant difference in sample frequencies between the first digits observed in the sample and those generated by the Monte Carlo simulation by using a "goodness of fit" test using the $\chi^2$ distribution. The null and alternative hypotheses are as follows:

**$H_0$:** The observed first digit frequencies approximate a Benford distribution

**$H_1$:** The observed first digit frequencies do not approximate a Benford distribution

This statistical test will not reveal whether or not a data misrepresentation has actually been committed. All it does is establish at a desired level of confidence, that the accounting data has not been manipulated (if H0 cannot be rejected).

However, given that $H_1$ is accepted and $H_0$ is rejected, it could imply any of the following events:

I. There is no manipulation - occurrence of a Type I error i.e. $H_0$ rejected when true.

II. There is manipulation *and* such manipulation *is definitely* misrepresented.

III. There is manipulation *and* such manipulation *may or may not be* misrepresented

IV. There is manipulation *and* such manipulation *is definitely not* misrepresented

**3. Neutrosophic Extension:**

Neutrosophic probabilities are a generalization of classical and fuzzy probabilities and cover those events that involve some degree of indeterminacy. It provides a better approach to quantifying uncertainty than classical or even fuzzy probability theory.

Neutrosophic probability theory uses a subset-approximation for truth-value as well as indeterminacy and falsity values. Also, this approach makes a distinction between "relative true event" and "absolute true event" the former being true in only some probability sub-spaces while the latter being true in all probability sub-spaces. Similarly, events that are false in only some probability sub-spaces are classified as "relative false events" while events that are false in all probability sub-spaces are classified as "absolute false events". Again, the events that may be hard to classify as either 'true' or 'false' in some probability sub-spaces are classified as "relative indeterminate events" while events that bear this characteristic over all probability sub-spaces are classified as "absolute indeterminate events".

While in classical probability **$n\_sup \leq 1$**, in neutrosophic probability **$n\_sup \leq 3^+$** where $n\_sup$ is the upper bound of the probability space. In cases where the truth and falsity components are complimentary, i.e. there is no indeterminacy, the components sum to unity and neutrosophic probability is reduced to classical probability as in the tossing of a fair coin or the drawing of a card from a well-shuffled deck.

Coming back to our original problem of financial misrepresented of a data set, let E be the event whereby a Type I error has occurred and F be the event whereby a misrepresented set is actually detected. Then the *conditional neutrosophic probability* **NP (F | $E^c$)** is defined over a probability space consisting of a triple of sets (T, I, U). Here, T, I and U are probability sub-spaces wherein event F is t% true, i% indeterminate and u% untrue respectively, given that no Type I error occurred.

The sub-space T within which t varies may be determined by factors such as past records of the misrepresented data set in the organization and effectiveness of internal control systems. On the other hand, the sub-space U may be determined by factors like personal track records of the employees in question, the position enjoyed and the remuneration drawn by those employees. For example, if the magnitude of the embezzled amount is deemed too frivolous with respect to the position and remuneration of the employees involved. The sub-space I is most likely to be determined by the mutual inconsistency that might arise between the effects of some of the factors determining T and U.

### 4. Conclusion:

No doubt then that the theory of neutrosophic probability opens up a new vista of analytical reasoning for the techno-savvy forensic accountant. In this paper, we have only posit that a combination of statistical testing of audit samples based on Benford's law combined with a neutrosophic reasoning could help the forensic accountant in getting a better fix on the quantitative possibility of actually dealing with misrepresented financial data set. This is an emerging science and thus holds a vast potential of future research endeavours the ultimate objective of which will be to actually come up with a reliable, comprehensive computational methodology to track down the misrepresented financial data. We believe our present effort is only one initial step in that direction.

# Intentionally and Unintentionally.
# On Both, A and Non-A, in Neutrosophy

Feng Liu
Department of Management Science and Engineering
Shaanxi Economics and Trade Institute (South Campus)
South Cuihua Road, Xi'an, Shaanxi, 710061, P. R. China
E-mail: liufeng49@sina.com

Florentin Smarandache
Department of Mathematics
University of New Mexico, Gallup, NM 87301, USA
E-mail: smarand@unm.edu

**Abstract:** The paper presents a fresh new start on the neutrality of neutrosophy in that "both A and Non-A" as an alternative to describe Neuter-A in that we conceptualize things in both intentional and unintentional background. This unity of opposites constitutes both objective world and subjective world. The whole induction of such argument is based on the intensive study on Buddhism and Daoism including I-ching. In addition, a framework of contradiction oriented learning philosophy inspired from the Later Trigrams of King Wen in I-ching is meanwhile presented. It is shown that although A and Non-A are logically inconsistent, but they are philosophically consistent in the sense that Non-A can be the unintentionally instead of negation that leads to confusion. It is also shown that Buddhism and Daoism play an important role in neutrosophy, and should be extended in the way of neutrosophy to all sciences according to the original intention of neutrosophy.

## 2. Objective world and subjective world

The common confusion about the objective world is: it **is** just what we see and feel. This is however very wrong. In fact, this is rather a belief than an objective reflection, and varies among different people, because none of us can prove it. In his paper "To be or not to be, A multidimensional logic approach" Carlos Gershenson [2] has generalized proofs:
- Everything is and isn't at a certain degree. (i.e., there is no absolute truth or false);
- Nothing can be proved (that it exists or doesn't) (i.e., no one can prove whether his consciousness is right);
- I believe, therefore I am (i.e., I take it true, because I believe so).

It is something, but not that figured in our mind. This is the starting point in Daoism (Liu [2]).

Daodejing begins with: "Dao, daoable, but not the normal dao; name, namable, but not the normal name." We can say it is dao, but it doesn't mean what we say. Whenever we mention it, it is beyond the original sense.

Daodejing mainly deals with the common problem: "What/who creates everything in the world we see and feel?" It is dao: like a mother that bears things with shape and form. But what/who is dao? It is just unimaginable, because whenever we imagine it, our imagination can never be it (we can never completely describe it: more we describe it, more wrong we are). It is also unnamable, because whenever we name it, our concept based on the name can never be it.

Daoism illustrates the origin of everything as such a form that doesn't show in any form we can perceive. This is the reason why it says, everything comes from nothingness, or this nothingness creates everything in forms in dynamic change. Whatever we can perceive is merely the created forms, rather than its genuine nature, as if we distinguish people by their outer clothes. We are too far from understanding the nature, even for the most prominent figures like Einstein.

Therefore **Name and Non-Name coexist** pertaining to an object:

Object = both Name and Non-Name.

Then what should we do **subjectively**? Very simple: both intentionally and unintentionally. Intentional conception relates to all the connotation and extension pertaining to Name, and unintentional one to Non-Name.

- **There are alternative interpretations on Non-Name: unintentionally and negatively. This is crucial in our confusion.**

This is the contradiction between creativity and implementation, as is stated below.

3. **Neutrosophy**

Neutrosophy is a new branch of philosophy that studies the origin, nature, and scope of neutralities, as well as their interactions with different ideational spectra.
It is the base of *neutrosophic logic*, a multiple value logic that generalizes the fuzzy logic and deals with paradoxes, contradictions, antitheses, antinomies.

**Characteristics** of this mode of thinking:
- proposes new philosophical theses, principles, laws, methods, formulas, movements;
- reveals that world is full of indeterminacy;
- interprets the uninterpretable;
- regards, from many different angles, old concepts, systems:
showing that an idea, which is true in a given referential system, may be false in another one, and vice versa;
- attempts to make peace in the war of ideas,
and to make war in the peaceful ideas;
- measures the stability of unstable systems,
and instability of stable systems.

Let's note by <A> an idea, or proposition, theory, event, concept, entity, by <Non-A> what is not <A>, and by <Anti-A> the opposite of <A>. Also, <Neut-A> means what is neither <A> nor <Anti-A>, i.e. neutrality in between the two extremes. And <A'> a version of <A>.

<Non-A> is different from <Anti-A>.

**Main Principle:**
Between an idea <A> and its opposite <Anti-A>, there is a continuum-power spectrum of neutralities <Neut-A>.

**Fundamental Thesis of Neutrosophy:**
Any idea <A> is T% true, I% indeterminate, and F% false, where T, I, F $\subset$ ] $^-$0, 1$^+$ [.

**Main Laws of Neutrosophy:**
Let <α> be an attribute, and (T, I, F) $\subset$ ] $^-$0, 1$^+$ [$^3$. Then:
- There is a proposition <P> and a referential system {R}, such that <P> is T% <α>, I% indeterminate or <Neut-α>, and F% <Anti-α>.
- For any proposition <P>, there is a referential system {R}, such that <P> is T% <α>, I% indeterminate or <Neut-α>, and F% <Anti-α>.
- <α> is at some degree <Anti-α>, while <Anti-α> is at some degree <α>.

## 4. Creativity and implementation

We can model our mind in the alternation of yin and yang that is universal in everything (Feng Liu):
- Yang pertains to dynamic change, and directs great beginnings of things; yin to relatively static stage, and gives those exhibited by yang to their completion.
    In the course of development and evolution of everything yang acts as the creativity (Feng Liu) that brings new beginnings to it, whereas yin implements it in forms as we perceive as temporary states. It is in this infinite parallelism things inherit modifications and adapt to changes.
    - On our genuine intelligence — **creativity** (Liu [4])

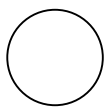

In the query about the figure on the left, whenever we hold the answer as a circle, we are inhibiting our creativity. Nor should we hold that it is a cake, a dish, a bowl, a balloon, or the moon, the sun, for we also spoil our creativity in this way. Then, what is it?

"It is nothing."

Is it correct? It is, if we do not hold on to the assumption "it is something". It is also wrong, if we persist in the doctrine "the figure is something we call nothing." This nothing has in this way become something that inhibits our creativity. How ridiculous!

**Whenever we hold the belief "it is …", we are loosing our creativity. Whenever we hold that "it is not …", we are also loosing our creativity.** Our true intelligence requires that we completely free our mind — neither stick to any extremity nor to "no sticking to any assumption or belief". This is a kind of genius or gift rather than logic rules, acquired largely after birth, e.g., through

Buddhism practice. **Note that our creativity lies just between internationality and uninternationality.**

**Not (it is) and not (it is not),**
**It seems nothing, but creates everything,**
**Including our true consciousness,**
**The power of genius to understand all.**

- The further insight on contradiction compatible **learning philosophy** inspired from the Later Trigrams of King Wen of I-ching shows that:

    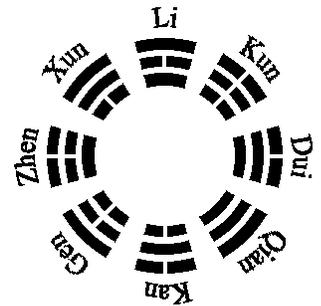

    When something (controversial) is perceived (in Zhen), it is referred (in Xun) to various knowledge models and, by assembling the fragments perceived from these models, we reach a general pattern to which fragments attach (in Li), as hypothesis, which needs to be nurtured and to grow up (Kun) in a particular environment. When the hypothesis is mature enough, it needs to be represented (in Dui) in diverse situations, and to expand and contradict with older knowledge (in Qian) to make update, renovation, reformation or even revolution in knowledge base, and in this way the new thought is verified, modified and substantialized. When the novel thought takes the principal role (dominant position) in the conflict, we should have a rest (in Kan) to avoid being trapped into depth (it would be too partial of us to persist in any kind of logic, to adapt to the outer changes). Finally the end of cycle (in Gen).

    I-ching [in Chinese: Yi Jing] means: *Yi* = change, *Jing* = scripture. It mainly deals with the creation, evolution (up and down) of everything in such perspective that everything is an outer form of a void existence, and that everything always exists in the form of unity (compensation, complementation) of opposites…

    This philosophy shows that **contradiction acts as the momentum or impetus to learning evolution. No controversy, no innovation. This is the essentially of neutrosophy** (Florentin Smarandache).

    In the cycle there is unintentionally implied throughout it:
    - Where do the reference models relating to the present default model come from? They are different objectively.
    - How can we assemble the model from different or even incoherent or inconsistent fragments?
    - If we always do it intentionally, how does the hypothesis grow on its own, as if we study something without sleep?
    - How can our absolute intention be complemented without contradiction?
    - Is it right that we always hold our intention?

> There is only one step between truth and prejudice — when the truth is overbelieved regardless of constraint in situations, it becomes prejudice.

- Is there no end for the intention? Then, how can we obtain a concept that is never finished? If there is an end, then it should be the beginning of unintentionally, as yin and yang in Taiji figure.

  > Unintentionally can be alternately derived when the intention is repeated over and over so that it becomes an instinct. **It is even severe that we develop a fallacy instinct.**

## 5. Completeness and incompleteness: knowledge and practice

There always is contradiction between completeness and incompleteness of knowledge. In various papers presented by Carlos Gershenson he proves this point. Same in Daoism and Buddhism. This contradiction is shown in the following aspect:
  a) People are satisfied with their knowledge relative to a default, well-defined domain. But later on, they get fresh insight in it.
  b) They face with contradictions and new challenges in their practice and further development.

This reflects in our weakness that:
- Do we understand ourselves?
- Do we understand the universe?

What do we mean by knowledge, complete whole or incomplete? Our silliness prompts us to try the complete specifications, but where on earth are they (Gershenson [1])? Meanwhile, our effort would be nothing more than a static imitation of some dynamic process (Liu [3]), since **human understands the world through the interaction of the inter-contradictory and inter-complementary two kinds of knowledge: perceptual knowledge and rational knowledge - they can't be split apart.**

- In knowledge discovering, there are merely strictly limited condition that focus our eyes to a local domain rather than a open extension, therefore our firsthand knowledge is only relative to our default referential system, and extremely subjective possibly.
- Is it possible to reach a relatively complete piece at first? No, unless we were gods (we were objective in nature).
- Then we need to perceive (the rightness, falseness, flexibility, limitation, more realistic conception, etc.) and understand the real meaning of our previous knowledge — how: only through practice (how can we comprehend the word "apple" without tasting it?)
- Having done that, we may have less subjective minds, based on which the original version is modified, revised, and adapted as further proposals.
- Again through practice, the proposals are verified and improved.
- This cycle recurs to the infinite, in each of which our practice is extended in a more comprehensive way; the same to our knowledge.
  - Discover the truth through practice, and again through practice verify and develop the truth. Start from perceptual knowledge and actively develop it into rational knowledge; then start from rational knowledge and actively guide revolutionary practice to change both the subjective and the objective world. Practice, knowledge, again practice, and again knowledge. This form repeats itself in endless cycles, and with each cycle the content of practice and knowledge rises to a higher level.
- Through practice, we can
  - verify our knowledge
  - find the inconsistency, incompleteness of our knowledge, and face new problems, new challenges as well
  - maintain critical thought

Therefore **knowledge is based on the infinite critics and negation (partial or revolutionary) on our subjective world.** It is never too old to learn.

## 6. Conclusion

Whenever we say "it is", we refer it to both subjective and objective worlds.

We can creatively use the philosophical expression **both A and non-A** to describe both subjective world and objective worlds, and possibly the neutrality of both.

Whenever there is "it is", there is subjective world, in the sense that concepts always include subjectivity. So our problem becomes: is "it" really "it"? A real story of Chinese Tang dynasty recorded in a sutra (adapted from Yan Kuanhu Culture and Education Fund)shows that:

> Huineng arrived at a Temple in Guangzhou where a pennant was being blown by wind. Two monks who happened to see the pennant were debating what was in motion, the wind or the pennant.
> Huineng heard their discussion and said: "It was neither the wind nor the pennant. What actually moved were your own minds." Overhearing this conversation, the assembly (a lecture was to begin) were startled at Huineng's knowledge and outstanding views.

- When we see pennant and wind we will naturally believe we are right in our consciousness, however it is subjective. In other words, what we call "the objective world" can never absolutely be objective at all.
- Whenever we believe we are objective, this belief however is subjective too.
- In fact, all these things are merely our mental creations (called illusions in Buddhism) that in turn cheat our consciousness: There is neither pennant nor wind, but our mental creations.
- The world is made up of our subjective beliefs that in turn cheat our consciousness. This is in fact a cumulative cause-effect phenomenon.
- Everyone can extricate himself out of this maze, said Sakyamuni and all the Buddhas, Bodhisattvas around the universe, their number is as many as that of the sands in the Ganges (Limitless Life Sutra).

# Logic: a Misleading Concept.
# A Contradiction Study toward Agent's Logic Ontology

Feng Liu
Department of Management Science and Engineering
Shaanxi Economics and Trade Institute (South Campus)
South Cuihua Road, Xi'an, Shaanxi, 710061, P. R. China
E-mail: liufeng49@sina.com

Florentin Smarandache
Department of Mathematics
University of New Mexico, Gallup, NM 87301, USA
E-mail: smarand@unm.edu

**Abstract:** The paper presents a fresh new comprehensive ideology on Neutrosophic Logic based on contradiction study in a broad sense: general critics on conventional logic by examining the essence of logic, fresh insights on logic definition based on Chinese philosophical survey, and a novel and genetic logic model as the elementary cell against Von Neumann oriented ones based on this novel definition. As for the logic definition, the paper illustrates that logic is rather a tradeoff between different factors than truth and false abstraction. It is stressed that the kernel of any intelligent system is exactly a contradiction model. The paper aims to solve the chaos of logic and exhibit the potential power of neutrosophy: a new branch of scientific philosophy.

1. Background

Although it is commonly believed that intelligence is a social activity, and it is therefore represented in multiagent forms, but its kernel, the logic of agents, remains controversial with its static, monolateral or homogeneous forms.

This reflects in their behaviors as: our agents appear social in outer forms but autarchic in nature, for this kind of multiagent system can never deal with controversies, critics, conflicts or something with flexibility. Our multiagent system has become a sort of software engineering or system engineering of fresh forms, failing to implement our presumed social intelligence.

In the long-term exploration, one realizes that the problem takes its root in the misleading definition of logic. Even the simplest logic such as "The earth turns around the sun" and "I'll visit him if it doesn't rain and he is in" can lead to ambiguous or contradictory actions of agent (Liu [7]). Limited to the length, I'll present in this paper only a framework to launch our discussion, as follows:
- Fact: a belief rather than truth

- Logic: dependent of situations, not absolute
- Logic is negating itself
- Logic is only one perspective of learning, not an independent entity
- As a part of learning, logic is dynamic
- As a part of learning, logic is multilateral
- Logic is always partial
- Illusion and creativity

Many scientists argue about the need to model human intelligence in the general level. The argument lies in our vague understanding of intelligent system (Liu [6]). Intelligent system should be, in our opinion, **a tradeoff machine in order to adapt to its environment**. Then a specific model becomes such a tradeoff between ideal philosophic model and practical system model, in the hierarchy from philosophic layer down to a specific application or situation constraint implementation. I'll show this philosophy in step ward way, as follows.

2. **Neutrosophy**

Neutrosophy is a new branch of philosophy that studies the origin, nature, and scope of neutralities, as well as their interactions with different ideational spectra.
It is the base of *neutrosophic logic*, a multiple value logic that generalizes the fuzzy logic and deals with paradoxes, contradictions, antitheses, antinomies.

**Characteristics** of this mode of thinking:
- proposes new philosophical theses, principles, laws, methods, formulas, movements;
- reveals that world is full of indeterminacy;
- interprets the uninterpretable;
- regards, from many different angles, old concepts, systems: showing that an idea, which is true in a given referential system, may be false in another one, and vice versa;
- attempts to make peace in the war of ideas, and to make war in the peaceful ideas;
- measures the stability of unstable systems, and instability of stable systems.

Let's note by <A> an idea, or proposition, theory, event, concept, entity, by <Non-A> what is not <A>, and by <Anti-A> the opposite of <A>. Also, <Neut-A> means what is neither <A> nor <Anti-A>, i.e. neutrality in between the two extremes. And <A'> a version of <A>.

<Non-A> is different from <Anti-A>.

**Main Principle:**
Between an idea <A> and its opposite <Anti-A>, there is a continuum-power spectrum of neutralities <Neut-A>.

**Fundamental Thesis of Neutrosophy:**
Any idea <A> is T% true, I% indeterminate, and F% false, where T, I, F $\subset$ ] $^-$0, $1^+$ [. Here ] $^-$0, $1^+$ [ is a non-standard unit interval, with $^-$0={0-ε, ε is a positive infinitesimal number} and $1^+$={1+ε, ε is a positive infinitesimal number}.

**Main Laws of Neutrosophy:**

Let <α> be an attribute, and (T, I, F) ⊂ ] $^-$0, 1$^+$ [$^3$. Then:
- There is a proposition <P> and a referential system {R}, such that <P> is T% <α>, I% indeterminate or <Neut-α>, and F% <Anti-α>.
- For any proposition <P>, there is a referential system {R}, such that <P> is T% <α>, I% indeterminate or <Neut-α>, and F% <Anti-α>.
- <α> is at some degree <Anti-α>, while <Anti-α> is at some degree <α>.

3. **Fact: a Belief rather than Truth**

We start with an ancient problem based on the following contradiction:
- The sun turns around the earth.
- The earth turns around the sun.

Of cause nearly everyone of us would answer: the later is absolute right. Note that this is merely a belief, because in Copernicus's age the majority believed in the former. Has anyone proved nowadays whether the former is incorrect? If yes, he must have assumed that the sun is relatively fixed. Unfortunately this is also his belief, because none of us has ever proved the absoluteness of his consciousness: when we see something, is it really something or just we believe that there is something (we really touch something or we really believe it is something we touched)? Or more specifically, is it an object or just we hold long this same belief? Do we really exist as in form we see or just we believe so? I have to introduce a heard experiment to show this point.

> A blindfold person is told to be experimented with an iron burnt hot. And through a chronic preparation before him, the iron is burnt fervid, and he is told that the iron is gradually moved closer and closer to him.
> "Yes, I am feeling hotter and hotter, …… really hot, extremely, ……"
> The gradual process goes on and on, until suddenly, he is instructed to have his skin burnt.
> "Oh……", his skin really burnt.
> When he opened his eyes, there is nothing but the scorch in him—there is no fire nor iron, but merely his **imagination—it is strong enough to cause the effect.**

I experienced another experiment in which four of us were pointing to a carefully set small wooden stool while rotating around it. According our mutual will, the stool turned itself in the same direction we turned!

Another fact (shown in a qigong journal quite a number of years ago, the following is based on our memory) shows the same thing:

> There is a qigong (commonly believed as some mental or physical exercise in order to gather the "energy" from nature, qi (there are such a kind of substance in Chinese medicine which is unseen but really affects our body), or the concentrative power of will to maintain health from disease, it is not a feasible way to us) expert in China who, through chronic practice, can "brake" a steel saw blade with nothing but his will, and he had been succeeding in it nearly every time, even in many qigong reports.
> Once he re-showed the same talent to the huge audience with great curiosity. He ordered: "break", but unexpectedly, the blade remain exact the same as before, and the following tries turned out to be the same failures. The atmosphere became extremely unfavorable.

Fortunately however, the chairman of the qigong report is experienced, and asked the audience to cooperate: more you are confident, more successful the experiment.

Magically, the expert broke the blade with a single command.

Conclusion:

We have to confess from the three experiments that **fact is really our belief**—if a single belief is not powerful enough to convince us, the mutual belief, especially of all the human beings, is definitely strong enough to illude ourselves. While this cause-effect goes on and on, we are unconsciously trapped in the inextricable web of deceit designed by ourselves.

Only wise man can see through this kind of deceit, e.g., Master Huineng in Chinese Tang dynasty when he saw an argument about a pennant aflutter: whether the wind was moving or the pennant.

**"Neither. What actually moved were your own minds."** (Liu [9], see also Yan Kuanhu [3])

Everyone can become wise when understands this cause-effect, which is the basic point of Buddhism (Chin Kung [1])

4. **Logic: Dependent of Situations, not Absolute**

Take the logic 1+1=2 for example. Is it correct? Consider

black+white=?, explosive+fire=?, warm+cold=?, theory+practice=?, and
yin+yang = ?

Does the idiom "Blind People Touching an Elephant" really refer to blind men and elephant?

- Blind People Touching an Elephant, a story from the Mahapra Janaparamita Sutra:

    The story shows that **the same elephant can be interpreted as such different things** as turnip, dustpan, pestle, bed, jar and rope **by different blind people** who touch it in turn. The first one touches the tusk, the second the ear, the third the foot, the fourth the back, the fifth the belly, and the last the tail.

    **Based on their different beliefs, the same elephant conveys diverse logics**.

Conclusion:

Logic is more a kind of mental behavior than an objective understanding, i.e., it is more a belief rather than truth; this belief is based on "facts" which are also beliefs.

This belief is subject to dynamic changes with situations, and more general belief (general understanding) relative to more general situations could be too flexible to grasp (e.g., Dao=yin+yang), therefore **logic suggests varying explanations based on particularity of situations.**

5. **Logic is Negating Itself**

Logic comes as mental reflection and leads to new reflection. So, **it is not the problem of logic (validity) but the ways we reflect it**, otherwise it would become the Chinese room experiment.
- J. R. Searle shows this "Chinese room problem" in his paper "Minds, brains and programs":

    We set an Englishman which does not know Chinese, in a closed room, with many symbols of the Chinese language, and a book of instructions in English of how to manipulate the symbols when a set of symbols (instructions) is given. So, Chinese scientists will give him instructions in Chinese, and the **Englishman will manipulate symbols in Chinese, and he will give a correct answer in Chinese. But he is not conscious of what he did.** We suppose that a machine behaves in a similar way: it might give correct answers, but it is not conscious of what it is doing.

Another argument on validity of logic is based on a Chinese idiom: Cutting a Mark on the Boatside to Retrieve a Sword (Young):

    Once, a man of the State of Chu (ancient China) took a boat to cross a river. It so happened that his sword slipped off and fell into the water. Immediately he cut a mark on the side of the boat and assured himself: "This is where I have dropped my sword."

    By and by the boat came to the destination and stopped. The man plunged into the stream at the point indicated by the incised mark trying to retrieve the lost sword.

    The boat has moved on, but not the sword. To recover his sword this way— the man is indeed muddle-headed .

This prompts us to doubt whether logic is always applicable to other circumstances as we know situation is subjective to constant change. It can be successful in closed systems where every state is well defined, but how about open ones?

Daodejing (Wang Bi, Guo Xiang) begins with: "Dao, daoable, but not the normal dao." Referring to the natural law, we can say it is dao, but it doesn't mean what we say. Whenever we mention it, it is beyond the original sense.

In Daodejing the creator of everything is defined as dao: like a mother that bears things with shape and form. But what/who is dao? It is just unimaginable, because whenever we imagine it, our imagination can never be it (we can never completely describe it: more we describe it, more wrong we are). Daoism illustrates the origin of everything as such a form that doesn't show in any form we can perceive. Whatever we can perceive is merely the created forms, rather than its genuine nature, as if we know people by their outer looks rather than by their inner intentions. We are too far from understanding the nature.

**Daodejing suggests that logic in the most original extremity is shapeless in nature:** it is unbodied, invisible, inexpressible, or even intangible.
- We frequently have such a feeling in learning English as a foreign language that there is no fixed meaning but an intangible impression or feeling to a word: the meaning varies with situations, contexts or even ages so that we can never assure our comprehension. In fact, it is due to the unbridled usage in logic made by people of different ages and districts—there is always creativity implied in the word so that we can rely on nothing more than our own creativity.

- Wherever there is logic, there is also the corresponding comprehensive understanding to logic: also based on our creativity relative to the different existences of human beings. Due to the diversity of comprehension, interpretation and creativity to the same logic, there are varying versions of perception (or conception) to the same logic.

Whenever we say "it is" by logic, we are subjective—how can we assure our objectiveness? We may have developed it from some strictly limited domain under constrained conditions, or it is merely a haphazard, because we can only observe limited cases in our limited lives. How can we convey the logic to those with different backgrounds, even with the slightest difference? How can we then assure the determinacy of its truthness? Is that the reason that majority of people hold it (e.g., Darwin's evolution theory and the functional difference between left and right brain—both are controversial in fact)? In fact none of the logics can be proved, even of we exist or not (Gershenson).

There is no truth and false actually: there is because the outcome has to meet someone's desire—they are merely the attributes of a tradeoff. One false dead can be true in another perspective, e.g., eating much is good, because of the excellent taste and nourishment, but it is also bad when he gets weighted. Neutrosophy (Smarandache) shows that a true proposition to one referential system can be false to another.

Conclusion:

Validity of logic depends on the way we reflect it, not logic itself. Logic never proves itself.

**Logic is a matter of tradeoff (balance) between contradictory factors. There seems no absolute correctness or falseness independent of environment.**

There is dao, but not the kind we mentioned, **accordingly, there is logic but not what we specified**.

## 6. Logic is only One Perspective of Learning, not an Independent Entity

Logic comes from perception and leads to new perception. It is shown that **human understands the world through the interaction of the inter-contradictory and inter-complementary two kinds of knowledge: perceptual knowledge and rational knowledge——they can't be split apart.**

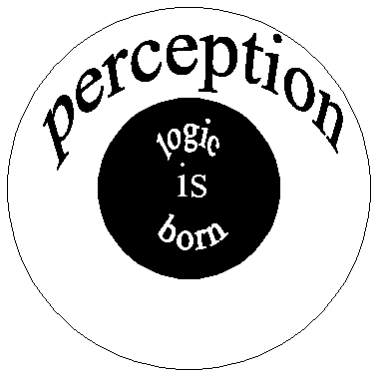
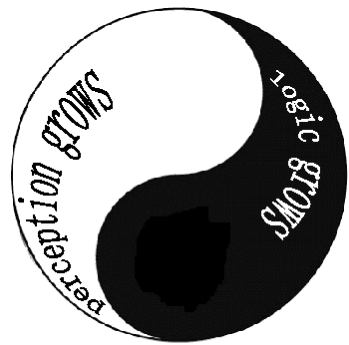

Logic is created through perception in which void, intangible feelings or impressions have been nurtured, brought up and developed into a mental model.

The logic born is nothing more than a subjective hypothesis at primitive stage—it is not within the sense of truthness and falseness.

However, the terms (symbols) and syntax (rules) is only understood by perception through practice. How can we imagine a bookworm who well reads books but has no experience?

Where are truth and false born? There are no such beliefs at the first stage of practice in fact (Daodejing):

> When beauty is abstracted (Peter A. Merel)
> Then ugliness has been implied;
> When good is abstracted
> Then evil has been implied.
>
> So alive and dead are abstracted from nature,
> Difficult and easy abstracted from progress,
> Long and short abstracted from contrast,
> High and low abstracted from depth,
> Song and speech abstracted from melody,
> After and before abstracted from sequence.

> So it is that existence and non-existence give birth the one to (the idea of) the other (James Legge); that difficulty and ease produce the one (the idea of) the other; that length and shortness fashion out the one the figure of the other; that (the ideas of) height and lowness arise from the contrast of the one with the other; that the musical notes and tones become harmonious through the relation of one with another; and that being before and behind give the idea of one following another.

> Accordingly (back to authors), the division between truth and false comes from practice: there is truth, because the outcome is desired, and vice versa to false. Furthermore, without false, where comes the truth? And without truth, where comes the false? Human has been unintentionally comparing, weighing, balancing and trading off between favorable and unfavorable, desired and undesired, based on his prompt subjective and objective situations, hence comes the abstraction (distinction).

> Therefore, **the distinction (or abstraction) of truth and false is nothing more than desires that are subject to constant change.**

Conclusion:

**Truth is born from false and false from truth. They are exactly the measurement of men's practice. This measurement is by no means isolated from practical situations.**

- Discover the truth through practice, and again through practice verify and develop the truth. Start from perceptual knowledge and actively develop it into rational knowledge; then start from rational knowledge and actively guide revolutionary practice to change both the subjective and the objective world.

Practice, knowledge, again practice, and again knowledge. This form repeats itself in endless cycles, and with each cycle the content of practice and knowledge rises to a higher level.

Accordingly, we need to represent learning in integral form:

***Learning=∫d(perception)· d(logic)= ∫d(perception(object))· d(logic(object))***

in time, space and situation domains, or in compound form:

***Learning=∫d(perception· logic)***

since perception and logic are interchangeable and inter-transformable (they melt each other), or they are decomposable in the same manner as learning.

### 7. As a Part of Learning, Logic is Dynamic

Logic depends on perception which is subject to dynamic change with environment, i.e., **the truth value swings**.

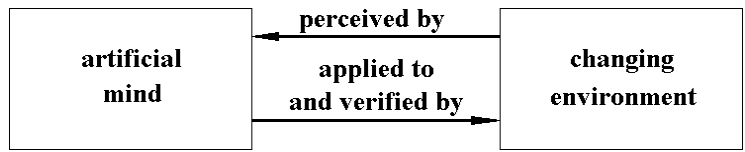

For the logic (Liu [7]):

"I'll visit him if it doesn't rain and he is in."

To avoid being trapped in an instant case that "the clouds promise impending rain" or "there is no answer at the moment I ring the door", we need to learn the long-term trend like "does it rain whole day" or "does he keep his promise", to make a feasible plan (wait or return).

Conclusion:

***Learning=∫d(perception· logic)***

in all its time, space and situation domains. **Whenever we persist in some instant or partial look, we loose the whole. Furthermore, whenever we are satisfied with one-sided view, we also loose it.**

### 8. As a Part of Learning, Logic is Multilateral

Human is normally too confident of himself, no matter how partial or monolateral he is, so that he always misleads himself.

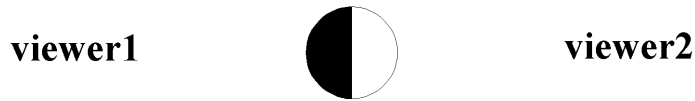

**A ball with a black hemisphere and a white hemisphere**

He may then solve such contradiction by standing on some other perspective. However, **he doesn't succeed until he reaches the opposite side**. This is where our yin-yang philosophy starts.

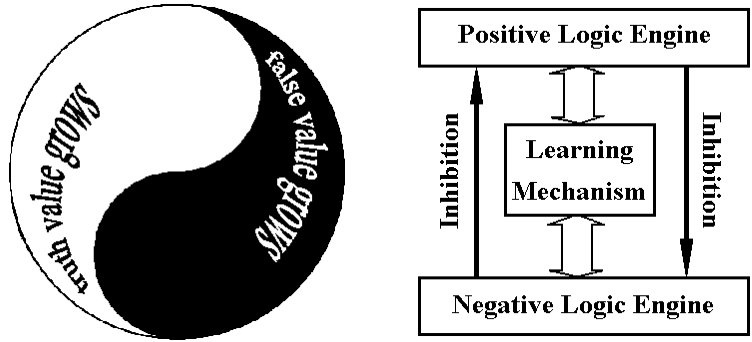

**Two kinds of partial reasoning are alternatively fired**

The incompleteness of human mind (by no means can we assert the perfectness of human being) indicates that **human always reasons in partial mode, i.e., positively partial at one time and negatively partial at another**, as indicated by taiji figure. There are plenty of reasons:
- As a holograph of the universe, human behaves in rhythmed way, e.g.: positive mode, negative mode, positive, negative… and so on, with each mode complementing and inhibiting the opposite one.
- Just because human sometimes stands on positive perspective and then the other, the truth value (it is, not it is) is in constant change, as shown in dynamic state to us.
- There remains a learning procedure hidden in the above bilateral logic: human has to balance the bilateral reasoning to adapt to his present or long-term situations, or meet his needs.
- It is from this inter-complementary and inter-inhibitory contradiction: the bilateral model, that multilateral system is generated, according to Chinese yin-yang philosophy or I-ching (in Chinese: Yijing, also known as the Book of Changes).

We can never base intelligence on the individual behavior, and this is the reason why we need group or society to exchange our views. This is also the underlying essence of multiagent approach that tries to simulate a society.

Based on our intensive exploration in Chinese philosophic perspective, a prototype of **logic cell is presented as an inter-complementary pair**: the positive and the negative logic engines represent positive partial reasoning and negative partial reasoning respectively, and a learning mechanism tries to assemble them to figure out the trend, based on its present and long-term situations, and make choices among various possibilities, e.g., to try a plan among a couple of possible plans. Suppose that each logic engine is decomposable in the same manner, which clones the entire system.

It starts with default logic—although believed absolutely valid, it is actually partial in nature as we know that human is by no means complete as long as he thinks in logic. This partial activity is warmed up and up until some time (one day) contradiction arises. This contradiction, growing up and up in previous partial mode, gives rise to negative reasoning, which later on inhibits the original logic. This inter-complementary process continues as the loop goes on and on, during which the contradiction tends to be neutralized. However, this is only a temporary balance when the two engines reach an agreement. New contradiction comes with the constant change in environment or situations outside. There still remains a chronic cycle hidden in the rhythm for long-term resolution. It is also important to note that this is a genetic proliferation in both time and space complexity, for it can clone all its subsystems.

More intensive study is being carried out in neutrosophy.

Conclusion:
A practical reference frame should originate from a single contradiction or yin-yang (see Daodejing or I-ching).
**Every existence is of bilateral character, or double characters, with each opposing and complementing the other to form a unity.**

$$d(object) = \frac{\partial(object)}{\partial (positive\ perspective)} d(positive\ engine) + \frac{\partial (object)}{\partial (negative\ perspective)} d(negative\ engine)$$

$$= \frac{\partial (object)}{\partial (yang)} d(yang) + \frac{\partial (object)}{\partial (yin)} d(yin) = \frac{d(object)}{d(reference\ frame)} d(reference\ frame)$$

Whenever we hold logic, we have already been standing on a default perspective. Is there universal logic? No, unless we reconceptualize it in an opposite perspective, e.g., Daoist or Buddhist view.

Conclusion:
When we hold logic, we have already believed that it is something. This belief in turn inhibits our negative consciousness that it may be something else (to some degree) or it can be another thing simultaneously (to some extent). We are in this way trapped. So:
"It is never too old (for a machine) to learn."

## 9. Concluding Remarks: Illusion and Creativity

It has long been illustrated in Daodejing (Wang Bi, Guo Xiang) that whenever we capture dao as the natural law, universal method or logic, etc., what we capture can never be it. Therefore:
**Although we can learn logic, we can never capture it.**
**Whenever we do, what we capture is merely a distortion.**
A famous poem from "Topic to Xilin Wall" by Su Shi, a great poet in the Chinese Song Dynasty:

>A great mountain by vertical and horizontal view,
>Far, near, high, low, and each not same.
>I can't see the true face of Lushan,
>Because I am just in there.

We cannot criticize ourselves just because we habitually and absolutely believe in ourselves. Accordingly, we cannot keep a critical mind to our so called truth just because we habitually and absolutely believe in so called truth. Is there really some kind of (absolute) truth on earth?

>One time in Tang dynasty China, the Fifth Patriarch of Buddhism announced to his disciples that everyone write a verse to show his insight of the Buddhist wisdom.
>
>At this, the most eligible one presented on the wall the verse:
>>Our body be a Bodhi tree,
>>Our mind a mirror bright,
>>Clean and polish frequently,
>>Let no dust alight.
>
>Just as a choreman in the mill of the temple, Huineng answered it with his own:
>>There is no Bodhi tree,
>>Nor stand of a mirror bright,
>>Since all is void,
>>Where can the dust alight?

**Whenever we hold the belief "it is …", we are loosing our creativity. Whenever we hold that "it is not …", we are also loosing our creativity.** Our genuine intelligence requires that we completely free our mind - neither stick to any extremity nor to "no sticking to any assumption or belief" (Liu [8]).

As we mentioned previously, whenever there is truth, there is also false that is born from/by truth—this abstraction (distinction) is fatal to our creativity.

Meanwhile, our creativity is nothing similar with things created (i.e., in the sense if truth and false). It must be void in form (no definite form), something like dao, since whenever we hold it, it is not our creativity. Nor does it mean to destroy everything (there is nothing to destroy nor such action, if there was, it is no longer void).

>**There is nothing to destroy, nor anything to create.**
>**If there is, it is rather our illusion than our creativity.**

Because everything believed existing, true or false, is nothing more than our mental creation, there is no need to pursue these illusions, as illustrated in the Heart Sutra [3]:

>When Bodhisattva Avalokiteshvara was practicing the profound Prajna Paramita, he illuminated the Five Skandhas and saw that they are all empty, and he crossed beyond all suffering and difficulty.
>
>Shariputra, form does not differ from emptiness; emptiness does not differ from form. Form itself is emptiness; emptiness itself is form. So too are feeling, cognition, formation, and consciousness.
>
>Shariputra, all Dharmas are empty of characteristics. They are not produced, not destroyed, not defiled, not pure; and they neither increase nor diminish. Therefore, in emptiness there is no form, feeling, cognition, formation, or consciousness; no eyes, ears, nose, tongue, body, or mind; no sights, sounds, smells, tastes, objects of touch, or Dharmas; no field of the eyes up to and including no field of mind consciousness; and no ignorance or ending of ignorance, up to and including no old age and death or ending of old age and death. There is no suffering, no accumulating, no extinction, and no Way, and no understanding and no attaining.

    Because nothing is attained, the Bodhisattva through reliance on Prajna Paramita is unimpeded in his mind. Because there is no impediment, he is not afraid, and he leaves distorted dream-thinking far behind. Ultimately Nirvana! All Buddhas of the three periods of time attain Anuttara-samyak-sambodhi through reliance on Prajna Paramita. Therefore know that Prajna Paramita is a Great Spiritual Mantra, a Great Bright Mantra, a Supreme Mantra, an Unequalled Mantra. It can remove all suffering; it is genuine and not false. That is why the Mantra of Prajna Paramita was spoken. Recite it like this:

    Gaté Gaté Paragaté Parasamgaté
    Bodhi Svaha!

Conclusion:

    **Everyone can extricate himself out of the maze of illusion**, said Sakyamuni and all the Buddhas, Bodhisattvas around the universe, their number is as many as that of the sands in the Ganges (Limitless Life Sutra, Chin Kung).

# THE INTANGIBLE ABSOLUTE TRUTH

Gh. C. Dinulescu-Campina

Str. Bucea, No. 4, Bl. 38, Ap. 31
2150 Campina, Prahova, Romania

In the work "The Modelling of Rationality" on the basis of my own MESER licence, I have raised a new spiritual doctrine sustained by scientific and logical hypotheses. The perception of the soundness of the mentioned concept proceeds both from the Leibnizian principle concerning the sufficient reason, and from Einstein's principle regarding internal perfection and the external acknowledgement of a new theory but, like any responsible "creators", I felt the need to also consider the expression of the feeling of uncertainty, mine first.
Although I had found many external confirmations in our great forerunners' ideas and theories, I have not had a proven substantiation yet (which is not by all means necessary with philosophical hypotheses) of the hypotheses that I have forwarded, until I got acquainted with the ideas of the mathematician and philosopher Florentin Smarandache - the creator of Neutrosophy - as a branch of Philosophy, that studies the origin, the character, the aim and the interactions of the "objects" from the idealistic spectre.
> I've found out that the Neutrosophy Theory, belonging to the mentioned thinker, based on a non-Manichean logic - that is, the trivalent logic- sets up as the scientifically demonstrated fundament for the great majority of the hypotheses I have set forth in "The Modelling of Rationality".

Essentially, Professor Smarandache's Neutrosophy stipulates that for any idea <A> there is also an idea <Anti-A> that does not mean <Non-A>. The fundamental thesis of Neutrosophy is: if <A> is t% true and f% false, as bivalent extremes, it is necessarily i% indeterminate (=achievable, to outline its probabilistic connotation), to the effect that, $t+i+f \leq 300^+$ (or $t\%+i\%+f\% \leq 3^+$) which gives a slightly altered meaning to some common concepts such as, for example, the one of complementarity. To this effect, the complementary of t is not f, but i+f, while the complementary of f is not t, but t+i.
Florentin Smarandache's theory of Neutrosophy suggests also the fact that any hypothesis has a nature of extremeness (it also allows an anti-hypothesis) which is not bad because the law t+i+f=100 must be considered dialectically, where both t and f tend to be decreasing (without annuling each other) to the advantage of i. Far from the idea that any hypothesis should not have a nature of extremeness, just such a nature is desirable to generate polemics which, in case of confrontation, draws nearer t and f to one another, aiming at the neutral equilibrium of the t+f+i=100 relationship, that provides the opportunity of accomplishment. (As regards the opportuneness of polemics, I would like to mention that the author of neutrosophy hasn't yet accepted the "realizable" alternative as "indeterminate", nor the impossibility for t and f to make null one another.)

The theory of Neutrosophy makes obvious the relative nature of the truth and the false, only the neutral nature tending to the absolute owing to its force of accomplishment.

*THANKS TO THE SPECIFICATIONS THAT ARE STIPULATED IN SMARANDACHE'S NEUTROSOPHY, SOME HYPOTHESES OF THE MESER CONCEPT SUCH AS: THE COMPLEMENTARITY BETWEEN THE SACRED AND THE PROFANE, BETWEEN THE DIVINE CREATION AND THE INTRA-SPECIFIC EVOLUTION, THE NON-CONTRADICTION BETWEEN SCIENCE AND RELIGION, BETWEEN MATERIALISM (SUBSTANTIALISM) AND IDEALISM, BETWEEN GNOSTICISM AND AGNOSTICISM, PROVE TO BE RATIONAL AND THEREFORE REAL, WHILE THE PARADOXES BECOME JUSTIFIED.*

Directly related to the intangibility of absolute truth, and tackling the issues of the aim of knowledge, according to the neutrosophical fashion, the MESER concept identifies two modalities: scientific knowledge - that specialised knowledge "more and more from that <<less and less>> and philosophic, encyclopaedic knowledge "less and less from that <<more and more>>". If the first modality of knowledge is limited especially by the possibilities of communication, the second one is also limited by the insufficient power of comprehension of the human mind. The equilibrium between the two directions which, in the last analysis, signifies the way to the truth, is determined by the divine laws of dissociation, purification (the selection and the erasing of the seals) and those of monadic recomposition - laws that ascertain for knowledge as a whole, a social character, expressed by the syntagm "more and more from that <<more and more>>, rendered by the well-known paradox "the more you learn, the less you know."

After all, the fundamental law of Neutrosophy is a successful attempt for resolving the paradox of knowledge and confirm the thesis that the absolute truth is intangible not in a derogatory way but in an optimistic one, approved and revealed by (and through) the will of God.

Being operative even in the case of particular interpretations, as is the case of the present one, Smarandache's neutrosophy confirms (according to Einstein theory) its validity, be it only for the fact that it suggests new methods and modalities for evaluation, new interpretative perspectives.

*On time*

In the paradigmatic construction of the spatiotemporality, the MESER concept has recourse to the necessary hypothesis according to which substantial reality – the one disseminated through the divine will within into material reality and anti-material reality possesses that "retrograde movement" from right to left, in respect to the spiritual reality. Independently of the behaviour of any rational (or not) materialized entity, the above mentioned retrograde movement meaning for each of those the outlining of its "embodied" existence, represents what is defined as that impalpable philosophical category called TIME.

One can assume that this movement related to the spiritual reality as a benchmark is characterized by uniformity (constant "speed") – hence the perception of superior rational entities (humanoids) on that uniform "flow" of time and on its reversibility.

Physical movement – that movement of any entity in respect to substantial reality as a benchmark, with parameters (speed, for example) that cannot be compared to those of the hypothetical retrograde movement of substantial reality – cannot change that apodictic character of the uniform time flow impression. From the time when physical movement, though, acquires parameters comparable to those of the retrograde movement (such as one of the so called cosmic speeds), one can speak of a different modality (slowed down) of time flow.

Speculations on time expansion (substantiated by the famous Einsteinian theory of relativity), already apodictic by virtue of their scientific support, concede that each entity bears its own time. A spaceship wandering through the immensity of substantial reality, using that already feasible cosmic speed, is a bearer of its time – a time characterized by an incredible and yet paradoxical dilatation in respect to the time of the departure station. Upon return of such a spaceship from a voyage of only a couple of years, one could notice that at the departure station several generations have already "passed". The voyagers of such a spaceship do not notice the dilatation of their own time (the spaceship"s time) but, on the contrary, the contraction of the time of the station.

Each rational entity fails to perceive the change of its own time, it only perceives the dilatation and contraction respectively, of the time of another entity in respect to one's own time.

Admitting that time is the overall movement of substantial (material) reality in respect the eminently spiritual (nonsubstantial) reality, but also the movement of any material entity in respect to the substantial reality (physical movement) we infer that the former, i.e. time, does not exist within the non-substantial reality. What is more, one can state with minimal risks, that time and space are "mere voyagers" within the spiritual reality.

The absence of time (and space) from the spiritual reality does not involve the absence of movement from that reality. On the contrary, the spiritual reality supposes, of necessity, the movement of nonsubstantial entities, that a-causal movement (from the physical point of view), the permanent stochastic movement of entities (entelechies), eminently spiritual, governed by haphasard (as a perfect law) – the only workable law, both within the spiritual reality and within the substantial one (material or anti-material).

As regards its triple polarity (spatiality, duration and order) one can point out that the spatiotemporality is governed by order (resulted from haphazard), time and space being the "gifts" with which God endowed the spatiotemporality when He created the matter and anti-matter, by their dissemination from the spiritual reality.

In order to tackle the tricky and delicate problem of what is called reversibility of time, some considerations stand out:

– I have defined time as that retrograde movement of substantial reality, having as a benchmark the spiritual reality where the material, and antimaterial, respectively, reality "floats" osmotically. According to the MESER axiom for the spatiotemporality, that retrograde movement was taken into consideration out of the need to justify the "voyage" of any entity, specifically of the "live" entity, between the primary limit and the secondary limit of material reality, that is from birth to that "extreme and insurmountable possibility" as Otto Poggeler defines corporeal death, or, in another manner, to justify "consumption" as a duration of that factual experience called DASEIN by M. Heidegger.

− As an element of the spatiotemporality, any rational (or not) entity has, in respect to spiritual reality, a "relative stability" (besides the inner Brownian movement of those in the same Brownnian style due to one's own will − the free will), having no reasons or capabilities to cross the manifest (substantial) reality in between its two limits.

− The "retrograde movement" and the "relative stability" confirm the supposition that any entity is a bearer of its own time within the material reality and that the former, time, does not exist within the spiritual reality.

− If material reality as a whole has a retrograde motion, any entity (and specifically the live one) has a contrary movement. To this effect, one can deem that the entities, whichever they are, "travel" uniformly towards the future within the substantial reality and undergo a perpetual present in respect to the spiritual reality, as long as they are embodied (materially split).

The foregoing considerations lead to the conclusion that the hypothesis of the retrograde movement, as a time defining hypothesis, as well as that of the relative stability do not grant to an entity the possibility to travel in a retrograde fashion (into the past), thus any possibility to explain an assumed temporal reversibility is ruled out.

By invoking the Einsteinian theory of relativity according to which temporal dilatation by travelling with cosmic speeds, or temporal contraction, as perceived by the rational entities left behind at the "station" (from which the spacecraft departed and to which it returned), are possible, there is no way to speak of a hypothetical reversibility of time. Indeed, for a rational entity at the station, a travelling entity is by no means "made younger" and in its turn it sees itself in no other way than "made older". The same travelling entity has almost null chances to find any living, and by no means younger, rational entity from among those left at the station.

In order to detect in the temporal "dilatation-contraction" a reason "pro" time reversibility, one should image an experiment during which the "arrival" from, should precede the "departure" for, the cosmic journey. This could be "possible" if the departure and arrival stations were appropriately placed next to the "line of date change", if the entire journey were to last less than 24 hours (time-keeping at the station) and if the so-called "conventional time" were involved. In real time, this is nevertheless not possible, no matter what other additional condition were involved.

Finally, we come to the stage of assessing what means the hypothetical spiritual journey in time (in this case, into the past) − that one made possible by the convenient travel, through that plausible time tunnel, of the soul (the eminently spiritual part of the rational entity), with the possibility to access a certain moment of the past by that entity. This access to the past does not mean a reverse "restoration" of time, but, in the best of cases, only a re-covering of a duration from the accessed moment towards a possible different future (that can be "past" in respect to the commencement of the journey through the tunnel). Obviously, neither the spiritual journey in time provides any reasons to consider that time could be reversible. All (known) modalities of temporal modification that have been invoked − retrograde movement, travel with cosmic speeds, spiritual journey in time − are nevertheless conductive to the idea that: "Fugit irreparabile tempus".

According to the neutrosophic theory (F. Smarandache: *Neutrosophy*) based on trivalent logic, any scientific hypothesis features an extreme character to the effect that it necessarily has also an anti-hypothesis; moreover, to achieve complementarity, the two extremes must also incorporate that indecisive part. By virtue of the above mentioned principle, that is eminently dialectic, the thesis on time irreversibility, as any other one, must be deemed merely as a possible truth, becoming such (or the opposite) only after the acquisition, by reason of some laws, of that "neutral equilibrium" that supposes a possible reformulation and necessarily, an extension of the respective theory in as much as possible, closely connected domains of definition. From among all involved theories for substantiating the thesis on time irreversibility, the one that lends itself mostly to reformulation, and that might lead to a contrary conclusion, seems to be exactly the one with the soundest scientific rationale – the Einsteinian theory of relativity, that Einstein himself did not consider as final.

If in 1905 the theory of relativity had an unprecedented impact on the scientific world, even before the brilliant Niels Bohr could formulate his objections (as Einstein used to behave in respect to Bohr), Einstein had an insight that his theory was based on much too restrictive "preliminary conditions". In 1916 he was "compelled" to reformulate his hypotheses within the framework of the "generalized relativity theory", renaming the initial one as "limited relativity theory".

As the relativistic theory of Einstein is based on the hypothesis on the maximum absolute value of the "C" constant, and as of late (November 2000) three scientists (a Romanian and two Americans) have demonstrated that some phenomena take place at speeds exceeding that of light, it is predictable that the theory of relativity might be reformulated and completed with that part that could be called the "theory of absolute relativity" dedicated to the phenomena occurring in the "mega-cosmos", "giga-cosmos" or "angstrom-cosmos" – hypothetical designations of some definition domains different from the macrocosmos and microcosmos, already "covered" by the Einsteinian theory. The extension and possible reformulation of the theory of relativity could supply substantiating elements on a thesis on "temporal reversibility", but only if "external confirmations" (according to the Einsteinian principle) and theoretical and practical requirements (according to the Leibnizian principle of sufficient reason) for such an approach are provided. In fact, anything is possible, the only absolute truth (figuratively) being that expressed by the syntagm: "nothing is final".

On death

"Dasein is factual Being , power of Being, as anticipation of the extreme possibility that can't be surpassed – death."

Otto Poggeler

According to Heidegger in "Sein und Zeit", the possibility of Being has its origin into an ultimate possibility that cancels, *post factum* even the possibility of Being itself. With his term of *dasein*, Heidegger always brings to the fore existence as a fatality, without explicitly denying its fate as an obligation, either imperative or assumed, nevertheless in vain.

Implicitly, Heidegger is preoccupied, concerning the existence, by the sentimental side of it, without insisting on its role in showing rationality as a great divine work, thus deliberately avoiding to refer to a certain objective aspect.

The MESER concept estimates that sensations as well as feelings, as attributes of existence, are not a purpose but a means, which ensures creating products of rationality through Being. The products of rationality – which are the purpose of Being, are not thwarted, and neither is Being itself, by death; they remain as monadic seals and only the divinity decides through its "Last Judgement" if the seal is to be "erased" or redistributed towards other new spiritual entities, formed as a primary fund of the intellect of future substantially splittable rational entities, which are in turn destined to create new rational products.

While the Heideggerian *dasein* sees death as a fatality which thwarts the possibility of Being, by obsessively bringing it out, the MESER concept (Speculative Existential Model of the Rational Entity) considers death as an objective necessity, whose imminence must not be pointed out in a paranoiac manner and which doesn't cancel the Being but submits it periodically to God's judgement in order to be purified and harmonized with the purpose of His Great Work.

As a subsequent theory to dasein, the MESER concept perceives it (with all due consideration) as a model of sentimental approach (subjective and therefore psychologistic) of the Being as a purpose ending with death; while the MESER concept evaluates itself as an objective, rationalist model, which without denying the purpose itself of Being, qualifies it as a path towards the construction of the divine work which is universal rationality.

    Although they seem to diverge, the two models are contingent and neither one won't be affected if the balance of truth inclines towards one or another.

Even though the two models would be, judging by bivalent logic, the extremes, each one prevails based on an indisputable accurate logic.

According to the basic law of the Smarandachian neutrosophism (F. Smarandache, University of New Mexico:Neutrosophy) $t\% + i\% + f\% \leq 3^+$, however t (the undoubted truth) and f (the undoubted falseness) both have the tendency to decrease (without annulling) one for the advantage of the other or for the advantage of i (the irresoluteness, the achievable), having effects on creating a possible model, as much as operable, of Being. Here $t = \sup T$ (truth), $i = \sup I$ (irresoluteness), $f = \sup F$ (falseness).

    Dasein and MESER are compatible, each of them taking advantage on the other one.

MESER supplies the basic scheme for a bygone Being, from the rational entity point of view, as well as for a periodic and permanent Being, referring to rationality as Being in general; while dasein refers mostly to the sentimental aspect of Being in its substantial sequence. In other words, from the escatological point of view, dasein includes the theory of death as final point of a subjective Being, (which is true), while the MESER concept refers to death as an intermezzo, as another beginning of Being in general.

# Name, Denominable and Undenominable
## - On Neither <A> Nor <Anti-A>

Feng Liu
Dept. of Management Science and Engineering,
Shaanxi Economics and Trade Institute (South Campus),
South Cuihua Road, Xi'an, Shaanxi, 710061, P. R. China
E-mail: youchul@fmmu.edu.cn

**Abstract:** Neutrosophy's underlying construction is far more sophisticated than we can imagine. I present in this paper a critical analysis on the logical description of <Neut-A> based on multicultural joint venture, and reach a contradictory argument that the axiom is rather a paradox than a valid definition. Starting with the fundamental issue in Daodejing, the paper carries out a widespread discussion on conflicts in denominating things, from yinyang philosophy to dialectics, from relativity of being to self-negating effect of concepts, exhibiting the genuine essence of philosophy against distortion. Discussion of feasible description of <Neut-A>  (neutrosophy) is also presented, followed by a brief distinction between human intelligence and machine intelligence. The paper aims to help scientists reach the genuine nature hidden in the ideology of neutrosophy.

**Keywords:** Dao, Genuine Nature, Contradiction, Identity, Self-negation, Partial Negation, Neutrosophy

## 1. Neutrosophy: a joint venture

It is not until recently that I came across the study of neutrosophy introduced by Florentin Smarandache (1995). It seems a very brave challenge to a number of developed sciences and technologies. However, from its intension and method of approach, I realized that it touches the most arcane, abstruse, and mysterious philosophies such as Daoism and Buddhism, and the toughest problems in the universe as difficult as uncovering the universe. I am afraid how western intelligents can handle such mysteries.

As a Chinese, I feel obliged to spread out our exploration based on the multicultural joint venture, with focus on the paradox "neither <A> nor <Anti-A>", as illustrated in step ward arguments, as shown below.

## 2. Name, denominable, but not the normal name

Daodejing (Wang Bi, Guo Xiang) begins with: "Dao, daoable, but not the normal dao; name, namable, but not the normal name." We can say it is dao, but it doesn't mean what we say. Whenever we mention it, it is beyond the original sense.

Daodejing mainly deals with the common problem: "What/who creates everything in the world we see and feel?" It is dao: like a mother that bears things with shape and form. But what/who is dao? It is just unimaginable, because whenever we imagine it, our imagination can never be it (we can never completely describe it: more we describe it, more wrong we are). It is also unnamable, because whenever we name it, our concept based on the name can never be it.

Daoism illustrates the origin of everything as such a form that doesn't show in any form we can perceive. This is the reason why it says, everything comes from nothingness, or this nothingness creates everything in forms in dynamic change. Whatever we can perceive is merely the created forms, rather than its genuine nature, as if we distinguish people by their outer clothes. We are too far from understanding the nature, even for the most prominent figures like Einstein.

You may then ask whether Laozi was genius enough to express it. Definitely no. It is true that he was aware of the problem, or we can assume that he really understood it, however, he could never describe it. Although Carlos Gershenson has presented this argument in his paper as incomplete language, I am still afraid whether he can catch my notion that **the most complete and perfect language is no language**.

- An example is described in a story (Lanier Young, 1991) in the Mahapra Janaparamita Sutra:
    - The Blind Men Trying to Size Up the Elephant: Once there was a king who ordered his minister to bring in an elephant and let some blind men touch the animal one by one. After every one of them had their turn, the king asked them what they thought the elephant was like. The one who had touched its tusk said it was like a turnip; the next had touched its ear, and said it was a dustpan; the third its foot, and said it was a pestle; the fourth its back, and said it was a bed; the fifth its belly, and said it was a jar; and the last its tail, and said it was a rope….
  
  We can imagine that our perceptions are just as partial as those of the blind people, then how can we name things that are believed known to everyone but actually as mysterious as the elephant to the blind people? Do we understand, for example, 1+1=2? We always believe so although it has never been proved by the mathematical world. Even when it were proved, how could we explain black+black or black+white?
- Carlos Gershenson [1] presents: Not only silliness, but all adjectives can only be <used|applied> in a relative way, dependant of a context. Language is relative as well. How can we speak about absolute being, then? We can and we cannot. We speak about it, but in that moment its absolute is relative. For us, it is and it is not-incomplete. But that we cannot completely speak about it, it is not a reason to stop speaking about it (as Wittgenstein would early suggest in his Tractatus Logicus Philosophicus), because we can incompletely represent its completeness... As Wittgenstein himself (but not most of his followers...) realized, following the ideas in the Tractatus, we would not be able to speak about anything... (languages are incomplete). Language is used inside a context. Depending of this context the language will be different.
- Can language be completely transferred by telesthesia? Although I acknowledge it is a better way of communication, but it also depends. For example, how can we understand the hidden, underlying or implied meaning of the transferred "words"? Is language transferable? I am afraid there is not a definite answer. One reason might be that the same language can be interpreted diversely by different people. This should be the reason why human failed to communicate properly with those in different languages, or with other species on earth or in the space, even with Jesus.
- Telesthesia truly exists universally in the Pure Land described in Buddhism: People there communicate with each other by heart rather than words. Further more, they can also "speak" to people everywhere outside the "Pure Land" World which was founded by Amitabha. A recent VCD video from Taiwan shows that a young man really came to this world after his medical death, because he really made frequent communication, after "death", with his sister still in Taiwan as a young Buddhist nun, not only did he answer all the questions of his sister, but also made an unbelievable promise which was later on testified astonishingly by his family on earth.
- Just because there are no perfect words to express the most sincere truth that uncovers our genuine nature, we can say nothing than Amitabha to communicate. Although we know that in Buddhism it is almost equal to saying nothing, but it implies saying everything.

As the conclusion, **everything can be named, but never absolutely proper. It is a name, but never a perfect name.**

**3. Name is always subjective, relative to the perception and perspective of observer**

Has anyone been confused about "Am myself really myself"? An old person, when gazing at the albums of his childhood, he always point to the photos and says: "It is me." How ridiculous!
- At first, he is pointing to some paper cards rather than humans.
- Secondly, what he is pointing to is an image of himself, not really himself.
- Thirdly, provided that he were really pointing to a younger himself, but it were really different from what he looks like presently——they are different themselves objectively, or there are an infinite number of themselves objectively.
- Even at the present age, a smiling himself and an anxious himself definitely look differently objectively.
- They are all himself because humans subjectively take it for granted and firmly believe so.
- But, are all these pertaining to his body really himself? Definitely no, because he will begin another life (in the next life cycle) after his medical death. He never dies actually.
    - Incarnation, samsara, wheel of life, transmigration of souls, or eternal cycle of birth and death, this is a basic phenomenon of every living being in the universe including those in the heaven, except in the Buddhist Pure Lands where everyone has escaped from his destiny, according to Buddhism.
- As the result, we don't actually understand who we are, just use the names subjectively.

Another example lies in question: "What on earth is the following figure?"

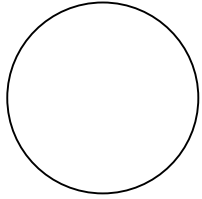

- To well educated students, it is a circle.
- But to uneducated kids, it can be: "a cake, a dish, a bowl, a balloon…, even the moon, the sun".

Let's reach our conclusion that **name is merely our mental creation. It is rather a belief than an objective being, and varies among different people**.
- We always believe "it is" but can't prove it.
    - In his paper "To be or not to be, A multidimensional logic approach" Carlos Gershenson [2] has generalized the proofs:
        - Everything is and isn't at a certain degree. (i.e., there is no absolute truth or false);
        - Nothing can be proved (that it exists or doesn't) (i.e., no one can prove whether his consciousness is right);
        - I believe, therefore I am (i.e., I take it true, because I believe so).
- In fact, this belief of "it is" is always critical (Buddhism).
    - In Buddhist saying, all such beliefs are created by ourselves, for:
    "I am human."
    I am because I always hold this belief, so persistently that I nearly forget I can be Buddha as well.

Multidimensional logic has been surpassed by infinite-valued logic, then fuzzy logic, and ultimately by neutrosophic logic.

**4. Name itself implies anti-name**

Whenever there is a name, it can never be a perfect name. Does it mean we are cheated by or trapped in those created by ourselves? It does, and it doesn't as well. First, there is only relative name, no absolute name. Second, **name actually acts as a tradeoff to unify the diversity of concepts**——whenever there is name, there is contradiction as well.
- We can name something as black, but to distinguish it from white. We can name a human, but to distinguish it from others.
- Names are useful to distinguish things, but don't absolutely describe natures, as stated above.
- Just because there is no absolute name or universal name, people would denominate things in their different perspectives or perceptions.
    - As an example, a modulator/demodulator is named a mouse, just from its casing that resembles a toy mouse.
- Despite all these problems above, we need a common name, however, to communicate. Therefore, we have to make balance, in the sense of acceptance or rejection, among the diversity.
    - Carlos Gershenson [1] points out that: There will always absolute-be injustice, because this one is relative. Since different people have

different contexts (or we can use the word Seelenzustand (soul state), to refer to the personal context, to distinguish from a general context)... So, since people have different Seelenzustandes, we cannot speak of absolute justice, so things will be just for the people with power... The less-catastrophic panorama (and most naive...) would be that the people in the power would have the less-incomplete Seelenzustandes, trying to contain and understand as many Seelenzustandes as they can, so, if they are just, in spite their relativity, they will be just as well for all the people whose Seelenzustandes they contain.

- Accordingly, contradiction is a universal phenomenon that can never be avoided.

The universality or absoluteness of contradiction has a twofold meaning. One is that contradiction exists in the process of development of all things, and the other is that in the process of development of each thing a movement of opposites exists from beginning to end.

> Contradiction is the basis of the simple forms of motion (for instance, mechanical motion) and still more so of the complex forms of motion.

- Despite the anti-name property in space domain we have just discussed, this contradiction also exists in time domain.
  - Mao Zedong also states that:
    We Chinese often say, "Things that oppose each other also complement each other." That is, things opposed to each other have identity. This saying is dialectical and contrary to metaphysics. "Oppose each other" refers to the mutual exclusion or the struggle of two contradictory aspects. "Complement each other" means that in given conditions the two contradictory aspects unite and achieve identity. Yet struggle is inherent in identity and without struggle there can be no identity.
    In identity there is struggle, in particularity there is universality, and in individuality there is generality. To quote Lenin, ". . . there is an absolute in the relative."
  - He implies such a cycle: conflict——identity——new conflict——new identity… to the infinite, in each cycle of which the conception undergoes a partial negation of its original stage to a higher level, and in this infiniteness of negations we make our revolutionary progresses in knowledge, as to negate the original sense of "it is". Hence comes the expression "nothingness" to replace the original meaning "it is something".
    - The Buddhist terms: emptiness, void of the world of senses might just come out of the endless negation of our partial consciousness that has been believed absolutely valid. Once we become aware of it, we are awake. And once we keep the genuine consciousness (without even the

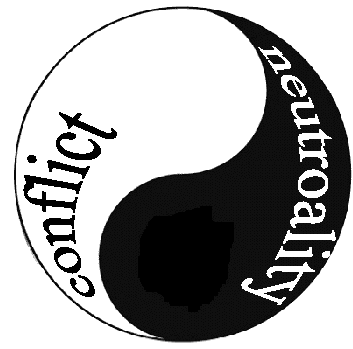

- slightest partialness) in every fraction of moment and forever, we are Buddhas.
- As the result, neutrality comes as the outcome of conflict, and in turn, conflict comes as the outcome of neutrality too. As shown in the taiji form.
  - The more we hold on to our original belief "it is" (although partial), the more mightily conflict arises, since we persist in a more incomplete concept, or fragment, to represent the complete. The same to the coming cycles. This is reflected in neutrosophy as the law of inverse effect (F. Smarandache, 1998).

As conclusion, **antagonism and neutrality are cause and effect to each other.**

### 5. Representing the <Neut-A>

<Neut-A> comes as the consequence of the contradiction between <A> and <Anti-A>, therefore we can say it is neither <A> nor <Anti-A>, but is it all this simple?
- Once we finish **<Neut-A>** as **<neither <A> nor <Anti-A>>** in the first cycle, then in the following cycle we come to the less incomplete concept as **<Neut-<Neut-A>>**. However, <u>**there must be an element of consistency between <Neut-A> and <Neut-<Neut-A>>,**</u> i.e.,
   between  **<Neut-A>**  *and*  **<neither <Neut-A> nor <Anti-<Neut-A>>>**
  <u>**that acts as the gene of reproduction**</u>, i.e., there must be a consistency
       between  **<B>**  *and*  **not <B> and not <C>**
  But where is it? Apparently they are not logically consistent! Once there is nothing indeterminate that can pass down to the next "generation", is it a feasible philosophy?
- <Neut-A> would alternately be expressed as both <A> and <Anti-A>, since it shares characteristic of both <A> and <Anti-A> to some extent (I prefer "extent" than "degree").
  - Although yin and yang are opposite in taiji figure, they are unbreakable friends in giving birth to novel form of development.
- The point lies in the confusion between absolute-be and partial-be (or absolute-not and partial-not). We can confirm or negate a being absolutely and partially as well, but I am afraid we can never be genius enough to do it absolutely unless in a well defined domain in which our perceptions are relatively complete.
- How can we properly express this partial approval/negation? In percentages? Then how can we deal with such partial operations in fuzzy and neutrosophic sets? It is quite another task.
- Provided that a genius giant had successfully settled it in pure mathematics, there would be no need then, I am afraid, to employ analog means. How incredible!
- A feasible alternative, I suggest, would be to put dynamic weight on concept instead of statistic percentages, to combine neural technology. Or more specifically, to create a pattern neutrosophically from the threshold ideology in neural network approaches.

As the conclusion, it is never too old for a machine to learn.

6. **On our genuine intelligence——creativity**

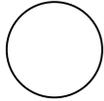 In the previous query about the figure on the left, whenever we hold the answer as a circle, we are inhibiting our creativity. Nor should we hold that it is a cake, a dish, a bowl, a balloon, or the moon, the sun, for we also spoil our creativity in this way. Then, what is it?

"It is nothing."

Is it correct? It is, if we do not hold on to the assumption "it is something". It is also wrong, if we persist in the doctrine "the figure is something we call nothing." This nothing has in this way become something that inhibits our creativity. How ridiculous!

**Whenever we hold the belief "it is …", we are loosing our creativity. Whenever we hold that "it is not …", we are also loosing our creativity.** Our true intelligence requires that we completely free our mind——neither stick to any extremity nor to "no sticking to any assumption or belief". This is a kind of genius or gift rather than logic rules, acquired largely after birth, e.g., through Buddhism practice.

**Not (it is) and not (it is not),
It seems nothing, but creates everything,
Including our true consciousness,
The power of genius to understand all.**

# Growing and Anti-growing Metabolism in Neutrosophy

Feng Liu
Dept. of Management Science and Engineering,
Shaanxi Economics and Trade Institute (South Campus), South Cuihua Road
Xi'an, Shaanxi, 710061, P. R. China
E-mail: youchul@fmmu.edu.cn

**Abstract:** The paper presents the introductory discussion on the concept of neutralism in neutrosophy in a new perspective: metabolism, as represented by the two opposing forces—the growing force and the anti-growing force. The discussion is based on the following issues: dialectics, I-ching and a hypothesis on unified field theory. Metabolism is another crucial issue in neutrosophy, that may probably refresh the definition, logical structure and the ontology of <Neut-A> presented by Florentin Smarandache. This is shown in the discussion on the representation of metabolism, where the paper suggests a conjunction point between neutrosophy and machine learning including neural learning.

**Keywords:** Neutrosophy, Metabolism, Contradiction, Identity, Whirl Vortex Body, Unified Field, Yin-yang

## 2. Background

Neutrosophy is defined as a new branch of philosophy that studies the origin, nature, and scope of neutralities, as well as their interactions with different ideational spectra (Florentin Smarandache, 1995). However the term neutrality needs to be discussed, explained or redefined in multicultural background, since neutrosophy covers a large area of studies and has to be based on the occidental and oriental joint venture (Liu [4]).

- This mode of thinking includes:
    1. proposes new philosophical theses, principles, laws, methods, formulas, movements;
    2. reveals that world is full of indeterminacy;
    3. interprets the uninterpretable;
    4. regards, from many different angles, old concepts, systems: showing that an idea, which is true in a given referential system, may be false in another one, and vice versa;
    5. attempts to make peace in the war of ideas, and to make war in the peaceful ideas;
    6. measures the stability of unstable systems, and instability of stable systems.
- Methods of neutrosophic study involves:
    - mathematization (neutrosophic logic, neutrosophic probability and statistics, duality),
    - generalization, complementarity, contradiction, paradox, tautology, analogy,
    - reinterpretation.

- Neutrosophy touches the most arcane, abstruse, and mysterious philosophies such as Daoism and Buddhism, and the toughest problems in the universe as difficult as uncovering the universe.

Many issues on the explanation and description concern space domain, where I have put my first analysis in a paper based on Daoism and Buddhism (Liu [4]). Here I am going to show some points pertaining to time domain.

3. **Metabolism in the nature**

As known to everybody, metabolism is a basic phenomenon of nature: something is dying as the premonitor of another growth, and something growing as the outcome of dying of another.

- When plant is dying, it yields seeds to another growth that adapts more to the geological situations like climate, soil and water.
- Animal is alive when one part of the cells dying and another part growing, which inherit its will.
- Human beings generate themselves in the same way to inherit families, knowledge and faith.
- As in society, e.g., in Chinese dynasties, new forces grow as the dying of old ones.
- It is also said in Buddhism that the earth undergoes such cycles in its evolution: the birth, the growing, prosperity, decrepitude and vanish, just as the existence of plants in the sequence of spring, summer, autumn and winter.

More issues can be found in a forthcoming encyclopedia of Buddhism (the only one based on modern society) which I am going to put into English.

- The alternation of yin and yang is universal in everything:
  - Yin-yang is the basic form of everything and every course.
    It is yin and yang that is meant by Dao (I-ching)).
  - There is Qian, symbolizing Heaven, which directs the great beginnings of things, and Kun, symbolizing Earth, which gives them to their completion (I-ching, James Legge).
    Qian pertains to yang and Kun to yin. In the course of development and evolution of everything Qian (yang) acts as the creativity (creativity is mentioned in my conference paper, (Liu [4])) that brings new beginnings to it, whereas Kun (yin) implements it in forms as we perceive as temporary states. It is in this infinite parallelism things inherit modifications and adapt to changes (author's note).
  - There is always contradiction between them, which needs to be neutralized in a new form.
    Yin acts as the preservative force (inhibiting yang: when something is shaped, it should have its own inertia against immediate change) that resists modification, reformation or revolution exhibited by yang, therefore, there should be conflict or struggle between these two forces (author's note):

- There are two states of motion in all things, that of relative rest and that of conspicuous change. Both are caused by the struggle between the two contradictory elements contained in a thing. When the thing is in the first state of motion, it is undergoing only quantitative and not qualitative change and consequently presents the outward appearance of being at rest. When the thing is in the second state of motion, the quantitative change of the first state has already reached a culminating point and gives rise to the dissolution of the thing as an entity and thereupon a qualitative change ensues, hence the appearance of a conspicuous change. Such unity, solidarity, combination, harmony, balance, stalemate, deadlock, rest, constancy, equilibrium, solidity, attraction, etc., as we see in daily life, are all the appearances of things in the state of quantitative change. On the other hand, the dissolution of unity, that is, the destruction of this solidarity, combination, harmony, balance, stalemate, deadlock, rest, constancy, equilibrium, solidity and attraction, and the change of each into its opposite are all the appearances of things in the state of qualitative change, the transformation of one process into another. Things are constantly transforming themselves from the first into the second state of motion; the struggle of opposites goes on in both states but the contradiction is resolved through the second state. That is why we say that the unity of opposites is conditional, temporary and relative, while the struggle of mutually exclusive opposites is absolute. (Z. Mao)

However, there is neutralism in them:

- The fact is that no contradictory aspect can exist in isolation. Without its opposite aspect, each loses the condition for its existence. Just think, can any one contradictory aspect of a thing or of a concept in the human mind exist independently? Without life, there would be no death; without death, there would be no life. Without "above", there would be no "below"; without "below", there would be no "above". Without misfortune, there would be no good fortune; without good fortune, there would be no misfortune. Without facility, there would be no difficulty; without difficulty, there would be no facility. Without landlords, there would be no tenant-peasants; without tenant-peasants, there would be no landlords. Without the bourgeoisie, there would be no proletariat; without the proletariat, there would be no bourgeoisie. Without imperialist oppression of nations, there would be no colonies or semi-colonies; without colonies or semi-colonies, there would be no imperialist oppression of nations. It is so with all opposites; in given conditions, on the one hand they are opposed to each other, and on the other they are interconnected, interpenetrating, interpermeating and interdependent, and this character is described as identity. In given conditions, all contradictory aspects possess the character of non-identity and hence are described as being in contradiction. But they also possess the character of identity and hence are interconnected. This is what Lenin means when he says that dialectics studies "how *opposites* can be . . . *identical*". How then can they be identical? Because each is the condition for the other's existence. This is the first meaning of identity (Z. Mao).

> As an example, although landlords and tenant-peasants are opposites, they have to rely on each other to coexist (author's note).

- But is it enough to say merely that each of the contradictory aspects is the condition for the other's existence, that there is identity between them and that consequently they can coexist in a single entity? No, it is not. The matter does not end with their dependence on each other for their existence; what is more important is their transformation into each other. That is to say, in given conditions, each of the contradictory aspects within a thing transforms itself into its opposite, changes its position to that of its opposite. This is the second meaning of the identity of contradiction (Z. Mao).

  > As an example, landlords can fall into tenant-peasants once they loose their fortune and power, and vice versa when tenant-peasants gain fortune and power (author's note). Or quite Mao Zedong's, "by means of revolution the proletariat, at one time the ruled, is transformed into the ruler".

- Neutralism is derived when yin and yang combine each other into compounds as illustrated in I-ching in the creation of trigrams, or when yin and yang inter-consist each other in an infinite course: yin, yang, yin in yang, yang in yin, yang in (yin in yang), yin in (yang in yin), …, at last, they reach a balance as the new form or new generation.

4.  **The metabolic universe**

As stated in his book Gravitational Measurement in Time and Space — A Hypothesis of Life and Sole, Jiang Xiufu has structured a kind of entirely new and total scientific theory. Starting from the different theories in principle, he has revised the physical theories developed since the 17$^{th}$ century completely and bodily. Also, from Aristotle's mechanics to Einstein's special theory of relativity and from classical electrodynamics to modern quantum mechanics, the book has generalized the thoughts of unified field theory in a broadly comprehensive way and made a thoroughgoing change in the concepts of today's natural sciences. Part of the summary is shown bellow:

> **Universe is made up of matter**: Whirl vortex body is the basic form of matter existence. As viewed from its attribution of nature, all the motions are the position changes of matter gains in space according to order (time). Whirl vortex bodies are in the frictions among matter gains as well as in non-central collision, thus causing them to move along the inclined direction and also making themselves begin rotate "left" or "right". Thus the turbid, disordered, pointless and random rotating motions begin to contract and rotate violently in the dominant "left" or "right" whirl vortexes and congregate to form universe various steps and whirl vortex body of different sizes… Therefore, the nearer to the center, the larger matter density and the faster rotation velocity, the whirl vortex bodies have. And, the whorl vortex bodies with small sizes have larger matter density and larger matter gains around with the rotating linear velocity. The maximum internal energy is congregated in electrons and atoms. In the super-density environment of electrons and atoms, the time-space scale structures contracts to their limits. As a result, all the energy exchanging process takes places and ends at an instance.

Whirl vortex steps and stairs in the matter world are in finite in time and space so that we regard the universe as self generating and developing system of the "left" and 'right" vortex body in each step or stair. "Left whirling matter" and "right whirling matter"----matter and anti-matter are the opposite of two sides existed within the matter world. And the interdependence of the two sides upon each other as well as both opposite and compensate each other. They interdependent upon each other, inter-supply each other with that they need and interconsume what they have. Accordingly, they themselves promote their own to grow and die as well as change and develop in the eternal cycles.

**Space is the form of matter existence.**

**Time is the attribution of matter motion.**

**The world unifications are their materiality.**

**The law of the unity of opposite is the basic law of the matter world:** The objects in physics are in continuous changes. On the one hand, they inter-act with their surroundings and they have the history of growth and death; on the other hand, the occurrence of the natural phenomenum of any kind is constrained by the internal inherent contradictions of any matter. Therefore, physics should repel all the rigid conceptions. There are no absolute rest and uniform straight line motion; there are no motions without resistance, there are no absolute "positive and negative as well as left and right"; there are no absolute time and space…; For this reason, starting from the fixed facts, it is necessary to study physical logic and physical concepts and their connotations and denotations of how to express the history of transition, conversion and interaction in time and space, and the history of constraining motion progress of conditions. Accordingly, using the dialectical viewpoint to study physics is the only one of scientific thinking methods.

**Spirit is the supreme form of matter motion.**

This metabolism of matter exchange is illustrated as the basic form from existence of matter (the principle of metabolism of whirl vortex body, of roles of recoil force, gravitational field, electrical field, magnetic field, the principle of electromagnetic inter-compensation, of roles of temperature field, of existence and evolution of whirl vortex body), wave mechanics, electrodynamics, motion medium electrodynamics, trajectory of stars, thermodynamics, to origin and evolution of life, thinking morphology and supernatural sense. It is his most impressive dead to rederive the theory of relativity in the simplest mathematical manner: finding the average (I met him once in Xi'an and hence got the book).

There are reasons why he failed to implement his hypothesis: he is not specialized in physics, or the specialized experts are too materialistic. In fact, our physics uncovers no more than ten percent of matter scientists are aware of, said by one of the most famous experts in the world, a Nobel Prize winner. Furthermore, what we are studying (including this hypothesis the author believes) are merely infant sciences comparative to those in the Pure Lands, said by one of the most famous Buddhist masters living in the world.

## 5. Representing the metabolism

It will be one of the toughest tasks through out the world to maintain feasible representations. However it is worthwhile to launch the discussion.

For an agent, we should distinguish between the subagents (e.g., sub elements, sub sets, etc., assuming that the system is decomposable.): what is growing and what is anti-growing; what (aspects) in the environment contribute to the growth and anti-growth respectively and the limits of each.

> For example, how can one pass an exam. Provided that he is an excellent student and his teacher is really good. This is his positive aspect in chief. However, he is very easily interfered and distracted by things outside, e.g., friends, games, music, movie, TV, etc. So the problem is: there are two sources (two kinds) in his metabolism of knowledge, one is active and the other is anti-active. Now we should focus on these contradictory aspects: positive and anti-positive.
>
> First, we should analyze the positive aspect by means of learning: whether the teacher improved his teaching, or whether the student puts more weight on his study, in a dynamic figure.
>
> Second, analyze the anti-positive aspect: whether he touches less the disturbing source, or whether he puts less weight on this aspect, also in dynamic figure.
>
> Third, simulate the neutral figure or curve by comparison, and/or through learning examples in history, as to predict the future behavior.
>
> This is the general clue presented in my paper Dynamic Modeling of Multidimensional Logic in Multiagent Environment (Liu [3]), where "the learning curve" is just what I mean by neutral representation.

This could probably be the conjunction point between neutrosophy and neural science.

As conclusion, we would very likely base neutrosophy on dynamic learning background, just as an idiom says: it is never too old to learn.

## 6. The final remark

I have just presented my discussion on neutrosophy in time domain. However we may ask whether there is neutrality between space and time, or whether there is a relativity theory that unifies time and space.
- Limited to my knowledge so far, I can present the only one I know, it is Buddhism — which combines time and space domains into a single "holograph", and combines the universe in it as well.
- Unfortunately, we normally cannot reach that point before we come to the Buddhist Pure Lands (as described in Buddhism, also mentioned in (Liu [4])), where everyone can see/foresee any moment in any life cycle of every being in the universe.
- Far more answers can be found in the encyclopedia mentioned in section 2.

# A Logic of Complex Values

Chris Lucas
CALResCo Group
52 Mount Road, Middleton, Manchester M24 1DZ, U.K.
clucas@calresco.org
http://www.calresco.org/

*Abstract*

Our world is run by a logic that has no room for values, by a scientific methodology that disdains the very notion. In this paper we try to redress the balance, extracting many modern scientific findings and forms of philosophical reasoning from the field of complex systems, to show that values can and should be made part of an enhanced normative logic derived from Neutrosophy. This can then be employed to quantitatively evaluate our beliefs based on their dynamic effects on a full set of human values.

**Keywords and Phrases:** *Complex Systems, Axiology, Neutrosophic Logic, Intrinsic Values, Synergy, Dynamical Fitness, Attractors, Connectivity, Holarchy, Teleology, Agents*

**1. Introduction**

As we move into the 21st Century it is opportune to take a look at where we are as humans and where we are going. The achievements of the 20th Century in material matters are clear, yet this period has done little if anything to improve humanity as a species, we still spend inordinate amounts of time on war and similar destructive practices. The cry from all quarters seems to be that our lifestyles are becoming unsustainable, that we have lost our values. Yet our world is run by a logic that has no room for values, by a scientific methodology that disdains the very notion. Here we will try to redress the balance, extracting many modern scientific findings and forms of philosophical reasoning from the field of complex systems, to show that values can and should be made part of an enhanced normative logic derived from Neutrosophy, which can then be employed to quantitatively evaluate our beliefs based on their dynamic effects on a full set of human values.

We start by looking in Section 2 at how we define values, and outline in Section 3 an existing science of values. Section 4 introduces complex systems science and we look in Section 5 at the need for a logic of wholes rather than parts. Section 6 looks at non-Aristotelian logics and what a logic of values would require, whilst Section 7 relates paradoxism to complexity notions. In Section 8 we look at the ideas of Neutrosophic logic and relate them in Section 9 to synergy and unpredictability. Finally Section 10 looks at what is still needed to fulfil the goals laid out in this introductory paper. An Appendix relates existing logics to the teleological fitness focus that we adopt.

## 2. Putting Values into Science

Like most scientists growing up in a 'value-free zone' it has taken much time for me to realise the self-deception involved in this stance. A recent paper [Lucas2000c] looks at this in detail, and outlines a metascience that includes explicit values to complement the implicit ones already included within the 'scientific method'.

Following on from an heavy involvement in the pursuit of complex systems science I realised that what we called 'fitness' in our experiments (the same 'reproductive' concept that is used in evolutionary biology) was actually a simplified overall value - survival. From this it seemed clear that this could be broken down into many related values (water, food, warmth etc.) all of which were necessary for survival or reproduction and all ignored in the one-dimensional biological practice of 'population genetics'. These form a complex of what I shall call 'primal' values and necessitate a teleological (agent driven) form of science that bridges the objective/subjective divide. In this view these 'values' are simply ends derived from our evolutionary past, and what we can call 'needs' are our genetic predispositions to actions (means) that can meet these ends, by enabling us to generate fitness enhancing trajectories through life. Thus we can regard 'values' as static descriptions of an organism's goals and 'needs' as their dynamic equivalents.

From this insight two others followed, firstly that the number of our needs or values increases with the complexity of the organism and with experience (adding flexibility to better cope with environmental change) - successive stages adding higher values - what I've called 'social' and 'abstract' value complexes (inspired by similar work by [Maslow] and others). Each again contains many dimensions or variables. Secondly, that all these needs need to balance dynamically, they do not exist in isolation but interact in nonlinear ways, both amongst themselves and with their environmental contexts (i.e. we cannot simply add their 'fitnesses' or fulfilments, they are epistatic).

## 3. Axiological Science

Looking around the Internet to identify work on values led me to discover Robert Hartman and his Science of Axiology. Study of this neglected area of science suggests that three types of valuing are possible. Firstly, we have the binary logic (classification) type of valuing which concentrates on existence or being (ontology). Here we use an either/or judgement, for example we ask is there a tree or not; is he an Arab or not ? This is the 'label' style most common to our descriptive science and philosophy, where we classify objects by lists of such traits or attributes. The second type of value is that of quantity, where we define parameters by size or ranking (cardinal or ordinal numbers), e.g. how heavy is it; is he taller ? This is the 'quantification' mode of valuing common to mathematical science and economics, and is an extrinsic or external (relative fitness) approach. For both of these 'objective' dimensions we assume (implicitly) a subject (or group of subjects) that are 'doing' the valuing of something outside themselves (which may, as in psychology, be another person's 'subjective' experience temporarily regarded as an 'object' to be studied).

The third type of valuing is that currently outside most science, and that is the valuing of wholes, or uniqueness. Here we treat all the values as a unique set of self-supporting attributes that contribute to a one-off absolute fitness or 'value-in-itself'. This is the valuing we see in such approaches as love or art, where the subject is not broken down into parts but ideally is accepted as an self-justifying whole. This is an intrinsic or internal perspective and it is this type of overall viewpoint that we will concentrate on here in trying to adapt Neutrosophy to valuation, where we will go beyond 'subject' and 'object' perspectives entirely and identify ourselves with the system from the inside as it were - we 'become' the system and can then look towards 'our' self-development as a contextually situated or 'embodied' entity.

The first thing that strikes about the [Hartman] approach is that it is an attempt to add logic to values by using forms of [Cantor]'s transfinite sets (Aleph$_0$, Aleph$_1$ & Aleph$_2$). This is ingenious but creates a number of problems which recent work by axiologists such as [Edwards], [Forrest] and [Moore] tries to address, the latter by adopting a finite version of axiology based upon quantum theory. One of the focuses of such work concentrates on ranking combinations of the 3 value dimensions, for example where we employ intrinsic valuations of systemic entities (e.g. love of logic), which is denoted as $S^I$, or systemic disvaluations of extrinsic entities (e.g. excluding the poor), denoted $E_S$. There are 18 combinations of these, the positive aspects are called 'compositions' by Hartman, the negative aspects 'transpositions'. Another problem is in differentiating between partial and full valuations, i.e. where we have 3 values met in a system containing only 3 values, versus 3 values met in a system which has 5 (denoted here as 3in3 v 3in5). Whilst many aspects of this tradition will be included in what follows, we shall instead adopt a different perspective based upon recent developments in the sciences of complex systems.

## 4. Dynamical Systems and Complexity

Over the last two decades of the 20th Century considerable work has taken place in bringing together and taking forward the study of complex systems, defined as systems comprising many autonomous parts and interacting in complex ways. The foundations of this transdisciplinary science go back to work in cybernetics (feedback and homeostasis) and general systems theory by people such as [Wiener], [Ashby] and [von Bertalanffy] in the mid 20th Century, and it was later influenced by work on 'the pattern that connects' by [Bateson]. Following the advent of inexpensive computers this has joined with mathematical dynamical systems theory to form a science of complex systems in their own right, with many topologically isomorphic specialisms appearing (e.g. neural networks, cellular automata, artificial life, evolutionary computation, production systems), for an overview see [Lewin] or [Waldrop]. All these areas consider non-equilibrium systems, where their histories and future trajectories through the possibilities open to them become important, and this demands a more dynamic form of time-critical logic, a logic of change.

Analysing such systems makes use of three critical concepts of complex systems, the first is 'state space' which defines all the possible combinations of the system variables, e.g. 4 binary variables can combine in 16 ways, so state space (the possibilities) here comprises 16 points. It can be seen that this space escalates exponentially with both the number of variables and number of states available for each variable. The second concept is that of an 'attractor', and here the connectivity of the system causes a subset of state space to be preferred, the system will self-organize over time (circular causation) to concentrate on a small area of state space which we call the attractor. Within this perspective there are three types of attractor: point, cyclic and strange (or chaotic). Simplifying a little, the first assumes a static ontology (systemic value); the second a changeable one, one dimension plotted against time (extrinsic value); and the third occupies all the dimensions - system wide (intrinsic value).

In highly complex (high dimensional) systems however there will be multiple attractors present and this brings in the third concept, the idea of 'edge-of-chaos'. Here the system is found to move spontaneously away from either stability or chaos towards a dynamically semi-stable intermediate state (equivalent to the phase-boundary in physics) which comprises a power-law (fractal) distribution of both temporal fluctuations and spatial structure. Analysis of this unstable boundary needs some mathematical sophistication, yet mathematics and complex systems enjoy a difficult relationship. In truth there is no sign as yet of a maths of complex systems that can adequately deal with such complications, just a number of partial mappings that treat abstracted aspects of the whole. I look at this in my introduction "Quantifying Complexity Theory" [Lucas2000b]. In these scenarios the specific connectivity proves crucial, and we can adopt many perspectives, from simple on/off (systemic or Boolean), through weighted (extrinsic or neural) to integral (intrinsic or chaotic). In complex systems science we find that a middle connectivity (generating the 'edge-of-chaos') is needed to obtain the maximum fitness from such value combinations [Kauffman].

## 5. The Need for an Intrinsic Logic

By equating values and complexity here I suggest that a maths for one may also form a maths for the other, so we can attempt to evolve such a synthesis. From a psychological standpoint we must I think discard infinities, based upon the practical problems for humans in using such concepts computationally (and also in the difficulty in getting people to relate to such mathematical paradoxes). So what have we left ? The Hartman 3 value stages is a good start. Systemic values can be regarded as philosophical abstractions, my third value complex. Extrinsic values could be regarded as socially driven scientific measurements, my second value complex, whilst Intrinsic values could be regarded as primal values - the survival of the 'organism-in-itself'. It is interesting to note that this way of looking at the matter reverses the standard idea in intellectual circles that 'logic' is at the top of the 'pyramid' of our faculties. We see instead that by intellectualising we remove value, firstly by neglecting the whole for the part, then by reifying the categories into an either/or logic which goes on to form the basis of our system of social values embodied in a legalism of dualist 'right' and 'wrong'. We invalidly reduce intrinsic values to the lower systemic type in many aspects of our daily lives, as seen in our one-dimensional prejudices and our behaviours of conflict and competition.

Given that logic is so highly regarded academically, one way out of this conundrum would be to adopt an intrinsic form of logic, in other words one that can accept extrinsic and systemic logics as special cases. But can we find one ? Fuzzy Logic, which previously many have been inclined to consider, copes well with extrinsic values, we can specify a half apple logically FL(0.5,0.5), but it seems to deal all too poorly with value combinations, a half AND a half = a half (using Zadeh's intersection operator), i.e. a similar problem to that found with adding infinities. Moore makes the point that fuzzy logics require time, unlike systemic logics propositions are not universally 'true', logical values change over time (a fresh apple FL(1,0) $\rightarrow$ eaten apple FL(0,1) ). I'd add the idea that intrinsic logic requires context, it is dependent upon all the other values present, thus needs to be an interactional form of logic, a connectionist logic. In complex systems the idea of 'uniqueness' relates to our current specific position on the 'fitness landscape', in other words it is dependent upon our history. This can be said to comprise a number of systemic dimensions (axes) and a number of displacements (vectors), but the 'now' is defined by how these all interact - what I would here call intrinsic fitness. This network approach has significant advantages, since the number of extrinsic paths possible through the network (the set of possible value relations) grows exponentially with size and connectivity - a true measure of such intrinsic value would therefore approach infinity (especially in humans), as desired by intuitive axiological approaches to the value of a human life.

A further complication however relates to the concept of synergy, the idea that the whole is more than the sum of the parts (I look at this in more detail along with fitness in my introduction "Fitness and Synergy" [Lucas2000a]) - and this is I think precisely what we are really looking for here. It relates to the complex systems notion of emergence, the creating of new global properties (values) by the combination of parts. Thus no matter how many extrinsic items you add, you don't get the emergent next stage unless you go beyond aggregation and create suitable connectivity - a value add step. This idea also allows us to compare intrinsic sets, a more 'developed' person will have more values ('higher' needs or more discriminative ones), avoiding the ant=human problem of equating the value of all life forms. But the potential v actual issue is important here also, we must take into account that experience and education can convert potential (as in a baby) into actual (an enlightened adult) and this possibility must also affect valuations.

Making synergy more specific perhaps, we can imagine a 'precision' axis where we can specify a value in terms of bits, from 1 bit (binary) to infinite bits (irrational number). We can also imagine a 'depth' axis where we can specify the number of values in the system (again from 1 bit to infinite). This gives us an integral plane, or map at a single level. But such interacting entities generate new emergent levels, so we have a further 'height' dimension corresponding to these extra layers (visualise, say, atoms, molecules, cells, organisms, societies, ecosystems, planet). We would regard this 'height' as comprising a further 'holarchic' dimension beyond the three Hartman ones (precision, depth, plane - which correspond to extrinsic, systemic and intrinsic values, but all relating to a single holon [Koestler]). The size (volume) of this 'box' perhaps corresponds to the overall value of the hypersystem [Baas] we are considering. Like the other solutions offered by Forrest and Moore however it has a snag, the general idea of network analysis (and specifically emergence) has proved to be mathematically intractable using all the normal techniques - most complexity problems are related to graph theory and mathematically tend to be NP-complete, insoluble in polynomial time e.g.[Crescenzi & Kann].

## 6. Beyond Aristotle

For over 2000 years we have used and taught almost exclusively the classical logic of Aristotle (and its modern Boolean equivalent), in which all issues must be 100% true or 100% false and there is no other possibility. That this causes immense problems when applied to humans has been known for many decades. [Korsybski] identified this with the confusion of 'map' and 'territory', this means that we try to force a limited model onto an unlimited reality. From the discussion in this paper on values, we can see that this relates to forcing a 1 bit systemic value (one dimensional) onto a multi-bit extrinsic value (two dimensional) to which it cannot relate, and even worst to forcing the same dualist evaluations onto intrinsically valued humans (four dimensional), a 2 level 'category' or 'type' error - showing the need for a 'higher-order logic' of at least 3 levels, a mathematical meta-model [Palmer]. Given our general pre-occupation with Aristotelian logic, it is sobering to discover just how many non-standard logics already exist e.g. [Suber], so there are a number of less familiar possibilities that might be explored to find a suitable method for a logic of values. And what would we want from such a logic ? I'd suggest the following at least (taking into account Moore's criticisms of Hartman):

**a) Evaluate to a higher total value the more values (dimensions) exist, i.e. complexity matters**

**b) Value intrinsic systems more than extrinsic variables, and those more than systemic distinctions**

**c) Discard 'Law of Excluded Middle' - which prevents us specifying fuzzy truths (extrinsic values)**

**d) Provide adequate resolution to deal with real variables, the full variety encountered in life**

**e) Include a method of treating the 'many' as of higher/lower value than the 'one' (aggregation)**

**f) Allow for synergy, i.e. A + B can generate an emergent higher value C, a value-add step**

**g) Differentiate between possibility, probability and actuality, i.e. future choice & past history**

**h) Be context specific, i.e. allow truth to depend on time, space and interactions/connectivity**

**i) Differentiate positive-sum & negative-sum trajectories, i.e. dynamical fitness effects**

**j) Give an intuitively adequate rank ordering for $S_I$, $I^E$ and the rest of the combinations**

k) Add a measure of 'fulfilment' or personal development, i.e. 3in3 v 3in5, actuality v possibility

l) Allow for circular causality, the multiple interconnected paths of real systems thinking

m) Allow for obligations, the idea that we should not degrade the values of others (morality)

This list seems to go beyond most forms of logic and includes aspects of many types of logic whose implications and technicalities are a specialist task to unravel (see the Appendix for a look at how these relate to our teleologically based fitness viewpoint). Only a few attempts have yet been made to try to combine fuzzy thinking and the more teleologically oriented logics, e.g. [Gounder & Esterline], which brings us perhaps to a novel type of paraconsistent logic (for a quick overview of these see: [Priest & Tanaka]). The one I have been looking at especially is Neutrosophic Logic [Smarandache] which is unique in that it has three axes of logical validity. One is 'truth', one 'falsity' and one 'indeterminacy'. Now of course the latter immediately suggests a role for quantum theory, and also allows for those paradoxes and contradictions that troubled Frege, Russell and Gödel. Additionally in this logic we need not have normalised values (i.e. 0 to 1) we can have 1 AND 1 giving 3 (or anything else), thus synergy seems possible, i.e. the generation of new niches, new opportunities or alternatives. This logic was intended to bridge the gap between literature/arts and science, so is already in the same area as we are considering here, and discusses multiple-value sentences and ways of distinguishing between relative and absolute truths.

### 7. A Philosophy of Stress

Before I consider this as a logic, perhaps I'd better say something about Neutrosophy as a philosophy. The creator, mathematician Florentin Smarandache, was something of an anarchist, a Romanian fighting against the repressive communist regime of Ceausescu in the 1980's. Living a 'double-life' (the 'spin' culture of deceit now familiar to us all) helped him to recognise paradox as crucial, so he came up with a philosophy in which one could prove anything - and also disprove it ! He applied it widely to highlight contradictions - combinations of opposites in stress, and founded the literary movement known as 'paradoxism'. Despite the nihilism suggested, this does have much in common with spiritual ideas (the figure/ground or Yin/Yang) and with complexity science (where we balance static conscious 'rational' order and dynamic unconscious 'irrational' chaos), and so realise that as Smarandache said, "constants aren't and variables won't" - the two descriptions are contextual or transient [Lucas1997]). For humans, if we are too static then we stagnate and die, if we are too dynamic then we disintegrate and die, paradoxically we must be both somewhat ordered to survive and somewhat chaotic to grow. To be human is thus to be indeterminate, to live a contradiction. In an insight from Eastern philosophy, we are not 'either' order 'or' chaos, but 'both' and 'neither'.

We could regard these two axes (of chaos and of order) as those of 'generalisation' (artistic scope) - where we encompass everything but make no distinctions (mystical awareness or intrinsic value perhaps), and 'specialisation' (scientific content) - where we make 'cuts' or valuations across infinite reality, this would relate in dynamical systems terms to taking Poincaré sections (low dimensional projections from a higher dimensional whole). Thus both width and depth can be included, but not at the same time again echoing quantum complementarity and granularity [Smith & Brogaard] - we can see either the whole (dynamic wave) or the part (static particle). We can view the move from 'indeterminacy' to 'true/false' as the making of distinctions, the creation of opposite pairs or dualisms, i.e. systemic values (something akin to [Spencer-Brown]'s 'Laws of Form'), but each such division must exclude all the others in either/or logic. Thus our very act of classifying the world generates its own stresses, a problem not unknown even within conventional science e.g. [Kuhn], where new paradigms or syntheses are occasionally necessary to transcend the tensions of suppressed inconsistencies and contradictions.

## 8. Contextual Neutrosophic Logic

Neutrosophic logic itself allows <A>, <Not-A>, <Anti-A> and <Neut-A>. The first two are standard Aristotelian, <Anti-A> is Hegelian (included in <Not-A>) whereas <Neut-A> includes all the other possibilities, i.e. the set of distinctions ignored when looking at opposites (e.g. if <A> is 'white', <Anti-A> is 'black', <Neut-A> includes blue, red, yellow etc., <Not-A> is <Neut-A> + <Anti-A>). The values however are neither binary nor fuzzy but are intervals, allowing vagueness (e.g. it could be 30-40% 'white', 10-20% 'black' and 40-60% 'unspecified'). Another idea included is that of Multispace, where a set M of structure $S_1$ is said to contain also many subsets with different structures $S_2..S_k$ not included in $S_1$ - a sort of fractal hierarchy similar both to the layers mentioned earlier and to the structure at the 'edge-of-chaos'.

One of the main tenets of this form of logic is that for any combination of the three dimensions NL(T, I, F) a context or 'referential system' can be generated to make the statements valid. Thus 'truth' can not be applied to all possible worlds, and whether any statement is 'true' in our human world becomes an empirical matter and not an issue of logical analysis. This idea allows us I think to effectively distinguish between intrinsic and extrinsic/systemic value schemes, in that the set of worlds in which a value is 'true' changes with complexity, i.e. context. Within any intrinsic system, any extrinsic value or systemic distinction will fail 'truth' in many frames of reference, whilst the intrinsic value of 'existence-in-itself' will still hold true. For example, the exact systemic statement 'I see the clock showing 12:00" would fail to be truth a minute later, the fuzzy extrinsic truth "I see the time" may hold for many hours, while the intrinsic value "I see" should hold true for all my life. Thus we naturally perhaps can justify higher truth values logically both for intrinsic values v extrinsic and for extrinsic v systemic, if we include domain-specific temporal and (state) spatial context.

In this logic a systemic distinction (a division of the world into system/environment or figure/ground) has a value of one bit, no more, no less - either 'in' NL(1,z,0) or 'out' NL(0,z,1), where z relates to all the undifferentiated content of the two halves. We can go on to make more distinctions, more cuts through the whole. In the limit we obtain a binary set, corresponding to the number of distinctions made, infinite if we wish. An extrinsic value however, a one dimensional measurement, has a variable number of bits of precision - a vague value NL(x,z,y,) where x+y is the range and z includes the measurement uncertainty. In a Fuzzy Logic reduction z disappears and x+y normalises to 1. Again we can make more measurements within our whole, we obtain then a set of reals. From this perspective we see that systemic values are simply low resolution extrinsic ones, crude value distinctions that discard the precision dimension. It may perhaps be quite reasonable to equate, say, 8 systemic values with one extrinsic value of 8 bit resolution (within a linear viewpoint).

Now we take a further step, we group distinctions. We make associations between variables, we create an algebraic matrix, a mesh or network of interactions - an intrinsic system. Given that every systemic is a 1 bit extrinsic, this is a matrix of extrinsics in the limit. Again we see a new perspective, in that extrinsics are just crude intrinsics, low resolution views that discard most of the connectivity effects (the two topological dimensions), i.e. how values interrelate to support or oppose each other (the 'higher-order' causal terms usually ignored in science). To evaluate this stage logically we perhaps may usefully employ a fuzzy matrix logic [Yamauchi] but using neutrosophic triples. However due to the nonlinear and nonadditive nature of such value interactions we cannot adopt a simple global mathematics, applying standard matrix operators to the array. Each intersection pair now may require an individual connective or 'transition function', a local 'law' - reminiscent of a spreadsheet mode (because of the circular causality inherent in complex systems this would then 'hunt' for a solution - that attractor representing the output triple or intrinsic value). It is easy to envisage experimental changes to this array in the search for fitter dynamical solutions. This contextual perspective has much in common with the 'constrained generating procedure' form of emergence pioneered by [Holland], which extends our treatment into more general areas by implying that we must we formulate a logic that can generate further triples, i.e. be creative.

Thus we add another stage, the matrix in systems terms (if sufficiently complex) gives rise to emergent properties, a higher level of system, so we have a fourth value dimension, the hyperset of systems, which I earlier called an 'holarchic' value level - a nested heterarchy of intrinsics made up of systems or 'holons'. To take an example, in a rainforest the 'systems' of the geologist, botanist, zoologist, artist and mystic will all see (and value) different environmental 'systems'. These may be disjoint (if the experts are too single-minded) or may overlap, some may be more complex than others, they may differ in 'zoom' ratio (scope in space or time). This is the realm of combinatorics, where everything can be permutated from the set of wholes, factorial combinations of intrinsic modules at many levels. Evaluating this whole obviously causes immense practical difficulties, but we can of course treat relevant subsets as necessary (if we can identify them). Thus 'sustainability' in environmental terms would be such an holarchic valuation. In these cases we need to move from a 2D (plane) logical matrix to a 3D one which includes the various levels, a cubic matrix logic of neutrosophic triples seems required

## 9. Dynamical Beliefs

Looking more at the dynamics, we can relate this to tensegrity, the system of balance proposed for collections of interacting elements by [Buckminster Fuller] in Synergetics (700.00). Here a tension between a continuous 'pull' to truth (an 'attractor' in complex systems terms) and a discontinuous 'push' to falsity (a chaotic move or perturbation in those terms) relates to a balance between convergence and divergence - our edge-of-chaos, or semi-stable state. The indeterminacy relates then to the uncertainty as to whether a change will create or destroy value (or have no effect) - the 'butterfly effect' familiar in chaos theory [Gleick] - this also can include stochastic effects and measurement uncertainties. Trajectories can move in two directions therefore, which Smarandache relates to 'underhuman' and 'superhuman' behaviours - what in my terms I'd call 'dysergy' (negative-sum or sub-animal) and 'synergy' (positive-sum or full human), both of course relating to Hartman's idea of transposition and composition. This relates also to the idea of cancellation and reinforcement of waves in quantum and electromagnetic theory, and we can thus also regard values as being potentially 'in-phase' or 'out-of-phase' with each other. In more general value terms we can say that the three neutrosophic axes correspond to values that are 'good' (positive-effects and thus 'true' beliefs), 'bad' (negative-effects and thus 'false' beliefs) and 'groundless' (unpredictable and thus 'careless' beliefs) - each with respect to a particular situation.

The relation of beliefs to values is a subtle one. We each have a worldview in which we believe certain actions to be advantageous to us whilst others are not, and we tend to reify those theories that have in the past proved advantageous in meeting our needs as 'true', whilst those that have failed as 'false'. But we all have different experiences and contexts, so there is always a tension between these two poles, what I see as 'true' you may see as 'false'. This relates to our often limited vision, and to correctly evaluate the trajectories of an action based upon our beliefs we must take into account how our worldview meshes with those of the other entities with whom we interact. Bringing together our common beliefs in 'truth' generates what we call science, a consensus as to which theories have been tested as being generally effective. Yet even here we make errors, we do not take into account the full range of interactions involved, we reduce the intrinsic whole to extrinsic slices - just as we do individually. It is for this reason that we need to formalise a science of values, generating a logic of interactions that can identify where our narrower beliefs (scientific or more general) fail to be true in terms of the whole.

In such cases the values relate to the hypersystem. i.e. subsets of the whole system may have the opposite form, e.g. for a system of 100 people (simplifying each here to just one systemic value), what I believe is good for me NL(1,0,0) may be bad for you and for 8 others NL(0,0,1) and neutral for the remaining 90 NL(0,1,0). So, for the matrixed hypersystem of interacting values, the overall result would evaluate as NL(0.01, 0.9, 0.09), assuming standard arithmetic connectives (summing over a simple diagonal matrix of triples) and normalising. We see that in overall utilitarian terms the result is 8% negative, thus taking everything into account my belief is proved intrinsically false even if it was extrinsically true. We can also see the relative effects of this action, in that 90% of our social group are unaffected, thus we are not inclined to escalate the issue out of all proportion - a major problem in logics based upon only two axes. Note however that this simple example is unidirectional, it takes account only of the effects of my belief on the group, it doesn't include the effect of the beliefs of the other 99 members of the group on me or on each other. Given better knowledge of what all our values are and how they all interact, we can in principle derive a resultant fitness trajectory for the whole. This applies equally if the whole is just me and the parts are 100 different personal values. For more complex nonlinear connectives the same principle holds, although the mathematical difficulties will of course increase considerably.

## 10. Towards a Value Formalism

It is not our intention here to outline a fully working model, simply to establish the feasibility of so doing. In this section we largely follow [Krivov] (to whom the reader is referred for more technical details) in his attempt to generate a logic based general systems theory for multi-agent modelling. We adapt those ideas here to our neutrosophic value focus. We take as our starting point the definition, analogous both with the classical definition of a Model in predicate calculus and von Bertalanffy's definition of a System, of <Values, Connections>, in other words we have a set of values (function or process) plus a set of interconnections or relations (structure). Our 'agents' originate internal states or goals which must be taken into account, indeed this is our definition of needs - our drives to meet a set of internal values. To include such a teleology we add to our Finite Protocol Language modal operators and time, i.e. Operator(ValueComplex, Time, State) where Operator can be such as needs, prefers, believes, acts etc. and State is the status of the value complex (true/false for systemic, variable for extrinsic, set of included extrinsics for intrinsic, set of intrinsics for holarchic). A Model of the system contains a function stating how the structure and needs will evolve over time, i.e. $M(t+1) = F(M(t))$, this we refer to as a 'logic machine' (a predicate automata) and it incorporates the connectives or quantifiers that relate values to each other dynamically. The Relations thus have the form Relation($Value_1$, $Value_2$, Affect) where Affect is a connective that relates the affect on $Value_1$ of $Value_2$, and the overall function is our matrix.

To clarify the real world systems that we are modelling, we have the following progression for any 'agent':

**Value** → **Need** (for change in that value) → **Preference** (ranking of alternatives) → **Belief** (fulfilment theory) → **Action** (environmental output) → **Reaction** (feedback) → **Update** (belief and need changes).

All agents (and values) may of course act simultaneously, so our model is a multiobjective constraint satisfaction problem. The ranking of preferences can cause its own problems in the making of decisions [Ha], especially within interconnected nonlinear systems where 'ceteris paribus' (all else being equal) rarely holds. This brings in logical implication in that some values imply others to some extent (e.g. the ability to 'philosophise' implies meeting our primal needs), so a fuzzy implication operator is required. The timescales for the evolution of the various components (values, beliefs, preferences, needs, actions) vary, so it may be possible to model these separately if we wish (given computational resource constraints). Additionally, important aspects of preferences and beliefs, as well as actions, are determined by social norms, which brings in higher-level obligations and canalysation of state space (i.e. [Campbell]'s downward causation based constraints on alternatives, e.g. laws), and also by environmental issues (resource availability, costs). It is apparent that needs, preferences and beliefs can all be fuzzy variables and that we can have considerable uncertainty as to their standing. This aspect relates to the standard logical notion of "for all x" ($\forall$x), where uncertainty is zero, through "there exists" or "some x" ($\exists$x) which has variable uncertainty, to complete undecidability (which we may call Ix). Thus the neutrosophic axis I defines the improbability that our T/F axes are correct, in other words the truth value of the believability of our assertions is B(1-I), which relates to the approach taken in the k-calculus in qualitative decision theory [Ha].

Restricting ourselves just to values, there are in neutrosophic logic many (possibly infinite) forms of definable connective which leaves the possibility of finding a definitional set that matches what we wish to achieve with values. This is complicated to do in formalised logical notation (especially given the open ended set of possible systems, values and nonlinear interactions as here), but we can attempt to do so in more general terms (it may be possible, more formally, to use genetic algorithms to search the space of possible connectives to locate the optimum definitions, given an adequately defined set of goals). I'm primarily interested here however in looking at intuitively simple ways of combining multiple vectors, i.e. n-Tuples - extrinsic values of the form $E_i$(T,I,F) for i in the set 1..k. If we are to get more from less, i.e. synergy (or conversely less from more, i.e. dysergy) then we need it seems a multiplicative form of connective. For a simplified 2-value case a form that appeals is R = S(A AND B) where R is the resultant value and S is a synergy operator which can vary from zero to plus infinity (or take more nonlinear forms). This allows for both cancellation (S = 0) and reinforcement (S = +2), matching Moore's quantum wave theory, but also allows for other S values for more generality. The AND is our normal logical connective, defined in whatever way we choose for fuzzy truth values. This can easily be extended to cover the multivalued case, and if necessary we can have separate synergy operators for each interacting value, e.g. for 3 values: R = ($S_A$A AND $S_B$B AND $S_C$C ) etc. If we assume normal additivity (e.g. 2 values are greater than 1) then we can generate a fuzzy truth table as follows (A and B both assumed to be 1 here).

| S | R | Comment |
|---|---|---|
| 2 | 4 | Synergy, positive-sum, increase factor 2 (100%) |
| $\sqrt{2}$ | 2.818 | Partial augmentation, 90° in phase |
| 1 | 2 | Aggregation of values, zero-sum, standard maths |
| $1/\sqrt{2}$ | 1.414 | Diminution, partially out of phase |
| $1/\sum v_i$ | 1 | Normalised disvaluation, classical logic |
| 0 | 0 | Dysergy, negative-sum, 180° phase cancellation |

Two problems immediately spring to mind, firstly how do we generalise this to neutrosophic triples ? There may need to be interchanges between the T, F & I axes as we vary S, since this seems to convert between T and F. Secondly, can we generalise further to allow for the interaction of a number of input triples to result in the generation of a number of new output triples - as needed for the emergence of new values and levels ? Here we may possibly make use of [Stern]'s Matrix Logic, which uses two by two truth table operators comprising the values true, false, both (synergy) and neither (dysergy) spanning the logical levels of scalar, vector and matrix, and capable of generating autopoietic emergent systems. By generalising these values to fuzzy values and adding the indeterminacy axis we naturally seem to end up with a 3 x 3 neutrosophic matrix logic operator.

It will be noted that our approach to logic isn't static, we move through time - either discretely or continuously, and to cover the increasing generality we need at least a 4th-order predicate calculus of triples. Our viewpoint throughout blurs the distinction between logic and mathematics, and sees logic as simply a 1 bit version of mathematics, whilst mathematics is an infinite bit version of logic. Whether the two can be successfully merged dynamically with values, using neutrosophic logic, remains to be seen.

## 11. Conclusion

In this paper we have looked at bringing together three spheres of intellectual activity, firstly the idea of values or axiology, secondly the field of complex systems science and thirdly the area of non-Aristotelian logic. We have examined a number of connections between these fields and can conclude that a new form of paraconsistent logic (Neutrosophic Logic) may prove instrumental in welding together the needs of interactive humans whilst getting to grips with contradictory values and unpredictable complexity. There is much to be done in formalising in detail how this would work in practice, and in taking into account the complications added by evolution and contextuality, within a framework of circular causality and emergence. But the indications are promising that this could be achieved, given sufficient time and expertise.

With the availability of inexpensive computers and the growth in the use of multi-agent simulations we are now perhaps in a position where we can instigate an experimental form of normative logic, looking to use agent evolution to develop and evaluate various logical formalisms. In such systems the agents are generally taken to be autonomous (not globally controlled), teleological (having internal goals), contextual (interacting with other agents and their environment), heterogeneous (different from each other) and autopoietic (self-perpetuating). Whilst much work is yet needed to put such approaches on a firm footing, there are considerable commonalties between the stance taken in this paper and recent work within these areas.

By using a predicate calculus approach common to recent work on formalising complex systems, and both generalising this to a higher-order form suitable for treating multiple levels and including neutrosophic triples as primitives, we obtain a methodology of considerable scope, very suitable for use computationally and potentially applicable far beyond the area of values which has been our main concern in this paper. We should not minimise the difficulties however, we subsume here in our generalisation many specialist fields which all have their own share of difficulties and controversies. Netherthless, the need to better integrate disjoint mathematical technicalities with interconnected real life applications is clear, and to this end normative concepts can form a bridge that links these two worlds. One final observation is that given our susceptibility for error, in the recognition of values, in understanding their interactions, in defining the scope of our systems and in defining suitable connectives, then the adoption of a form of logic that permits uncertainties and supports paradox seems highly appropriate.

URL: http://www.cs.uregina.ca/~yyao/JCIS2000_inv/Yamauchi.ps

**Appendix: Relating Logics to Teleological Fitness**

In this appendix we look briefly at the various forms of logic that we think need to be integrated if we are to create a valid logic of values. We relate the main idea of each of these formalisms to our treatment of values in terms of teleological fitness. One of the key differences we should emphasise is that we replace impersonal logical generalisations by personalised contextual ones. Each statement is related to the viewpoint of a living organism (not always human), which must reason out (not always consciously) ways of meeting their needs within ongoing environmental situations. In human terms this external context is largely social, and so we understand the effects of our actions on our own needs and the fitness of our wider societies to be dependent upon the validity of our socially derived beliefs (we make no distinction formally between ethical or other type of value fitness). By bringing in the idea of personal agency we highlight the tension between global rules and local ones. In a specific context it is often the case that global rules are inapplicable (as they were designed for standardised situations which do not always hold) and local ones are in conflict (different values are mutually incompatible when taken together) so we may need to restrict the scope of our logical rules by imposing appropriate boundary conditions.

**Natural Deduction**: The definition of a logic without a rigid set of axioms. Our system has no axioms for generality (allowing us to use fuzzy truth values), all connectives are defined as local inference rules applied to arbitrary premises (scientific hypotheses). This open (sympoietic) contextual system contrasts with the closed (autopoietic) global systems of most formal treatments, and discards completeness for applicability.

**Classical Logic**: Truth value is derived from valid syntactic forms of argument. The premises are considered to be fixed truth values, and only Aristotelian truth values (1, 0) are permitted. We generalise this initially to allow fuzzy truth values for intersection (AND), union (OR) and the other connectives, and later add semantic considerations relating to wider values.

**Alethic Modal Logic**: Qualifies truth by adding possibility, actuality and necessity. The first we regard as encouraging creative alternatives, new paths through state space; the second denotes our current position in state space. The third implies that a value cannot exist or an act take place unless a condition is met (a critical path analogy).

**Deontic Logic**: Adds obligation, permission and forbidden operators. We regard the first as an historical social norm intended to avoid fitness reduction (which may be empirically invalid), the second as an allowable alternative (within the social structure) and the last as a denial of a freedom to pursue an alternative (due to implied socially unacceptable disvaluations). Each relation is contextual within a cultural worldview.

**Epistemic Logic**: Adds knowledge and beliefs. We adopt a coherentist approach to beliefs, based on the circular logic of complex systems, in which all beliefs support each other in a consistent worldview (but one grounded by empirical trial and error, so this is neither absolute nor relative). We can have three types of belief, 'true' - that the proposition will have positive fitness effects if acted upon, 'false' in that it will have negative effects if acted upon and 'indeterminate' where we don't know what result will pertain.

**Temporal Logic**: Adds future and past operators, which relates to values changing over time (e.g. with metabolism, experience and education) and validity being dependent upon the timescales involved. In our treatment time is regarded as a sequence of discrete steps.

**Dynamic Doxastic Logic**: Adds propositions and actions that implement changes in beliefs. We thus allow for both changes in base values and changes in the effectiveness of our beliefs in satisfying needs.

**Proairetic Logic**: Adds preferences, which we relate to rankings of values, and the rankings of alternatives in the achieving of them, plus the changes in ranking over time with actions. This logic includes utility measures, which we generalise to matrixed global fitnesses, incorporating nonlinear interactions of preferences.

**Quantum Logic**: Adds indeterminacy and probability, discarding the law of excluded middle. Measurement here relates to making a decision between alternatives, the act of choice transferring information (knowledge) from possibility to actuality - a trajectory through state space opening up a new set of resultant possibilities.

**Nonmonotonic Logic**: Allows new rules to restrict the validity of the existing system, permitting evolution of truth contexts. This allows us to reduce or increase the scope of our values and beliefs, depending upon new circumstances, and adds falsification of existing rules.

**Fuzzy Logic**: Adds partial fulfilment of truths, an infinite valued form of logic. We allow partial fulfilment of our values (e.g. what we called 3in5) and partial beliefs about them, i.e. partial set membership or completeness.

**Mereology**: Focuses on the relationship between parts and wholes (usually on a reductionist assumption that the whole equals the sum of the parts). Relates in our treatment to the difference between intrinsic systems and their component extrinsic and systemic values.

**Non-Adjunctive Logic**: Adds non-additivity, i.e. A AND B ≠ A & B. Here this allows contradictions of joined value systems and the possibility of dysergy and synergy (i.e. emergence - the whole not being equal to the sum of the parts).

**Higher-Order Logic**: Add propositions that act upon sets of predicates rather than atomic facts. Here it relates to intrinsic systems being defined as sets of extrinsic values, and to higher levels of value logic incorporating sets of intrinsics. Thus the 4 orders of logic necessary as a minimum are systemic, extrinsic, intrinsic and holarchic (although these are not regarded as disjoint but are overlapping frameworks).

**Neutrosophic Logic**: Adds independent true and false axes plus an indeterminacy one. We can have values and beliefs that have both true and false aspects (e.g. good and bad interactions with other sets of values), plus uncertainties as to their effects.

It can be seen that each of these logics expresses only a subset of life's possibilities, so it is not surprising that when applied to real human situations that problems occur and important value data is excluded. It is one purpose of our logic of values to highlight such problems, and to formulate a set of connectives for each context that minimises such reductions.

Note that since any analogue value can be expressed as a string of binary digits, it is in principle possible to generate a logic of values using standard Boolean logic. However the added complications of doing this, together with the already difficult nature of the task, suggest instead approaches that adopt more natural human ways of expressing comparative variability and reasoning. Additionally the advantages of having an axis of indeterminacy, with neutrosophic logic, allows us to keep in mind the partial nature (contextual incompleteness) of all logics, and the vast difference between infinite possibility space and finite actuality space.

# Comments to Neutrosophy

Carlos Gershenson
School of Cognitive and Computer Sciences
University of Sussex
Brighton, BN1 9QN, U. K.
C.Gershenson@sussex.ac.uk

*Abstract*

*Any system based on axioms is incomplete because the axioms cannot be proven from the system, just believed. But one system can be less-incomplete than other. Neutrosophy is less-incomplete than many other systems because it contains them. But this does not mean that it is finished, and it can always be improved. The comments presented here are an attempt to make Neutrosophy even less-incomplete. I argue that less-incomplete ideas are more useful, since we cannot perceive truth or falsity or indeterminacy independently of a context, and are therefore relative. Absolute being and relative being are defined. Also the "silly theorem problem" is posed, and its partial solution described. The issues arising from the incompleteness of our contexts are presented. We also note the relativity and dependance of logic to a context. We propose "metacontextuality" as a paradigm for containing as many contexts as we can, in order to be less-incomplete and discuss some possible consequences.*

## 1. Introduction

Upset because none of the logics I knew were able to handle contradictions and paradoxes, I created my own, which was able to do so (Gershenson 1998; 1999). Later I would find out that my Multidimensional Logic (MDL) fitted in the category of paraconsistent logics (Priest and Tanaka, 1996). Basically, I rejected the axiom of no contradiction, and instead of having truth values of propositions, I defined truth vectors, where each element was independent of the other. In this way, we can have a vector representing something that is true and false at the same time, or nor true nor false at the same time. MDL gave interesting results, and the new perspective influenced me to propose ideas in different branches of Philosophy, Complexity, and Computer Science (e.g. Gershenson, 1999; 2002).

In the autumn of 2001, I became aware of Neutrosophy (Smarandache, 1995). The fact that it reached similar results than the ones I had makes me think that there is a lower probability that we are completely mistaken. Neutrosophy covers a much wider area than my ideas, but nevertheless I believe it could be enriched by them, since there are non-overlapping parts. In this work I will expose how some of these ideas might help in making Neutrosophy less-incomplete.

The contents of this paper were mainly sparse comments made to Florentin Smarandache, and while giving them the shape of an article, the continuity of ideas was not an entire success. My apologies to the sensible reader.

## 2. To be *and* not to be...

When people speak about the being, since it is one of the most general things you can speak about, they often seem to speak about different things. I define two types of being: absolute and relative. The **absolute being** (a-being) is independent of the observer, infinite. The **relative being** (re-being) is dependant of the observer, therefore finite, and different in each individual. The re-being can approach as much *as we want to* towards the a-being, but it can never comprehend it. Objects a-are. Concepts, representations, and ideas re-are (Objects do not depend on the representations we have of them).

We can only suppose about what things a-are, we cannot be absolutely sure, we can only speculate, because they a-are infinite and we are not. We cannot say that something is absolutely true or false. We can only assert things in a relative way. We could assign truth values or vectors to them, but these would be relative to our **context**.

The **being** would be the conjunction of re-being and a-being. This is, without the distinction we are drawing between them. But isn't it confusing to speak about something which is absolute and relative at the same time? Yes, precisely. But that is what we do every day.

We can say that every effect a-has cause(s), but it does not mean that we can perceive them, so some effects do not re-have causes.

## 3. Reason, beliefs, and experience

Reason cannot prove itself. It would be as if a theorem would like to prove the axioms it is based upon. Reason, as all systems based in axioms, is incomplete (in the sense of Gödel (1931) and Turing (1936)). Beliefs are the axioms of reason and thought. Reason cannot prove the beliefs it is based upon. But how do these beliefs arise, then? Through experience. But experience needs of previous beliefs and reason to be assimilated, and reason also needs of experience to be formed, as beliefs need of reason as well. Beliefs, reason, and experience, are based upon each other. Which came first? All of them. We cannot have one without the others. **Contexts** are dynamic, and formed upon beliefs, reason, and experience. It is there where the re-being lies. Since the re-being is dependant of our context, it is also dependant of our beliefs, reasonings, and experiences. Contexts are dynamic because they are changed constantly as we have new experiences, change our beliefs, and our ways of reasoning.

For example, we cannot say if capitalism a-is good or bad (independent of a context). It re-is good for the people who get a profit from it, re-is neutral for the people who think that are not affected by it, and bad for people which suffer from it. There re-is a god for people who believe so, and there re-is not for people who do not believe so, but we cannot say if there a-is a god. We can speculate as much as we want to about everything which a-is, but we will never contain it, therefore we can only have *an idea* of it. It is not only that "we know without knowing" (Smarandache, 1995), but less-incompletely: we re-know, but a-know not.

Anything we want to assert has an implicit **context**, for which what we want to assert is consistent. But *every idea is valid in the context it was created*. Ideas cannot a-be right or wrong. They re-are right or wrong according to a context. And since in order to be created they need to be consistent with their context, they re-are always right according to their context. Then which idea is more valuable? We should try to see which context is more valuable. Well, we cannot say which a-is more valuable, we can only see that if a context contains others, the ideas created in it will be valid also in the other contexts, and not necessarily vice versa. We do not know which one a-is better, but the context which contains the others will be **less-incomplete**.

A problem arises on the horizon...

## 4. The silly theorem problem

For any silly theorem $T_s$, we can find at least one set of axioms such that $T_s$ is consistent with the system defined by the axioms. How do we know a theorem is not silly? Or how do we know if the axioms are not silly?

We can extend this same problem to contexts: there can re-be any silly idea $I_s$ so that there is at least one context for which it is consistent[1]. Empiricism comes to the rescue. Well, not only experience, but **evolution** more or less helps us get out of the problem. This is, contexts which support silly ideas, even when they re-are not silly in their own context, will not be able to contend with the contexts which are able to describe more closely (less-incompletely) what things a-are. The evolution of contexts consists in making them less-incomplete. In formal systems a similar evolutionary (tending to pragmatic) criterium could be seen: the axioms which support silly theorems are not useful, and that is why we then call them silly. If the axioms and the theorems derived from them are useful in a context, we do not consider them silly. But this is not a complete solution for the silly theorem problem. There is no complete solution for the silly theorem problem (Or by saying this *that* is the solution?).

But silly ideas (for our context) are useful as well, because our ideas become stronger when they neutralize them, and become less-incomplete if they contain them. There would not be smart ideas without silly ones.

We can take advantage of the silly theorem problem, in order to define axioms after the theorems, or contexts after ideas. For example, we might want such sets that A contains B, and B contains A, but without A and B being equal. We would just need to define the proper axioms. If this new set theory will be useful or silly is no reason for stopping, because we will not know if it a-is silly or not, and we will know if it re-is silly only after we create it.

---

[1] But the silliness is also relative to the context of the one who judges the idea...

## 5. Incomplete Language

Not only silliness, but all adjectives can only be <used|applied> in a relative way, dependant of a context. Language is relative as well. How can we speak about absolute being, then? We can and we cannot. We speak about it, but in that moment it a-is relative. For us, it is and it is not-incomplete. But that we cannot completely speak about it, it is not a reason to stop speaking about it (as Wittgenstein (1918) would early suggest in his Tractatus Logicus Philosophicus), because we can incompletely represent its completeness... As Wittgenstein himself (but not most of his followers...) realized, following the ideas in the Tractatus, we would not be able to speak about anything... (languages are incomplete). Language is used inside a context. Depending of this context the language will be different.

## 6. Incompleteness

If we cannot create complete systems, we should try to make them as less-incomplete as possible. Since our systems are incomplete because they are finite, there is no way of measuring the completeness of a system. It is like asking: how infinite is x, when x is a finite number? We could also have a huge system, but if it is not related to the a-being of something, it is not so useful as a small system which describes closely a part of what something a-is. So, we could say that a system **re-is** less-incomplete as it approaches more the a-being.

But, for example, what is the a-being of the number four? We define numbers, so we determine what they are. But this definition arises on our generalizations of what things a-are. There a-is no number four, but we could say that there the number four a-is in things.

The a-being is far from materialism. Materialism re-is, and we cannot say for sure if things a-are only materialistic or not, because we do not know what matter a-is.

So, getting back to the incompleteness of systems, they can approach the a-being as much as we want to, but we <do not|will never> know what the *a-being a-is*. We could measure the incompleteness of a system in a relative way, but I am not interested in defining such a method. But we can see that a system will be less-incomplete if it contains others, even the ones that (re-)"are wrong". This is because a system which tries not only to explain why things a-are, but also why other people thought it was something else will have a wider perspective. Because people a-are not mistaken (nor right). There cannot be errors inside their context (people do not create systems silly/mistaken for themselves). So, by containing as many contexts as possible, we are also containing as many systems as possible. The great attempts in science of unifying theories could be seen with these eyes. But for containing as many contexts as we can, we are required to leave all hope of non contradiction behind, but I prefer to do it for the sake of less-incompleteness. The non contradiction in systems is just a prejudice. The systems are incomplete, so the contradictions are caused by their own limits. Extending the system allows us to contain the paradoxes, and once we <understand|comprehend> paradoxes, they stop being contradictory.

Perhaps the first reason for entering happily the realm of paradoxes and contradictions is because we know that our systems are incomplete. This means that, even when they might be valid for our present context, as our context enlarges, our systems sooner or later will not be able to be consistent with all our context, as it has happened all through history. And also because we are aware that people who do not share crucial parts of our context, they will not agree with our systems. We are predicting the failure of our systems outside our context (as 1+1=2 in a decimal context, but 1+1=10 in a binary context) by perceiving its limits. This means that we are predicting that someone will say: "Your system is wrong" (related to his incompatible context). But by saying this, he proves that we were right. But people inside our context will say: "You (re)are right", which as our context evolves will prove to be "wrong" (or less-incompletely speaking, "less-incomplete"). This is a paradox. But this paradox makes us less-incomplete, because we are containing the contexts of the people who say "You are wrong". Our theory will re-be true and false at the same time (related to different contexts), but it will **always** re-be true and false at the same time, as compared with theories which now claim to "be" true, without admitting that outside their context they might re-be false.

We should not only seek for the truth or falsity or indeterminacy of an idea, because these are relative (absolute truth (and falsity) is unreachable (but we can approximate as much as we want to)) but for its **less-incompleteness**. This less-incompleteness is also relative, based on our beliefs (one could argue that a small context could be more close to the a-being than another one which contains it, but this cannot be discussed, just believed).

A less-incomplete theory should contain theories which are not wrong, but incomplete in our context (i.e. myths, religions, dogmas, etc...). To be against something is useful, but it is more useful to contain both (<A> and <Anti-A>'s). For example, a less-incomplete political theory should contain socialism and capitalism, despotism and anarchism, dictatorship and democracy, fascism and republic, communism and terrorism... At least, understanding why each one exists and because of what, and for what each one is more suitable; take the things WE need from them, reject the rest, add a bit of sugar, mix it for five minutes (or until the foam is assimilated), and you have another utopia! *If you know the rules of the game, you can change them...*

### 6.1. Neutrosophy and Incompleteness

Neutrosophy, as other systems, embraces the spirit just described: attempting to contain as most contexts as possible, even when they are contradictory, for attempting to be as less-incomplete as possible. But of course Neutrosophy itself predicts its own decay. It is not the *non plus ultra*. It is not finished[2]. But <we cannot|there is no need to> go further now because it fills successfully our contexts.

For example, we could add more values (concepts) to a (Neutrosophic) vector than True, Indeterminate, and False. We just need to define them... (e.g. {T, I, F, T∧I, T∧F, I∧F, T∧I∧F, ~[T|I|F]}) How useful it is this? Someone might ask the same about T, I and F. Why not only T & F, or only F. It depends on the context we are, on the things we need. Logics are just a tool. It depends on what we want to do that we need to make a different tool (of course, there is no "ultimate" tool). Since Neutrosophy is also based in axioms, it is also incomplete.

---

[2] Is there such a thing as a "finished idea"???

<A> is finite, but <Neut-A> is infinite. Therefore, the ideas will evolve infinitely. There are infinite <Anti-A>'s, each one related to a different context (Reference System in Smarandache (1995)). The values/ranges of truth, indeterminacy and falsity (T. I. F.), and are dynamic, relative to a context, and there are infinite contexts, so about any event, there are infinite number of Neutrosophic values/ranges, and any context a-is also infinite (it cannot be **completely** described, as well as the event) Which one is the more representative? The answer to this question is also related to a context! All the "answers" are related to a context, closed, finite... the interesting thing is that all questions seem to be open and infinite... so they can be answered in an infinitude of different (incomplete) ways.

If <A> is combined with <Anti-A>'s, let's suppose it evolves into <Neut-A>. But this is infinite, so it would not be completely <Neut-A>, but as noted by Smarandache (1995), <Neo-A>. Each step/cycle the idea is less-incomplete, but there will always be an infinitude of <Anti-A>'s, and an unreachable <Neut-A>... From a smaller context, <Neo-A> would look just as the old <A>?

## 7. Contextuality

Things re-are in dependence of their context. Since there is an infinitude of contexts, things can re-be in an infinitude of ways. It is only inside a specific context, which we can speak about the truth or falsity or indeterminacy of a proposition.

### 7.1. Context-dependant logic

Every proposition P can only have a truth value (or vectors) *in dependence* of a context C. This truth value/vector is *relative* to the context C. Propositions have only sense (in the sense of Frege (1892)) inside a context(s). Propositions have no sense without a context.

The truth values/vectors of a given proposition can change with context. So, for example "This proposition is false" has a value of 0.5 in Lukaciewicz logic, [1,1] in multidimensional logic, (1,1,1) in Neutrosophic Logic, and "?!" in Aristotelean logic. Or, the proposition "The king of France wears a wig" would be, in terms of multidimensional logic, nor true nor false ([0,0]) in the XX[th] century context, but true ([1,0]) in the context of the 1[st] of January of 1700.

We are just indicating the limits of logics. Logics are just tools. They are useful only inside a context. The context determines the logic. If a proposition goes beyond the context, the logic developed in the context will not be able to contain it (but not necessarily vice versa).

7.2. Derivations of contextuality

Since all <Anti-A>'s of <A> are related to a reference point, we can find all reference points so that all the elements of <Neut-A> are contradictory with <A> (and later with themselves...). This is, any element of <Neut-A> can be potentially <Anti-A>, you just need to have a reference point.

There will always a-be injustice, because this one is relative. Since different people have different contexts (or we can use the word Seelenzustand (soul state), to refer to the personal context, to distinguish from a general context)... So, since people have different

Seelenzustandes, we cannot speak of absolute justice, so things will be just for the people with power... The less-catastrophic panorama (and most naive...) would be that the people in the power would have the less-incomplete Seelenzustandes, trying to contain and understand as many Seelenzustandes as they can, so, if they are just, in spite their relativity, they will be just as well for all the people whose Seelenzustandes they contain.

If ideas are different in each Seelenzustand... well, they might have many similarities, but on the other hand there is the problem of language representing ideas.... one can quote the words of another in another context to communicate different ideas... but the "problem" of language is a different story... ("In Philosophy there are no problems, just opinions (like this one)")

All adjectives are relative. Thus, we can find "opposite" (related to a reference point) adjectives for the same object from different perspectives. (What is wrong from one perspective is right from another, what is beauty-ugliness, good-bad, complex-simple, complicated-simple, complicated-non-complicated, etc...). Are there adjectives which do not behave this way?

## 8. Metacontextuality

Following the ideas exposed above, we can argue that, in order to be less-incomplete, we need to strive for a *metacontextuality*, containing as many contexts as possible. We **believe** this is the only way we have to approach to the a-being more than we already have, but of course there might be other ways, and the proposed way will not be definite.

This does not mean that we will agree with "silly" ideas. This means that instead of just declaring them silly and forgetting about them, we will try to understand what led people to have such ideas, in order to try to see which perspective of the a-being they had. Then our perspective will be greater than if we would just ignore the ideas we do not agree with.

We are tied to a context, a relative one. This is because there re-are no basic elements from where we can build the rest of our world. There are only more complexity and indeterminism in the entrails of the subatomic particles and quarks. There a-is no essence in the universe, because everything is related. Everything a-is for and because of everything. Everything is based on everything. It is naive to try to justify the world once we are already on it. Everything is a condition of everything. Otherwise, it would not be AS IT IS.

Therefore everything can re-be seen in terms of everything (just as a Turing Machine (Turing, 1936) can represent all computations (and computations can represent Turing Machines)). There a-is no base. Everything can re-be a base.

Metacontextuality, as the one Neutrosophy and other currents strive for, by consequence predicts *tolerance*. This is, if we try to contain as much contexts as possible, we will be able to be less intolerant to contexts and Seelenzustandes that are neglected from absolutist non contradictory points of view. And this tolerance should be able to prevent many *conflicts*.

But may metacontextuality lead to *indifference*? This is, since everything can re-be depending on a specific context, does it matter which context we choose? I believe that this might be avoided if we put reason on its place, leaving place for experience and beliefs. Anyway this issue should not be disregarded.

## 9. Conclusions

I **believe** that the question "Do you **believe** in an absolute reality/truth?" is on the same level as "Do you believe in God?". This is, the question is completely metaphysical... And relative truths are incomplete. Should we keep searching for truth of things? I believe now it would be more useful to search for the **less-incompleteness** of things. This can be achieved by enlarging our contexts, containing as many contexts as we can.

Our limits are the ones we set to ourselves. We need only to take the blindfolding of our prejudices in order to attempt a **metacontextuality**. And Neutrosophy bravely does this. If the ideas exposed here, and the ones of Neutrosophy, are not assimilated by our society, this will be because the society does not need them. This is, it "functions" based precisely on the partial blindness of the individuals. Then, should we try to help everyone to open their eyes? We already are doing so, they will open their eyes only if they *want* to.

## 10. Acknowledgements

This work was supported in part by the Consejo Nacional de Ciencia y Tecnologia (CONACYT) of Mexico.

# CONTENTS

The Neutrosophy is a new branch of philosophy, introduced by Florentin Smarandache in 1980, which studies the origin, nature and scope of neutralities, as well as their interactions with different ideational spectra. Neutrosophy considers a proposition, theory, event, concept, or entity, <A> in relation to its opposite, <Anti-A>, and that which is not <A>, <Non-A>, and that which is neither <A> nor <Anti—A>, denoted by <Neut—A>. Neutrosophy serves as the basis for the Neutrosophic Logic (NL) or Smarandache's logic, which is a general framework for the unification of all existing logics. The main idea of NL is to characterize each logical statement in a 3D Neutrosophic Space, where each dimension of the space represents respectively the truth (T), the falsehood (F), and the indeterminacy (I) of the statement under consideration, where T, I, and F are standard or non-standard real subsets of $]^{-}0, 1^{+}[$. Moreover, in NL each statement is allowed to be over or under true, over or under false, and over or under indeterminate by using hyper real numbers developed in the non-standard analysis theory. The neutrosophical value $\mathcal{N}(A) = (T(A), I(A), F(A))$ in a frame of discernment (world of discourse) $\Theta$ of a statement <A> is then defined as a subset (a volume, not necessary connective; i.e. a set of disjoint volumes) of the neutrosophic space. Any statement <A>, represented by a triplet $\mathcal{N}(A)$, is called a *neutrosophic event* or *$\mathcal{N}$-event*. […]

This Smarandache's representation is close to the human reasoning. It characterizes and catches the imprecision of knowledge or linguistic inexactitude received by various observers, uncertainty due to incomplete knowledge of acquisition errors or stochasticity, and vagueness due to lack of clear contours or limits. This approach allows theoretically to consider any kinds of logical statements. For example, the fuzzy set logic or the classical modal logic (which works with statements verifying $T(A), I(A) \equiv 0, F(A) = 1-T(A)$, where T is a real number belonging to [0, 1]) are included in NL. The neutrosophic logic can easily handle also paradoxes. We emphasize the fact that in general the neutrosohic value $\mathcal{N}(A)$ of a proposition <A> can also depend on dynamical parameters which can evolve with time, space, etc.

[…] the foundations for a new theory of paradoxical and plausible reasoning [i.e. Dezert-Smarandache Theory, *ref. n.*] has been developed, that takes into account in the combination process itself the possibility for uncertain and paradoxical information. The basis for the development of this theory is to work with the hyper-power set of the frame of discernment relative to the problem under consideration rather than its classical power set since, in general, the frame of discernment cannot be fully described in terms of an exhaustive and exclusive list of disjoint elementary hypotheses. In such general case, no refinement is possible to apply directly the Dempster-Shafer Theory (DST) of evidence. In our new theory [i.e. DSmT, *ref. n.*], the rule of combination is justified from the maximum entropy principle and there is no mathematical impossibility to combine sources of evidence even if they appear at first glance in contradiction (in the Shafer's sense) since the paradox between sources is fully taken into account in our formalism. We have also shown that in general, the combination of evidence yields unavoidable paradoxes.

This theory has shown, through many illustrated examples, that conclusion drawn from it provides results which agree perfectly with the human reasoning and is useful to take a decision on complex problems where DST usually fails. […] this work has been devoted to the development of a theoretical bridge between the neutrosophic logic and this new theory, in order to solve the delicate problem of the combination of neutrosophic evidences. The neutrosophic logic serves here as the most general framework (prerequisite) for dealing with uncertain and paradoxical sources of information through this new theory.

<div align="right">Dr. Jean Dezert</div>